\documentclass[a4paper,reqno]{amsart}

\usepackage{verbatim} 
\usepackage{enumerate}
\usepackage{paralist}
\usepackage{tabto}
\usepackage{slashbox}
\usepackage{etoolbox}

\usepackage{tikz}
\usepackage{tikz-cd}

\pgfqkeys{/tikz/commutative diagrams}{%
  path/.style={/tikz/arrows=cm cap-cm cap}%
}

\usepackage{pdftricks}
\begin{psinputs}
\usepackage{pst-coil}
\end{psinputs}
\usepackage{amsthm} 
\usepackage[pdf,dvips,color,all,poly]{xy}
\usepackage{stmaryrd}
\usepackage{latexsym}
\usepackage{amssymb}
\usepackage{amsmath}
\usepackage{extarrows}
\usepackage{upgreek}
\usepackage{mathrsfs} 
\newtheorem{theorem}[equation]{Theorem}
\newtheorem{proposition}[equation]{Proposition}
\newtheorem{lemma}[equation]{Lemma}
\newtheorem{corollary}[equation]{Corollary}

\theoremstyle{definition} 
\newtheorem{definition}[equation]{Definition}
\newtheorem{remark}[equation]{Remark}
\theoremstyle{remark} 
\newtheorem{example}[equation]{Example}

\renewcommand{\and}{\quad\text{and}\quad}

\newcommand{\ZZ}{\mathbb{Z}}

\newcommand{\lra}{\longrightarrow}

\newcommand{\cc}{{\mathcal C}}
\newcommand{\cd}{{\mathcal D}}

\newcommand{\cj}{{\mathcal J}}
\newcommand{\cp}{{\mathcal P}}
\newcommand{\ck}{{\mathcal K}}
\newcommand{\cR}{{\mathcal R}}
\newcommand{\cs}{{\mathcal S}}

\newcommand{\DR}{\mathop{{\operatorname{D}}_{\!R}}}
\newcommand{\DK}{\mathop{{\operatorname{D}}_{\!K}}}

\newcommand{\cpp}{{\mathcal P}(Q_\sigma,W_\sigma,F)}

\newcommand{\G}{\Gamma_\sigma}

\newcommand{\com}[1]{}

\newcommand\spovermat[1]{\newline$\vphantom{\smash{\overbrace{\phantom{\begin{matrix}\end{matrix}}}^{#1}}}$}

\newcommand\bovermat[2]{%
  \makebox[0pt][l]{$\smash{\overbrace{\phantom{%
    \begin{matrix}#2\end{matrix}}}^{#1}}$}#2}

\renewcommand{\epsilon}{\varepsilon}

\renewcommand{\mod}{\mathrm{mod}\,}
\DeclareMathOperator{\proj}{proj}
\DeclareMathOperator{\add}{add}

\DeclareMathOperator{\im}{Im}
\DeclareMathOperator{\soc}{soc}
\DeclareMathOperator{\M}{M}

\DeclareMathOperator{\Id}{Id}

\DeclareMathOperator{\coker}{coker}
\DeclareMathOperator{\Hom}{Hom}
\DeclareMathOperator{\Ext}{Ext}
\DeclareMathOperator{\End}{End}
\DeclareMathOperator{\Aut}{Aut}

\DeclareMathOperator{\inj}{inj}
\DeclareMathOperator{\CM}{CM}
\DeclareMathOperator{\uCM}{\underline{CM}}
\DeclareMathOperator{\oCM}{\overline{CM}}
\DeclareMathOperator{\uHom}{\underline{Hom}}
\DeclareMathOperator{\uEnd}{\underline{End}}
\DeclareMathOperator{\umod}{\underline{mod}}
\DeclareMathOperator{\thick}{thick}
\DeclareMathOperator{\diag}{diag}
\DeclareMathOperator{\gldim}{gl. dim}
\DeclareMathOperator{\Tr}{Tr}

\newcommand{\tr}[1]{\,{\vphantom{#1}}^{\text{t}}\!{#1}}

\newcommand{\shift}[2]{\ifstrequal{#2}{1}{\Sigma #1}{\Sigma^{#2}#1}}

\newcommand{\vecv}[2]{
    \begin{bmatrix}
      #1 \\
      #2
    \end{bmatrix}
    }
\newcommand{\vech}[2]{
    \begin{bmatrix}
      #1 & #2
    \end{bmatrix}
    }
\newcommand{\svecv}[2]{
\scriptstyle \setlength{\arraycolsep}{2pt}
    \begin{bmatrix}
      #1 \\
      #2
    \end{bmatrix}
    }
\newcommand{\svech}[2]{
\scriptstyle \setlength{\arraycolsep}{2pt}
    \begin{bmatrix}
      #1 &  #2
    \end{bmatrix}
    }

\newcommand{\vecvt}[3]{
    \begin{bmatrix}
      #1 \\
      #2 \\
      #3
    \end{bmatrix}
    }
\newcommand{\vecht}[3]{
    \begin{bmatrix}
      #1 & #2 & #3
    \end{bmatrix}
    }
\newcommand{\mat}[4]{
    \begin{bmatrix}
      #1 & #2 \\
      #3 & #4
    \end{bmatrix}
    }

\newcommand{\mattd}[6]{
    \begin{bmatrix}
      #1 & #2 \\
      #3 & #4 \\
      #5 & #6
    \end{bmatrix}
    }

\newcommand{\matdt}[6]{
    \begin{bmatrix}
      #1 & #2 & #3 \\
      #4 & #5 & #6
    \end{bmatrix}
    }

\newcommand{\cci}[2]{\mathopen\llbracket #1,#2 \mathclose\rrbracket}
\newcommand{\coi}[2]{\mathopen\llbracket #1,#2 \mathclose\llbracket}
\newcommand{\oci}[2]{\mathopen\rrbracket #1,#2 \mathclose\rrbracket}
\newcommand{\ooi}[2]{\mathopen\rrbracket #1,#2 \mathclose\llbracket}
\newcommand{\delt}[1]{\delta_{#1}}

\newcommand{\notvdash}{\hspace{.16cm}\raisebox{2pt}{$\scriptstyle{/}$}\hspace{-.3cm}\vdash}
\newtheorem{claim}{Claim}

\newenvironment{sbmatrix}{\scriptstyle\setlength{\arraycolsep}{2pt}\begin{bmatrix}}{\end{bmatrix}}

\allowdisplaybreaks

\numberwithin{equation}{section}

\AtBeginEnvironment{figure}{\stepcounter{equation}}

\title[Ice quivers with potential arising from once-punctured polygons]{Ice quivers with potential arising from once-punctured polygons and Cohen-Macaulay modules}
\author[L. Demonet]{Laurent Demonet}
\address{Graduate School of Mathematics, Nagoya University, Furo-cho, Chikusa-ku, 464-8602 Nagoya, Japan}
\email{Laurent.Demonet@normalesup.org}
\author[X. Luo]{Xueyu Luo}
\address{Graduate School of Mathematics, Nagoya University, Furo-cho, Chikusa-ku, 464-8602 Nagoya, Japan}
\email{luoxueyu2012@gmail.com}
\begin{document}
\begin{abstract}
Given a tagged triangulation of a once-punctured polygon $P^*$ with $n$ vertices, we associate an ice quiver with potential such that the frozen part of the associated frozen Jacobian algebra has the structure of a Gorenstein $K[X]$-order $\Lambda$.
Then we show that the stable category of the category of Cohen-Macaulay $\Lambda$-modules is equivalent to the cluster category $\cc$ of type $D_n$.
 It gives a natural interpretation of the usual indexation of cluster tilting objects of $\cc$ by tagged triangulations of $P^*$. Moreover, it extends naturally the triangulated categorification by $\cc$ of the cluster algebra of type $D_n$ to an exact categorification by adding coefficients corresponding to the sides of $P$. Finally, we lift the previous equivalence of categories to an equivalence between the stable category of graded Cohen-Macaulay $\Lambda$-modules and the bounded derived category of modules over a path algebra of type $D_n$. 
\end{abstract}

\maketitle

\section{Introduction}
\label{intro}

In a previous paper \cite{DeLu}, we constructed ice quivers with potential arising from triangulations of polygons and we proved that the frozen parts of their frozen Jacobian algebras are orders. We proved that the categories of Cohen-Macaulay modules over these orders are stably equivalent to cluster categories of type $A$. The aim of this paper is to extend these results to tagged triangulations of once-punctured polygons to recover cluster categories of type $D$. We refer to \cite{DeLu} for a detailed introduction of the context and we will focus here on the tools we need specifically for this new case.



For every bordered surface with marked points, Fomin, Shapiro and Thurston introduced the concept of tagged triangulations and their mutations \cite{FST08}. Then, they associated to each of these triangulations a quiver $Q(\sigma)$ and showed that the combinatorics of triangulations of the surface correspond to that of the cluster algebra defined by $Q(\sigma)$.
Later in \cite{LF09}, Labardini-Fragoso associated a potential $W(\sigma)$ on $Q(\sigma)$. 
He proved that flips of triangulations are compatible with mutations of quivers with potential.
This was generalized to the case of tagged triangulations by Labardini-Fragoso and Cerulli Irelli in \cite{LF11, LF12}.

We refer to \cite{FMO,RT,POS,CM} for a general background on Cohen-Macaulay modules (or lattices) over orders. Recently, strong connections between Cohen-Macaulay representation theory and tilting theory, especially cluster categories, have been established \cite{AIRCC,Araya,Thanvanden,IyTa13,Matrixrepresentation,Matrixweight,KR}.
This paper enlarges some of these connections by dealing with frozen Jacobian algebras associated with tagged triangulations of once-punctured polygons from the viewpoint of Cohen-Macaulay representation theory.

Through this paper, let $K$ denote a field and $R=K[X]$.  We extend the construction of \cite{FST08}, and associate an ice quiver with potential $(Q_\sigma, W_\sigma, F)$ to each tagged triangulation $\sigma$ of a once-punctured polygon $P^*$ with $n$ vertices by adding a set $F$ of $n$ frozen vertices corresponding to the edges of the polygon and certain arrows (see Definition \ref{def:quiver}). We study the associated frozen Jacobian algebra 
\[  \Gamma_\sigma:=\cpp \]
(see Definition \ref{frozen JA}). Our main results are the following:
\begin{theorem}
[Theorem \ref{thm:frozen Jacobian algebra is order} and Theorem \ref{cor:lambdamod}]
Let $e_F$ be the sum of the idempotents of $\Gamma_\sigma$ at frozen vertices. Then 
\begin{enumerate}
\item the frozen Jacobian algebra $\Gamma_\sigma$ has the structure of an $R$-order (see Definition \ref{def:order} and Remark \ref{rem:order});
\item the frozen part $e_F \Gamma_\sigma e_F$ is isomorphic to the Gorenstein $R$-order
 \begin{equation}
    \label{the order}
 \Lambda:=  \begin{bmatrix}
    R' & R' & R' &\cdots&R' & \hspace*{-.4cm}X^{-1}(X,Y)\\
    (X,Y) & R'& R'&\cdots& R'&\hspace*{-.4cm}R'\\
    (X) & (X,Y)& R' &\cdots& R' &\hspace*{-.4cm} R'\\
    \vdots&\vdots&\vdots&\ddots&\vdots&\hspace*{-.4cm}\vdots\\
   (X) &(X)& (X)&\cdots& R' & \hspace*{-.4cm}R'\\
   (X) &(X)& (X)&\cdots  & (X,Y) &\hspace*{-.4cm}R'
   \end{bmatrix}_{n\times n,}
\end{equation}
 where $R' = K[X,Y]/(Y(X-Y))$ and each entry of the matrix is a $R'$-submodule of $R'[X^{-1}]$.
\end{enumerate}
\end{theorem} 

\begin{remark}
 In view of the isomorphism of $R$-algebras
 \begin{align*}
  R' & \cong R - R := \{(P,Q) \in R^2 \,|\, P-Q \in (X)\} \\
  Y & \mapsto (0,X),
 \end{align*}
 we have an isomorphism
 $$\Lambda \cong  \begin{bmatrix}
    R-R & R-R & R-R &\cdots&R-R & R\times R\\
    (X) \times (X) & R-R& R-R&\cdots& R-R&R-R\\
    (X)-(X) & (X) \times (X)& R-R &\cdots& R-R & R-R\\
    \vdots&\vdots&\vdots&\ddots&\vdots&\vdots\\
   (X)-(X) &(X)-(X)& (X)-(X)&\cdots& R-R & R-R\\
   (X)-(X) &(X)-(X)& (X)-(X)&\cdots  & (X) \times (X) &R-R
   \end{bmatrix}_{n\times n}$$
 ($(X)-(X)$ is the ideal of $R-R$ generated by $(X,X)$). 

 This order is part of a wide class of Gorenstein orders, called almost Bass orders, introduced and studied by Drozd-Kirichenco-Roiter and 
Hijikata-Nishida \cite{HiNi92, HiNi97} (see also \cite{Iya-order}). More precisely, $\Lambda$ is an \emph{almost Bass order of type (III)}.
\end{remark}

\begin{theorem}[Theorems \ref{cor:lambdamod}, \ref{thm:classif}, \ref{thm:cluster tilting} and \ref{th:cyeq}]
\label{introthm2}
 The category $\CM \Lambda$ has the following properties:
 \begin{enumerate}
  \item For any tagged triangulation $\sigma$ of $P^*$, we can map each tagged arc $a$ of $\sigma$ to the indecomposable Cohen-Macaulay $\Lambda$-module $e_F \Gamma_\sigma e_a$, where $e_a$ is the idempotent of $\Gamma_\sigma$ at $a$. This module does only depend on $a$ (not on $\sigma$) and this map induces one-to-one correspondences 
 \begin{align*}
  \left\{\text{sides and tagged arcs of $P^*$} \right\} & \longleftrightarrow \left\{\text{indecomposable objects of $\CM \Lambda$}\right\}/\cong \\
  \left\{\text{sides of $P$} \right\} & \longleftrightarrow \left\{\text{indecomposable projectives of $\CM \Lambda$}\right\}/\cong\\
  \left\{\text{tagged triangulations of $P^*$} \right\} & \longleftrightarrow \left\{\text{basic cluster tilting objects of $\CM \Lambda$}\right\}/\cong.
 \end{align*} 
  \item For the cluster tilting object $T_\sigma := e_F \Gamma_\sigma$ corresponding to a tagged triangulation $\sigma$, $$\End_{\CM \Lambda}(T_\sigma) \cong \Gamma_\sigma^{\operatorname{op}}.$$
  \item \label{intro2cy} The category $\uCM \Lambda$ is $2$-Calabi-Yau.
  \item If $K$ is a perfect field, there is a triangle-equivalence $\cc(K Q) \cong \uCM \Lambda$, where $Q$ is a quiver of type $D_n$ and $\cc(K Q)$ is the corresponding cluster category. 
 \end{enumerate}
\end{theorem}

\begin{remark}
 To prove Theorem \ref{introthm2} \eqref{intro2cy}, we establish that
 $$\CM \Lambda \cong \CM^{\ZZ/n \ZZ} \left(\frac{K[x,y]}{(x^{n-1} y - y^2)}\right),$$
 where $x$ has degree $1$ and $y$ has degree $-1$ (modulo $n$).
\end{remark}

Usually, the cluster category $\cc(K Q)$ is constructed as an orbit category of the bounded derived category $\cd^{\operatorname{b}}(K Q)$. We can reinterpret this result in this context by studying the category of graded Cohen-Macaulay $\Lambda$-modules $\CM^\ZZ \Lambda$:

\begin{theorem}[Theorem \ref{thm:cluster cat as stable cat for D type}]
  Using the same notation as before: 
  \begin{enumerate}
   \item The Cohen-Macaulay $\Lambda$-module $e_F \Gamma_\sigma$ can be lifted to a tilting object in $\uCM^\ZZ(\Lambda)$.
   \item There exists a triangle-equivalence $\cd^{\operatorname{b}}(K Q) \cong \uCM^{\ZZ}(\Lambda)$.
  \end{enumerate}
\end{theorem}

In Section \ref{Ice QP}, we introduce ice quivers with potential $(Q_\sigma, W_\sigma, F)$ associated with tagged triangulations $\sigma$ of a once-punctured polygon $P^*$. We also introduce combinatorial and algebraic elementary tools in subsection \ref{sec:notation}. Finally, we prove in this section that the frozen Jacobian algebra $\Gamma_\sigma$ associated with $(Q_\sigma, W_\sigma, F)$ is an $R$-order, as well as $\Lambda \cong e_F \Gamma_\sigma e_F$ which is independent of $\sigma$. In section \ref{sec:CM}, we classify Cohen-Macaulay modules over $\Lambda$, we compute homological properties of $\CM \Lambda$ and we establish the correspondence between tagged triangulations of $P^*$ and basic cluster tilting objects of $\CM \Lambda$. Thus, after proving that $\CM \Lambda$ is Frobenius stably $2$-Calabi-Yau, we conclude that $\CM \Lambda$ is stably triangle equivalent to a cluster category of type $D$. In section \ref{s:graded}, we deal with results about $\CM^\ZZ \Lambda$.

Notice that the naive generalization of these results to other surfaces do not hold in general as shown in Subsection \ref{cex} for a digon with two punctures.

\medskip\noindent{\bf Acknowledgements }
We would like to show our gratitude to the second author's supervisor, Osamu Iyama, for his guidance and valuable discussions. We also thank him for introducing us this problem and showing us enlightening examples.

\section{Ice quivers with potential associated with triangulations}
\label{Ice QP}
In this section, we introduce ice quivers with potential associated with tagged triangulations of a once-punctured polygon and their frozen Jacobian algebras. We show that in any case, the frozen Jacobian algebra has the structure of an $R$-order, and its frozen part is isomorphic to a given $R$-order $\Lambda$ defined in \eqref{the order}.

\subsection{Frozen Jacobian algebras}
\label{ss:Frozen Jacobian algebras}
We refer to \cite{DWZ1} for background about quivers with potential.
Let $Q$ be a finite connected quiver without loops, with set of vertices $Q_0=\{1,\ldots,n\}$ and set of arrows $Q_1$. 
As usual, if $\alpha \in Q_1$, we denote by $s(\alpha)$ its starting vertex and by $e(\alpha)$ its ending vertex.
We denote by $KQ_i$ the $K$-vector space with basis $Q_i$ consisting of paths of length $i$ in $Q$, and by $KQ_{i,{\rm cyc}}$ the subspace of $KQ_i$ spanned by all cycles in $KQ_i$. 
Consider the path algebra $KQ=\bigoplus_{i\ge0}KQ_i$. 
An element $W\in \bigoplus_{i\ge1}KQ_{i,{\rm cyc}}$ is called a {\em potential}.
Two potentials $W$ and $W'$ are called \emph{cyclically equivalent} if $W-W'$ belongs to $[KQ,KQ]$, the vector space spanned by commutators.
A \emph{quiver with potential} is a pair $(Q,W)$ consisting of a quiver $Q$ without loops and a potential $W$ which does not have two cyclically equivalent terms.

For each arrow $\alpha\in Q_1$, the cyclic derivative $\partial_\alpha$ is the linear map $\bigoplus_{i\ge1}KQ_{i,{\rm cyc}}\to KQ$ defined on cycles by
\[
\partial_\alpha(\alpha_1\cdots \alpha_d)=\sum_{\alpha_i=\alpha}\alpha_{i+1}\cdots \alpha_d \alpha_1\cdots \alpha_{i-1}.\]

\begin{definition}[\cite{BIRS}]
\label{frozen JA}
An \emph{ice quiver with potential} is a triple $(Q,W,F)$, where $(Q,W)$ is a quiver with potential and $F$ is a subset of $Q_0$.
Vertices in $F$ are called \emph{frozen vertices}.

   The \emph{frozen Jacobian algebra} is defined by \[\cp(Q,W,F)=KQ/\cj(W,F),\] where $\cj(W,F)$ is the ideal \[\cj(W,F)=\langle \partial_\alpha W \mid \alpha \in Q_1, \ s(\alpha)\notin F\ \mbox{ or }\ e(\alpha)\notin F \rangle\] of $KQ$. 
\end{definition}

\begin{figure}[t]
 \[{ \begin{tikzcd}[ampersand replacement=\&]
         \&             \&               \& 1\drar{\beta_1}                 \&                \&                 \&\\
         \&             \& 2\drar{\beta_2}\urar{\alpha_1}     \&                      \& 3\arrow{ll}[swap]{\gamma_1} \drar{\beta_3} \&    \&\\
       \& 4\urar{\alpha_2}   \&        \& 5\arrow{ll}[swap]{\gamma_2}\urar{\alpha_3} \&   \& 6\arrow{ll}[swap]{\gamma_3} \&\\ 
  \end{tikzcd}}
  \]
\caption{Example of iced quiver with potential}
\label{ausl}
\end{figure}
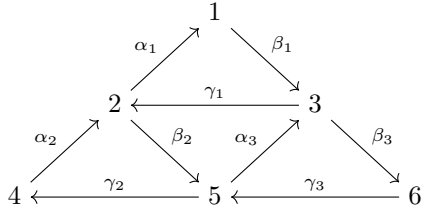

\begin{example}
Consider the quiver $Q$ of Figure \ref{ausl}
 with potential $W=\alpha_1 \beta_1 \gamma_1 + \alpha_2 \beta_2 \gamma_2 + \alpha_3 \beta_3 \gamma_3- \gamma_1 \beta_2 \alpha_3$
  and set of frozen vertices $F=\{4,5,6\}.$
  Then the Jacobian ideal is 
  \begin{align*}
  \cj(W,F) =\langle \beta_1 \gamma_1, \gamma_1 \alpha_1, \alpha_1 \beta_1 - \beta_2 \alpha_3, \beta_2 \gamma_2, \gamma_2 \alpha_2 - \alpha_3 \gamma_1, \beta_3 \gamma_3 - \gamma_1 \beta_2, \gamma_3 \alpha_3  \rangle.
  \end{align*}

Note that this ice quiver with potential appeared from preprojective algebras \cite{BIRS,GLS11}.
\end{example}

\subsection{Ice quivers with potential arising from triangulations}
\label{ss:Ice QP}
We recall the definition of triangulations of a polygon with one puncture and introduce our definition of ice quivers with potential arising from tagged triangulations of a polygon with one puncture.

\begin{definition}
Let $P$ be a regular polygon with $n$ vertices and $n$ sides. Fix a marked point in the center of $P$. Then the marked point is called a \emph{puncture} and the combination of $P$ and the puncture is called a \emph{(once-)punctured polygon $P^*$}.
We call \emph{interior of $P^*$} for the interior of the polygon $P$ excluding the puncture. We denote by $M$ the set of the $n$ vertices of the polygon and the puncture.
\end{definition}

\begin{definition}[Tagged arcs]
\cite[Definition 7.1]{FST08}
 A \emph{tagged arc} in the punctured polygon $P^*$ is a curve $a$ in $P$ such that 
 \begin{enumerate}
   \item the endpoints of $a$ are distinct in $M$;
   \item $a$ does not intersect itself;
   \item except for the endpoints, $a$ is disjoint from $M$ and from the sides of $P$;
   \item $a$ does not cut out an unpunctured digon. (In other words, $a$ is not contractible onto the sides of $P$.)
 \end{enumerate}
 Each arc $a$ is considered up to isotopy inside the class of such curves.

 Moreover, each arc incident to the puncture has to be tagged either \emph{plain} either \emph{notched}.
 
 In the figures, the plain tags are omitted while the notched tags are  represented by the symbol $\bowtie$.
\end{definition}

\begin{definition}[Compatibility of tagged arcs, {\cite[Definition 7.4]{FST08}}]
Two tagged arcs $a$ and $b$ are \emph{compatible} if and only if the following conditions are satisfied:
 \begin{enumerate}
   \item there are curves in their respective isotopy classes whose relative interiors do not intersect;
   \item if $a$ and $b$ are incident to the puncture and not isotopic, they are either both plain, either both notched.
 \end{enumerate}
\end{definition}

\begin{definition}[\cite{FST08}]
A \emph{tagged triangulation} of the punctured polygon $P^*$ is the union of the set of the sides of $P$ and any maximal collection of pairwise compatible tagged arcs of $P^*$.
\end{definition}

\begin{remark}
 The set of all tagged arcs in a punctured polygon is finite. Moreover, any tagged triangulation can be realized up to isotopy as a collection of tagged non-intersecting arcs.
\end{remark}

Let us now define ice quivers with potential arising from punctured polygons:
\begin{definition}
\label{def:quiver}
Let $P^*$ be a punctured polygon with $n$ sides and $\sigma$ be a tagged triangulation of $P^*$. For convenience, the $n$ sides of $P$ and all the tagged arcs of $\sigma$ are called the edges of $\sigma$. A \emph{true triangle} of $\sigma$ is a triangle consisting of edges of $\sigma$ such that the puncture is not in its interior. 

We assign to $\sigma$ two ice quivers with potential $(Q_\sigma, W_\sigma, F)$ and $(Q'_\sigma, W'_\sigma, F)$ as follows.

The quiver $Q'_\sigma$ is a quiver whose vertices are indexed by the edges of $\sigma$. Whenever two edges $a$ and $b$ are sides of a common true triangle of $\sigma$, then $Q'_\sigma$ contains an internal arrow $a \to b$ in the true triangle if $a$ is a predecessor of $b$ with respect to anticlockwise orientation centred at the joint vertex.
For every vertex of the polygon $P$, there is an external arrow $a \to b$ where $a$ and $b$ are its two incident sides of $P$, $a$ being a predecessor of $b$ with respect to anticlockwise orientation centred at the joint vertex.
Moreover, if the puncture is adjacent to exactly one notched arc and one plain arc of $\sigma$, we have the configuration shown in Figure \ref{digon}.

\begin{figure}
\begin{minipage}{4cm}
{
\begin{pdfpic}
\begin{pspicture}(0,-2.4)(3.61,2.28)
\psdots[dotsize=0.12](1.6,0.6)
\psbezier[linewidth=0.02](1.6,2.2)(2.8,1.8)(3.6,-1.0)(1.6,-2.2)
\psbezier[linewidth=0.02](1.6,2.2)(0.0,1.8)(0.0,-1.4)(1.6,-2.2)
\psline[linewidth=0.02cm](1.6,0.6)(1.6,-2.2)
\psbezier[linewidth=0.02](1.6,0.6)(2.4,0.2)(2.0,-1.4)(1.6,-2.2)
\usefont{T1}{ptm}{m}{n}
\rput(2.051455,-0.095){$\bowtie$}
\usefont{T1}{ptm}{m}{n}
\rput(0.24145508,0.305){$a$}
\usefont{T1}{ptm}{m}{n}
\rput(3.051455,0.305){$b$}
\usefont{T1}{ptm}{m}{n}
\rput(1.411455,-0.895){$i$}
\usefont{T1}{ptm}{m}{n}
\rput(2.211455,-0.895){$j$}
\usefont{T1}{ptm}{m}{n}
\rput(1.6,-2.4){$P_k$}
\end{pspicture} 
\end{pdfpic}
}
\end{minipage}
\begin{minipage}{6.5cm}\begin{tikzcd}[ampersand replacement=\&]
\&\& j \arrow{dll}[swap]{\beta}\&\& \\
a\arrow[bend left=90]{rrrr}{\eta} \&\&\& \& b \arrow{ull}[swap]{\alpha} \arrow{dll}[swap]{\gamma} \\
\&\& i\arrow{ull}[swap]{\delta}\&\& \\
\&\&\&\&\&
\end{tikzcd}\end{minipage}

We depict the part of $Q_\sigma$ quiver induced by a once-punctured digon. Notice that there are two other arrows linking $a$ and $b$ if $P$ is the digon itself.
\caption{Once-punctured digon}
\label{digon}
\end{figure}

Then, the quiver $Q_\sigma$ is obtained from $Q'_\sigma$ by removing external arrows winding around vertices of $P$ with no incident tagged arc in $\sigma$.

We say that a cycle of $Q_\sigma$ (resp. $Q'_\sigma$) is \emph{planar} if it does not contain any arrow of $Q_\sigma$ (resp. $Q'_\sigma$) in its interior  and each arrow appears at most once. Notice that for the definition of planar, the quivers are not abstract but embedded in the plane (each internal arrow being drawn inside the triangle it is constructed from, and each external arrow winding around the corresponding vertex outside the polygon). We have the following possible different kinds of planar cycles in $Q_\sigma$ and $Q'_\sigma$:
\begin{enumerate}
 \item \emph{clockwise triangles} which come from true triangles in $\sigma$;
 \item an \emph{anticlockwise punctured cycle} which consists of the arrows connecting arcs incident to the puncture.
 \item \emph{anticlockwise external cycles} which contain exactly one external arrow and each of which is centered at a vertex of $P$.
\end{enumerate}

We define $F$ as the subset of $(Q_\sigma)_0$ indexed by the $n$ sides of the polygon $P$.
The potential $W_\sigma$ (resp. $W'_\sigma$) is defined as 
\begin{align*}
& \sum{\text{clockwise triangles}} - \sum{\text{anticlockwise external cycles}} \\
 &  - \text{the anticlockwise punctured cycle}
\end{align*}
in $Q_\sigma$ (resp. $Q'_\sigma$).

When there is a once-punctured digon in the triangulation $\sigma$ as shown in Figure \ref{digon}, we have to adapt slightly the previous definition. The anticlockwise external cycle centred at $P_k$ which is taken in account is the one containing $\gamma \delta$. On the other hand, both $\eta \alpha \beta$ and $\eta \gamma \delta$ appear as clockwise triangles in $W_\sigma$ and $W'_\sigma$. In this case, there is no anticlockwise punctured cycle. An explicit case involving a once-punctured digon is  described in Proof of Lemma \ref{lemdig}.
\end{definition}

\begin{example}
 Let us consider the triangulation $\sigma$ of Figure \ref{exmp}. We drew the corresponding quivers ($\gamma$ is in $Q'_\sigma$ but not in $Q_\sigma$). We have
 $$W_\sigma = fgh + abc + ade - \alpha ag - \beta fbc \quad \text{and} \quad W'_\sigma = W_\sigma - \gamma h.$$
 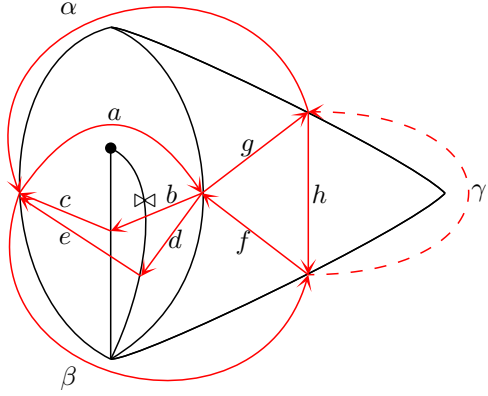
\begin{figure}
{
\begin{pdfpic}
\begin{pspicture}(0,-2.8)(6.5,2.8)
\psdots[dotsize=0.15](1.6,0.6)
\psbezier[linewidth=0.02](1.6,2.2)(2.8,1.8)(3.6,-1.0)(1.6,-2.2)
\psbezier[linewidth=0.02](1.6,2.2)(0.0,1.8)(0.0,-1.4)(1.6,-2.2)
\psline[linewidth=0.02cm](1.6,0.6)(1.6,-2.2)
\psbezier[linewidth=0.02](1.6,0.6)(2.4,0.2)(2.0,-1.4)(1.6,-2.2)
\psbezier[linewidth=0.02](1.6,-2.2)(2,-2.2)(5.5,-0.5)(6.0,0)
\psbezier[linewidth=0.02](1.6,2.2)(2,2.2)(5.5,0.5)(6.0,0)
\psline[linewidth=0.02cm,linecolor=red,arrowsize=0.05cm 8.0,arrowlength=1.0,arrowinset=0.6]{->}(2.8,0)(2.0,-1.1)
\psline[linewidth=0.02cm,linecolor=red,arrowsize=0.05cm 8.0,arrowlength=1.0,arrowinset=0.6]{->}(2.0,-1.1)(0.4,-0.05)
\psline[linewidth=0.02cm,linecolor=red,arrowsize=0.05cm 8.0,arrowlength=1.0,arrowinset=0.6]{->}(2.8,0)(1.6,-0.5)
\psline[linewidth=0.02cm,linecolor=red,arrowsize=0.05cm 8.0,arrowlength=1.0,arrowinset=0.6]{->}(1.6,-0.5)(0.4,0)
\psbezier[linewidth=0.02cm,linecolor=red,arrowsize=0.05cm 8.0,arrowlength=1.0,arrowinset=0.6]{->}(0.4,0)(1.2,1.2)(2,1.2)(2.8,0)
\psbezier[linewidth=0.02](1.6,-2.2)(2,-2.2)(5.5,-0.5)(6.0,0)
\psbezier[linewidth=0.02](1.6,2.2)(2,2.2)(5.5,0.5)(6.0,0)
\psline[linewidth=0.02cm,linecolor=red,arrowsize=0.05cm 8.0,arrowlength=1.0,arrowinset=0.6]{->}(2.8,0)(4.2,1.08)
\psline[linewidth=0.02cm,linecolor=red,arrowsize=0.05cm 8.0,arrowlength=1.0,arrowinset=0.6]{->}(4.2,1.08)(4.2,-1.08)
\psline[linewidth=0.02cm,linecolor=red,arrowsize=0.05cm 8.0,arrowlength=1.0,arrowinset=0.6]{->}(4.2,-1.08)(2.8,0)
\psbezier[linewidth=0.02cm,linecolor=red,arrowsize=0.05cm 8.0,arrowlength=1.0,arrowinset=0.6]{->}(4.2,1.08)(3.5,3.5)(-0.5,2.5)(0.4,0)
\psbezier[linewidth=0.02cm,linecolor=red,arrowsize=0.05cm 8.0,arrowlength=1.0,arrowinset=0.6]{->}(0.4,0)(-0.5,-2.5)(3.5,-3.5)(4.2,-1.08)
\psbezier[linewidth=0.02cm,linecolor=red,arrowsize=0.05cm 8.0,arrowlength=1.0,arrowinset=0.6,linestyle=dashed,dash=0.16cm 0.16cm]{->}(4.2,-1.08)(7,-1.08)(7,1.08)(4.2,1.08)
\usefont{T1}{ptm}{m}{n}
\rput(1.05,2.45){$\alpha$}
\rput(1.05,-2.45){$\beta$}
\rput(6.45,0){$\gamma$}
\rput(1.65,1.05){$a$}
\rput(2.4,0){$b$}
\rput(1,-.1){$c$}
\rput(2.45,-.6){$d$}
\rput(1,-.6){$e$}
\rput(3.35,-.65){$f$}
\rput(3.4,.6){$g$}
\rput(4.35,0){$h$}

\usefont{T1}{ptm}{m}{n}
\rput(2.051455,-0.095){$\bowtie$}
\end{pspicture} 
\end{pdfpic}
}
  \caption{Quivers associated with a tagged triangulation}
  \label{exmp}
 \end{figure}
\end{example}

From now on, for any tagged triangulation $\sigma$ of the punctured polygon $P^*$, denote $\cpp$ by $\G$.

\begin{remark}
 \begin{enumerate}
  \item We can also realize $\G$ as
 $$\G \cong \frac{KQ'_\sigma}{\mathcal{J}'(W'_\sigma)} \quad
 \text{where} \quad \mathcal{J}'(W'_\sigma) := \langle \partial_\alpha W'_\sigma \mid \alpha \in Q_1, \ \alpha \text{ is not external} \rangle.$$

 We will use both definitions freely depending on convenience.
  \item  Notice that $(Q'_\sigma, W'_\sigma)$ and $(Q_\sigma, W_\sigma)$ are not necessarily reduced, in the sense that oriented $2$-cycles can appear in the potential, because some vertices of the polygon have no incident tagged arcs in $\sigma$, or because the puncture has exactly two non-isotopic incident tagged arcs. Thus, it is possible that non-admissible relations appear.
  \item All arrows of $Q_\sigma$ (resp. $Q'_\sigma$) either appear once in $W_\sigma$ (resp. $W'_\sigma$), either twice with opposite signs. Thus all relations derived from the potential are either commutativity relations (of the form $w = w'$ for two paths $w$ and $w'$ of length at least $1$) or $0$ relations (of the form $w = 0$ for a path $w$ of length at least $1$). 
 \end{enumerate}
\end{remark}

\subsection{Notation and preliminaries} \label{sec:notation}

The vertices of the polygon $P$ are called $P_1$, $P_2$, \dots, $P_n$ counter-clockwisely. When we do computations on the indices of vertices of $P$, we compute modulo $n$. If $r,s \in \cci{1}{n} := \{1,2, \dots, n\}$, we denote
 $$d(r, s) := {\begin{cases}
s-r, & \text{if } s \ge r, \\
s-r+n, & \text{if } s<r.
\end{cases}}$$

We also denote $\cci{r}{s} = \{r,r+1, \dots, s\}$, $\coi{r}{s} = \cci{r}{s} \smallsetminus \{s\}$, $\oci{r}{s} = \cci{r}{s} \smallsetminus \{r\}$ and $\ooi{r}{s} = \cci{1}{n} \smallsetminus \cci{s}{r}$ (notice that $\ooi{r}{r} = \cci{1}{n} \smallsetminus \{r\}$). If $A$ is a condition, we define $\delt{A}$ to be $1$ if $A$ is satisfied and $0$ if $A$ is not satisfied.

For $r, s, t \in \cci{1}{n}$, we will use freely the identities $\oci{r}{s} = \coi{r+1}{s+1}$, $\oci{r}{s} = \cci{1}{n} \smallsetminus \oci{s}{r}$, $\coi{r}{s} = \cci{1}{n} \smallsetminus \coi{s}{r}$ and $\delt{r \in \ooi{s}{t}} = \delt{s \in \coi{t}{r}} = \delt{t \in \oci{r}{s}}$.

If $r, s \in \cci{1}{n}$, we denote by $(P_r, P_s)$ the arc going from $P_r$ to $P_s$ turning counter-clockwisely around the puncture (thus, $(P_r, P_{r+1})$ is a side of the polygon and $(P_{r+1}, P_r)$ is not except if $n = 2$). We denote by $(P_r, *)$ the plain tagged arc from $P_r$ to the puncture and by $(P_r, \bowtie)$ the notched tagged arc from $P_r$ to the puncture. 

From now on, we always denote $a = (P_{a_1}, P_{a_2})$ if $a$ is not incident to the puncture and $a = (P_{a_1}, *)$ or $a = (P_{a_1}, \bowtie)$ if $a$ is incident to the puncture (in the latter case, we fix by convention $a_2 = a_1$). We need to fix some geometrical definition. To make it precise, we suppose that $P$ is a regular polygon inscribed in the unit circle and that the puncture is at the origin of the plane. Then, $\vec a$ is the vector from $P_{a_1}$ to $P_{a_2}$ if $a_1 \neq a_2$ and it is the unit vector tangent at $P_{a_1}$ to the unit circle in the clockwise direction if $a_1 = a_2$.

If $a$ and $b$ are tagged arcs of $P^*$, we define
$$\ell^\theta_{a,b} := d(a_1,b_1) + d(a_2, b_2) + n \left|\delt{a_1 \in \ooi{b_1}{a_2}} - \delt{b_2 \in \ooi{b_1}{a_2}}\right|.$$

Lemmas \ref{lem:easyangle} and \ref{lem:subadd} are elementary observations about $\ell^\theta$. Proofs are computational and given for the sake of completeness. We suggest the reader to skip them first.

\begin{lemma} \label{lem:easyangle}
 If $a$ and $b$ are two sides or tagged arcs of $P^*$, the angle from $\vec a$ to $\vec b$ is
 $$\frac{\pi}{n} \ell^\theta_{a,b},$$
 up to a multiple of $2 \pi$.
\end{lemma}

\begin{proof}
 First of all, in complex coordinates, if $a_1 \neq a_2$,
 \begin{align*}\vec a &= \exp\left(2 \pi i\frac{a_2}{n}\right) - \exp\left(2 \pi i\frac{a_1}{n}\right) \\ &= \exp\left(\pi i\frac{a_1 + a_2}{n}\right) \left(\exp\left(\pi i\frac{a_2-a_1}{n}\right) - \exp\left(\pi i\frac{a_1-a_2}{n}\right)\right)  \\ &= \exp\left(\pi i\frac{a_1 + a_2}{n}\right) 2i \sin\left( \pi \frac{a_2-a_1}{n} \right) \end{align*}
 so the argument of $\vec a$ is
 $$\frac{\pi}{n}\left(a_1+a_2+\frac{n}{2}+n \delt{a_1 \geq a_2}\right)$$
 (note that this formula works also if $a_1 = a_2$). So the angle from $\vec a$ to $\vec b$ is
 $$\frac{\pi}{n}\left(b_1 - a_1+ b_2 - a_2 + n \left(\delt{b_1 \geq b_2} - \delt{a_1 \geq a_2} \right)\right).$$
 As $\ell^\theta_{a,b}$ is clearly invariant by rotation of the polygon, as well as the angle from $\vec a$ to $\vec b$, we can suppose that $a_2 = 1$ and the angle from $a$ to $b$ becomes
 \begin{align*}&\frac{\pi}{n}\left(d(a_1, b_1) - n \delt{a_1 > b_1} + d(1, b_2)  + n \left(\delt{b_1 \geq b_2} - 1 \right)\right)
  \\ = & \frac{\pi}{n}\left(d(a_1, b_1) + d(1, b_2) - n \delt{a_1 > b_1} - n \delt{b_1 < b_2} \right)
  \\ = & \frac{\pi}{n}\left(d(a_1, b_1) + d(1, b_2) - n \delt{a_1 \in \ooi{b_1}{1}} - n \delt{b_2 \in \ooi{b_1}{1}} \right)
 \end{align*}
 which is clearly congruent to $\pi \ell^\theta_{a,b} / n$ modulo $2 \pi$.
\end{proof}

Another important point is that $\ell^\theta$ is subadditive:

\begin{lemma} \label{lem:subadd} \label{caldif}
 If $a$, $b$ and $c$ are three sides or tagged arcs of $P^*$, then $\ell^\theta_{a,b} + \ell^\theta_{b,c} \geq \ell^\theta_{a,c}$. More precisely, 
 \begin{itemize}
  \item if $a$ is a side of $P$,
   $$\ell^\theta_{a,b} + \ell^\theta_{b,c} = \ell^\theta_{a,c} + 2 n \left( \delt{c_2 \in \ooi{c_1}{b_2}}\delt{b_1 \in \cci{b_2}{c_1}} + \delt{a_1 \in \oci{b_1}{c_1}}  \right);$$
  \item if $b$ is a side of $P$,
   $$\ell^\theta_{a,b} + \ell^\theta_{b,c} = \ell^\theta_{a,c} + 2 n \left( \delt{a_1 \in \oci{a_2-1}{c_1}} \delt{c_2 \in \oci{a_2-1}{c_1}} + \delt{b_1 \in \ooi{c_1}{a_2-1}}  \right);$$
  \item if $c$ is a side of $P$,
   $$\ell^\theta_{a,b} + \ell^\theta_{b,c} = \ell^\theta_{a,c} + 2 n \left(\delt{a_1 \in \ooi{b_1}{a_2}}  \delt{b_2 \in \cci{a_2}{b_1}} + \delt{c_2 \in \coi{a_2}{b_2}} \right).$$
 \end{itemize}

\end{lemma}

\begin{proof}
 We have:
 \begin{align*}
  & \ell^\theta_{a,b} + \ell^\theta_{b,c} - \ell^\theta_{a,c} \\
  =\, & d(a_1, b_1) + d(b_1, c_1) - d(a_1, c_1) +  d(a_2, b_2) + d(b_2, c_2) - d(a_2, c_2) \\ & +n\left(\left|\delt{a_1 \in \ooi{b_1}{a_2}} - \delt{b_2 \in \ooi{b_1}{a_2}}\right| + \left|\delt{b_1 \in \ooi{c_1}{b_2}} - \delt{c_2 \in \ooi{c_1}{b_2}}\right| \right. \\&\left.- \left|\delt{a_1 \in \ooi{c_1}{a_2}} - \delt{c_2 \in \ooi{c_1}{a_2}}\right|\right) \\
  =\, & n\left(\delt{b_1 \in \ooi{c_1}{a_1}} + \delt{b_2 \in \ooi{c_2}{a_2}} + \left|\delt{a_1 \in \ooi{b_1}{a_2}} - \delt{b_2 \in \ooi{b_1}{a_2}}\right| \right. \\&\left. + \left|\delt{b_1 \in \ooi{c_1}{b_2}} - \delt{c_2 \in \ooi{c_1}{b_2}}\right| - \left|\delt{a_1 \in \ooi{c_1}{a_2}} - \delt{c_2 \in \ooi{c_1}{a_2}}\right|\right).
 \end{align*}
 If this quantity was negative, we would have $b_1 \in \cci{a_1}{c_1}$ and $b_2 \in \cci{a_2}{c_2}$ and:
 \begin{itemize}
  \item either $a_1 \in \ooi{c_1}{a_2}$ and $c_2 \in \cci{a_2}{c_1}$. As $c_2 \notin \ooi{c_1}{b_2}$, we get $b_1 \notin \ooi{c_1}{b_2}$ so $b_1 \in \cci{b_2}{c_1}$. It is then easy to deduce that $a_1 \in \ooi{b_1}{a_2}$ and $b_2 \notin \ooi{b_1}{a_2}$ and it contradicts the hypothesis.
  \item either $a_1 \in \cci{a_2}{c_1}$ and $c_2 \in \ooi{c_1}{a_2}$. As $a_1 \notin \ooi{b_1}{a_2}$, we get $b_2 \notin \ooi{b_1}{a_2}$ so $b_2 \in \cci{a_2}{b_1}$. It is then easy to deduce that $b_1 \notin \ooi{c_1}{b_2}$ and $c_2 \in \ooi{c_1}{b_2}$ and it contradicts the hypothesis.
 \end{itemize}
 In any case, we reached a contradiction.

 Notice that for any four $i,j,k,l \in \cci{1}{n}$, we have the identities
 \begin{align*}
  \left|\delt{i \in \ooi{k}{j}} - \delt{l \in \ooi{k}{j}}\right| &= \delt{i \in \ooi{k}{j}} \delt{l \in \cci{j}{k}} + \delt{i \in \cci{j}{k}} \delt{l \in \ooi{k}{j}} \\
   &= \delt{i \in \ooi{k}{j}} \delt{l \in \cci{j}{k}} + (1-\delt{i \in \ooi{k}{j}}) (1 - \delt{l \in \cci{j}{k}}) \\
   &= 2 \delt{i \in \ooi{k}{j}} \delt{l \in \cci{j}{k}} + 1 - \delt{i \in \ooi{k}{j}} - \delt{l \in \cci{j}{k}} \\
  \text{and} \quad \left|\delt{i \in \ooi{k}{j}} - \delt{l \in \ooi{k}{j}}\right| &= \delt{i \in \ooi{k}{j}} \delt{l \in \cci{j}{k}} + \delt{i \in \cci{j}{k}} \delt{l \in \ooi{k}{j}} \\
   &= \delt{i \in \ooi{k}{j}} (1-\delt{l \in \ooi{k}{j}}) + (1-\delt{i \in \ooi{k}{j}}) \delt{l \in \ooi{k}{j}} \\
   &= \delt{i \in \ooi{k}{j}} + \delt{l \in \ooi{k}{j}} - 2\delt{i \in \ooi{k}{j}} \delt{l \in \ooi{k}{j}}.
 \end{align*}

 If $a$ is a side of the polygon, we have $a_2 = a_1+1$ and the previous difference becomes (up to a factor $n$):
 \begin{align*}
   & \delt{b_1 \in \ooi{c_1}{a_1}} + \delt{b_2 \in \ooi{c_2}{a_1+1}} + \left|\delt{a_1 \neq b_1} - \delt{b_2 \in \ooi{b_1}{a_1+1}}\right| \\ &+ \left|\delt{b_1 \in \ooi{c_1}{b_2}} - \delt{c_2 \in \ooi{c_1}{b_2}}\right| - \left|\delt{a_1 \neq c_1} - \delt{c_2 \in \ooi{c_1}{a_1+1}}\right| \\
   =\, & \delt{a_1 \in \oci{b_1}{c_1}} + \delt{b_2 \in \ooi{c_2}{a_1+1}} + \delt{a_1 \neq b_1} - \delt{b_2 \in \ooi{b_1}{a_1+1}} \\ &+ 2\delt{b_1 \in \cci{b_2}{c_1}} \delt{c_2 \in \ooi{c_1}{b_2}} + 1 -\delt{b_1 \in \cci{b_2}{c_1}} - \delt{c_2 \in \ooi{c_1}{b_2}} - \delt{a_1 \neq c_1} + \delt{c_2 \in \ooi{c_1}{a_1+1}} \\
   =\, & \delt{a_1 \in \oci{b_1}{c_1}} + 2\delt{b_1 \in \cci{b_2}{c_1}} \delt{c_2 \in \ooi{c_1}{b_2}} + 1 -\delt{b_1 \in \cci{b_2}{c_1}} - \delt{c_2 \in \ooi{c_1}{b_2}} \\&+ \delt{a_1 \in \oci{b_2-1}{c_2-1}} - \delt{a_1 = b_1} - \delt{a_1 \in \oci{b_2-1}{b_1-1}}  + \delt{a_1 = c_1} + \delt{a_1 \in \oci{c_2-1}{c_1-1}} \\
   =\, & \delt{a_1 \in \oci{b_1}{c_1}} + 2\delt{b_1 \in \cci{b_2}{c_1}} \delt{c_2 \in \ooi{c_1}{b_2}} -\delt{b_1 \in \cci{b_2}{c_1}} - \delt{c_2 \in \ooi{c_1}{b_2}} \\&+ \delt{a_1 \in \oci{b_2-1}{c_2-1}} - \delt{a_1 = b_1} + \delt{a_1 \in \oci{b_1-1}{b_2-1}} + \delt{b_1=b_2}  + \delt{a_1 = c_1} + \delt{a_1 \in \oci{c_2-1}{c_1-1}} \\   
   =\, & \delt{a_1 \in \oci{b_1}{c_1}} + 2\delt{b_1 \in \cci{b_2}{c_1}} \delt{c_2 \in \ooi{c_1}{b_2}}  -\delt{b_1 \in \oci{b_2}{c_1}} - \delt{c_2 \in \ooi{c_1}{b_2}} \\&+ \delt{a_1 \in \oci{b_1-1}{c_1-1}} + \delt{b_2 \in \ooi{c_1}{b_1} } + \delt{c_2 \in \ooi{c_1}{b_2}} - \delt{a_1 = b_1} + \delt{a_1 = c_1} \\   
   =\, & 2\delt{a_1 \in \oci{b_1}{c_1}} + 2\delt{b_1 \in \cci{b_2}{c_1}} \delt{c_2 \in \ooi{c_1}{b_2}}.  
 \end{align*}

 The other computations are analogous.
\end{proof}

Let us introduce the $K$-algebras that will play an important role in this paper.

As before, $R = K[X]$. We define the $R$-algebra $R' = K[X,Y]/(YX - Y^2)$. It is in fact an $R$-order of rank $2$ (see Definition \ref{def:order}), and we have the following $R$-isomorphism with a classical Bass order:
$$
 R' \rightarrow R - R := \{(P,Q) \in R^2 \,|\, P-Q \in (X)\}, \quad Y \mapsto (0,X).
$$

 The three indecomposable Cohen-Macaulay $R'$-modules and irreducible morphisms over $R$ appear in each of the two lines of the following commutative diagram:
\[{ \begin{tikzcd}[ampersand replacement=\&, text depth=0pt, text height=\heightof{$A'$},
  column sep={{{{10em,between origins}}}}]
  	(Y) \arrow[out=15,in=165]{r}{\iota} \& R' \arrow[out=15,in=165]{r}{\pi} \arrow[out=-165,in=-15]{l}[swap]{Y} \arrow[equal]{d} \& R'/(Y) \arrow[out=-165,in=-15]{l}[swap]{X-Y} \arrow{d}{X-Y}[swap]{\wr} \\
R'/(X-Y) \arrow[out=15,in=165]{r}[swap]{Y} \arrow{u}{Y}[swap]{\wr} \& R' \arrow[out=15,in=165]{r}[swap]{X-Y} \arrow[out=-165,in=-15]{l}{\pi'} \& (X-Y) \arrow[out=-165,in=-15]{l}{\iota'}
  \end{tikzcd}}
  \]
where $Y$ and $X-Y$ are multiplications, $\pi$ are projections and $\iota$ are natural inclusions.

Finally, we denote by $\cR'$ the algebra $K[u^{\pm 1},v]/(vu-v^2)$ where $R'$ is seen as a subalgebra of $\cR'$ through the inclusion
\begin{equation}
 R' \hookrightarrow \cR', \quad X \mapsto u^{2n}, \quad Y \mapsto v^{2n}. \label{cRp}
\end{equation}

\subsection{Frozen Jacobian algebras are $R$-orders}
\label{ss:Jacobian Alg as Order}
Let $(Q_{\sigma},W_{\sigma},F)$ be an ice quiver with potential arising from a tagged triangulation $\sigma$ as defined in Section \ref{ss:Ice QP} and $e_i$ be the trivial path of length $0$ at vertex $i$.
The main result of this section is that $\Gamma_\sigma := \cpp$ (Definition \ref{frozen JA}) is an $R$-order.

First, we introduce the definition of orders and Cohen-Macaulay modules over orders.
\begin{definition} 
\label{def:order}
Let $S$ be a commutative Noetherian ring of Krull dimension 1.
An $S$-algebra $A$ is called an {\emph{$S$-order}} if it is a finitely generated $S$-module and $\soc_S A=0$. 
 For an $S$-order $A$, a left $A$-module $M$ is called a (maximal) {\emph{Cohen-Macaulay $A$-module}} if it is finitely generated as an $S$-module and $\soc_S M=0$ (or equivalently $\soc_A M=0$). 
  We denote by $\CM A$ the category of Cohen-Macaulay $A$-modules. It is a full exact subcategory of $\mod A$.
\end{definition}  

\begin{remark}
\label{rem:order}
If $S$ is a principal ideal domain (\emph{e.g.} $S=R$) and $M \in \mod S$, then $\soc_S M=0$ if and only if $M$ is free as an $S$-module.
\end{remark}

We refer to \cite{FMO},\cite{RT},\cite{POS} and \cite{CM} for more details about orders and Cohen-Macaulay modules.


The main theorem of this subsection is the following.
\begin{theorem}
\label{thm:frozen Jacobian algebra is order}
The frozen Jacobian algebra $\Gamma_\sigma$ has the structure of an $R$-order.
\end{theorem}

The remaining part of this subsection is devoted to prove Theorem \ref{thm:frozen Jacobian algebra is order}. The strategy is to define a grading on $\Gamma_\sigma$, to prove that the centre $Z(\Gamma_\sigma)$ of $\Gamma_\sigma$ is $R'$ and to give its order structure as an $R'$-module. Notice that the center of Jacobian algebras coming from surfaces without boundary was computed by Ladkani in \cite[Proposition 4.11]{Lad}. 

We describe $\Gamma_\sigma$ in details in Proposition \ref{inductD}. Let us define a grading on $Q_\sigma$ (and $Q'_\sigma$).
 
 \begin{definition}[$\theta$-length]
 Let $a$ and $b$ be sides or tagged arcs of $P^*$ and $\alpha: a \rightarrow b$ be an arrow of $Q'_\sigma$. Let $\theta$ be the value of the oriented angle from $\vec a$ to $\vec b$ taken in $[0, 2\pi)$. We define the $\theta$-length of $\alpha$ by 
  $$\ell^\theta(\alpha) = \frac{n}{\pi} \theta.$$

 The $\theta$-length of arrows extends additively to a map $\ell^\theta$ from paths to integers, defining a grading on $KQ_\sigma$ (and $K Q'_\sigma$) which will also be denoted by $\ell^\theta$. 
   \end{definition}
 
 \begin{remark}
  Using Lemma \ref{lem:easyangle}, we see easily that for any arrow $\alpha: a \rightarrow b$, 
  $\ell^\theta(\alpha) = \ell^\theta_{a,b}$. Indeed, if $a$ and $b$ share a common endpoint, then $0 \leq \ell^\theta_{a,b} < 2n$.
 \end{remark}

 We now prove that for any tagged arcs or sides $a$ and $b$ of $P^*$, the possible $\theta$-lengths of paths from $a$ to $b$ in $Q_\sigma$ does not depend on the triangulation $\sigma$ containing $a$ and $b$. 

 \begin{proposition}
  \label{pro:independence}
  Let $\sigma$ and $\sigma'$ be two different triangulations of the punctured polygon $P^*$. For any two edges $a$ and $b$ common to $\sigma$ and $\sigma'$, the minimal $\theta$-length of paths from $a$ to $b$ in $Q_\sigma$ is the same as the one in $Q_{\sigma'}$. 
 \end{proposition}
\begin{proof}
Any two triangulations can be related by a sequence of flips such that each time we only change one arc in the related triangulation to get another one.
Therefore, without losing generality, we can assume that the two triangulations $\sigma$ and $\sigma'$ have the same arcs except one.
We show the possible differences between $\sigma$ and $\sigma'$ in Figure \ref{flip}. 

\begin{figure}\begin{center}
 \begin{minipage}{6cm}\begin{center}\begin{pdfpic}
\begin{pspicture}(0,-2)(5.5,2.5)
\psframe[linewidth=0.02,dimen=outer](4.08,1.2)(1.58,-1.3)
\psline[linewidth=0.02cm](1.6,1.16)(4.06,-1.28)
\psline[linewidth=0.02cm,linecolor=red,arrowsize=0.05cm 8.0,arrowlength=1.0,arrowinset=0.6]{->}(4.08,-0.1)(2.88,-0.12)
\psline[linewidth=0.02cm,linecolor=red,arrowsize=0.05cm 8.0,arrowlength=1.0,arrowinset=0.6]{->}(2.96,1.18)(4.0,-0.02)
\psline[linewidth=0.02cm,linecolor=red,arrowsize=0.05cm 8.0,arrowlength=1.0,arrowinset=0.6]{->}(2.9,-0.1)(2.94,1.14)
\psline[linewidth=0.02cm,linecolor=red,arrowsize=0.05cm 8.0,arrowlength=1.0,arrowinset=0.6]{->}(2.9,-0.2)(2.9,-1.24)
\psline[linewidth=0.02cm,linecolor=red,arrowsize=0.05cm 8.0,arrowlength=1.0,arrowinset=0.6]{->}(1.62,-0.12)(2.86,-0.14)
\psline[linewidth=0.02cm,linecolor=red,arrowsize=0.05cm 8.0,arrowlength=1.0,arrowinset=0.6]{->}(2.92,-1.26)(1.66,-0.14)
\usefont{T1}{ptm}{m}{n}
\rput(4.3,1.45){$P_k$}
\usefont{T1}{ptm}{m}{n}
\rput(1.45,1.45){$P_i$}
\usefont{T1}{ptm}{m}{n}
\rput(4.3,-1.6){$P_j$}
\usefont{T1}{ptm}{m}{n}
\rput(1.45,-1.6){$P_l$}
\psdots[dotsize=0.15](2.9,2)
\usefont{T1}{ptm}{m}{n}
\rput(3.1,2){$*$}
\psbezier[linewidth=0.02cm,linestyle=dashed,dash=0.16cm 0.16cm](1.58,1.2)(1.6,2.5)(4.06,2.5)(4.08,1.2)
\end{pspicture} 
\end{pdfpic}
\end{center}\end{minipage} $\leftrightarrow$ \begin{minipage}{6cm}\begin{center}\begin{pdfpic}
\begin{pspicture}(0,-2)(5.5,2.5)
\psframe[linewidth=0.02,dimen=outer](4.08,1.2)(1.58,-1.3)
\psline[linewidth=0.02cm](1.58,-1.3)(4.08,1.2)
\psline[linewidth=0.02cm,linecolor=red,arrowsize=0.05cm 8.0,arrowlength=1.0,arrowinset=0.6]{<-}(2.8,-0.08)(2.8,1.15)
\psline[linewidth=0.02cm,linecolor=red,arrowsize=0.05cm 8.0,arrowlength=1.0,arrowinset=0.6]{<-}(4.08,-0.1)(2.88,-0.12)
\psline[linewidth=0.02cm,linecolor=red,arrowsize=0.05cm 8.0,arrowlength=1.0,arrowinset=0.6]{<-}(1.62,-0.12)(2.7,-0.12)
\psline[linewidth=0.02cm,linecolor=red,arrowsize=0.05cm 8.0,arrowlength=1.0,arrowinset=0.6]{<-}(2.8,-1.25)(4.08,-0.12)
\psline[linewidth=0.02cm,linecolor=red,arrowsize=0.05cm 8.0,arrowlength=1.0,arrowinset=0.6]{->}(2.8,-1.25)(2.8,-0.12)
\psline[linewidth=0.02cm,linecolor=red,arrowsize=0.05cm 8.0,arrowlength=1.0,arrowinset=0.6]{<-}(2.8,1.15)(1.62,-0.12)
\usefont{T1}{ptm}{m}{n}
\rput(4.3,1.45){$P_k$}
\usefont{T1}{ptm}{m}{n}
\rput(1.45,1.45){$P_i$}
\usefont{T1}{ptm}{m}{n}
\rput(4.3,-1.6){$P_j$}
\usefont{T1}{ptm}{m}{n}
\rput(1.45,-1.6){$P_l$}
\psdots[dotsize=0.15](2.9,2)
\usefont{T1}{ptm}{m}{n}
\rput(3.1,2){$*$}
\psbezier[linewidth=0.02cm,linestyle=dashed,dash=0.16cm 0.16cm](1.58,1.2)(1.6,2.5)(4.06,2.5)(4.08,1.2)
\end{pspicture} 
\end{pdfpic}
\end{center}\end{minipage}
 
 \begin{minipage}{6cm}\begin{center}\begin{pdfpic}
\begin{pspicture}(0,-2)(5.5,1.8)
\psframe[linewidth=0.02,dimen=outer](4.08,1.2)(1.58,-1.3)
\psline[linewidth=0.02cm](1.6,1.16)(4.06,-1.28)
\psline[linewidth=0.02cm,linecolor=red,arrowsize=0.05cm 8.0,arrowlength=1.0,arrowinset=0.6]{->}(4.08,-0.1)(2.88,-0.12)
\psline[linewidth=0.02cm,linecolor=red,arrowsize=0.05cm 8.0,arrowlength=1.0,arrowinset=0.6]{->}(2.96,1.18)(4.0,-0.02)
\psline[linewidth=0.02cm,linecolor=red,arrowsize=0.05cm 8.0,arrowlength=1.0,arrowinset=0.6]{->}(2.9,-0.1)(2.94,1.14)
\psline[linewidth=0.02cm,linecolor=red,arrowsize=0.05cm 8.0,arrowlength=1.0,arrowinset=0.6]{->}(2.9,-0.2)(2.9,-1.24)
\psline[linewidth=0.02cm,linecolor=red,arrowsize=0.05cm 8.0,arrowlength=1.0,arrowinset=0.6]{->}(1.62,-0.12)(2.86,-0.14)
\psline[linewidth=0.02cm,linecolor=red,arrowsize=0.05cm 8.0,arrowlength=1.0,arrowinset=0.6]{->}(2.92,-1.26)(1.66,-0.14)
\usefont{T1}{ptm}{m}{n}
\rput(4.3,1.45){$P_k$}
\usefont{T1}{ptm}{m}{n}
\rput(1.45,1.45){$P_i$}
\usefont{T1}{ptm}{m}{n}
\rput(4.3,-1.6){$P_j$}
\psdots[dotsize=0.15](1.62,-1.26)
\usefont{T1}{ptm}{m}{n}
\rput(1.45,-1.45){$*$}
\end{pspicture} 
\end{pdfpic}
\end{center}\end{minipage} $\leftrightarrow$ \begin{minipage}{6cm}\begin{center}\begin{pdfpic}
\begin{pspicture}(0,-2)(5.5,1.8)
\psframe[linewidth=0.02,dimen=outer](4.08,1.2)(1.58,-1.3)
\psline[linewidth=0.02cm](1.58,-1.3)(4.08,1.2)
\psline[linewidth=0.02cm,linecolor=red,arrowsize=0.05cm 8.0,arrowlength=1.0,arrowinset=0.6]{<-}(2.8,-0.08)(2.8,1.15)
\psline[linewidth=0.02cm,linecolor=red,arrowsize=0.05cm 8.0,arrowlength=1.0,arrowinset=0.6]{<-}(4.08,-0.1)(2.88,-0.12)
\psline[linewidth=0.02cm,linecolor=red,arrowsize=0.05cm 8.0,arrowlength=1.0,arrowinset=0.6]{<-}(1.62,-0.12)(2.7,-0.12)
\psline[linewidth=0.02cm,linecolor=red,arrowsize=0.05cm 8.0,arrowlength=1.0,arrowinset=0.6]{<-}(2.8,-1.25)(4.08,-0.12)
\psline[linewidth=0.02cm,linecolor=red,arrowsize=0.05cm 8.0,arrowlength=1.0,arrowinset=0.6]{->}(2.8,-1.25)(2.8,-0.12)
\psline[linewidth=0.02cm,linecolor=red,arrowsize=0.05cm 8.0,arrowlength=1.0,arrowinset=0.6]{<-}(2.8,1.15)(1.62,-0.12)
\usefont{T1}{ptm}{m}{n}
\rput(4.3,1.45){$P_k$}
\usefont{T1}{ptm}{m}{n}
\rput(1.45,1.45){$P_i$}
\usefont{T1}{ptm}{m}{n}
\rput(4.3,-1.6){$P_j$}
\psdots[dotsize=0.15](1.62,-1.26)
\usefont{T1}{ptm}{m}{n}
\rput(1.45,-1.45){$*$}
\end{pspicture}
\end{pdfpic}
\end{center}\end{minipage}
 
 \begin{minipage}{6cm}\begin{center}\begin{pdfpic}
\begin{pspicture}(0,-2.7)(3.61,2.8)
\psdots[dotsize=0.15](1.6,0.6)
\psbezier[linewidth=0.02](1.6,2.2)(2.8,1.8)(3.6,-1.0)(1.6,-2.2)
\psbezier[linewidth=0.02](1.6,2.2)(0.0,1.8)(0.0,-1.4)(1.6,-2.2)
\psline[linewidth=0.02cm](1.6,0.6)(1.6,-2.2)
\psbezier[linewidth=0.02](1.6,0.6)(2.4,0.2)(2.0,-1.4)(1.6,-2.2)
\psline[linewidth=0.02cm,linecolor=red,arrowsize=0.05cm 8.0,arrowlength=1.0,arrowinset=0.6]{->}(2.8,0)(2.0,-1.1)
\psline[linewidth=0.02cm,linecolor=red,arrowsize=0.05cm 8.0,arrowlength=1.0,arrowinset=0.6]{->}(2.0,-1.1)(0.4,-0.05)
\psline[linewidth=0.02cm,linecolor=red,arrowsize=0.05cm 8.0,arrowlength=1.0,arrowinset=0.6]{->}(2.8,0)(1.6,-0.5)
\psline[linewidth=0.02cm,linecolor=red,arrowsize=0.05cm 8.0,arrowlength=1.0,arrowinset=0.6]{->}(1.6,-0.5)(0.4,0)
\psbezier[linewidth=0.02cm,linecolor=red,arrowsize=0.05cm 8.0,arrowlength=1.0,arrowinset=0.6]{->}(0.4,0)(1.2,1.2)(2,1.2)(2.8,0)
\usefont{T1}{ptm}{m}{n}
\rput(2.051455,-0.095){$\bowtie$}
\usefont{T1}{ptm}{m}{n}
\rput(1.6,2.5){$P_i$}
\usefont{T1}{ptm}{m}{n}
\rput(1.6,-2.5){$P_j$}
\usefont{T1}{ptm}{m}{n}
\rput(1.4,0.6){$*$}
\end{pspicture} 
\end{pdfpic}
\end{center}\end{minipage} $\leftrightarrow$ \begin{minipage}{6cm}\begin{center}\begin{pdfpic}
\begin{pspicture}(0,-2.7)(3.61,2.8)
\psdots[dotsize=0.15](1.6,0)
\psbezier[linewidth=0.02](1.6,2.2)(2.8,1.8)(3.6,-1.0)(1.6,-2.2)
\psbezier[linewidth=0.02](1.6,2.2)(0.0,1.8)(0.0,-1.4)(1.6,-2.2)
\psline[linewidth=0.02cm](1.6,0)(1.6,-2.2)
\psline[linewidth=0.02cm](1.6,0)(1.6,2.2)
\psbezier[linewidth=0.02cm,linecolor=red,arrowsize=0.05cm 8.0,arrowlength=1.0,arrowinset=0.6]{->}(1.65,-0.9)(1.8,-0.6)(1.8,0.6)(1.65,1.0)
\psbezier[linewidth=0.02cm,linecolor=red,arrowsize=0.05cm 8.0,arrowlength=1.0,arrowinset=0.6]{->}(1.55,0.9)(1.4,0.6)(1.4,-0.6)(1.55,-1.0)
\psline[linewidth=0.02cm,linecolor=red,arrowsize=0.05cm 8.0,arrowlength=1.0,arrowinset=0.6]{->}(1.6,1.1)(2.8,0)
\psline[linewidth=0.02cm,linecolor=red,arrowsize=0.05cm 8.0,arrowlength=1.0,arrowinset=0.6]{->}(2.8,0)(1.6,-1.1)
\psline[linewidth=0.02cm,linecolor=red,arrowsize=0.05cm 8.0,arrowlength=1.0,arrowinset=0.6]{->}(1.6,-1.1)(0.4,0)
\psline[linewidth=0.02cm,linecolor=red,arrowsize=0.05cm 8.0,arrowlength=1.0,arrowinset=0.6]{->}(0.4,0)(1.6,1.1)

\usefont{T1}{ptm}{m}{n}
\rput(1.6,2.5){$P_i$}
\usefont{T1}{ptm}{m}{n}
\rput(1.6,-2.5){$P_j$}
\usefont{T1}{ptm}{m}{n}
\rput(1.3,0){$*$}
\end{pspicture} 
\end{pdfpic}
\end{center}\end{minipage}
 \end{center}
 \caption{Possible flips}
 \label{flip}
\end{figure}
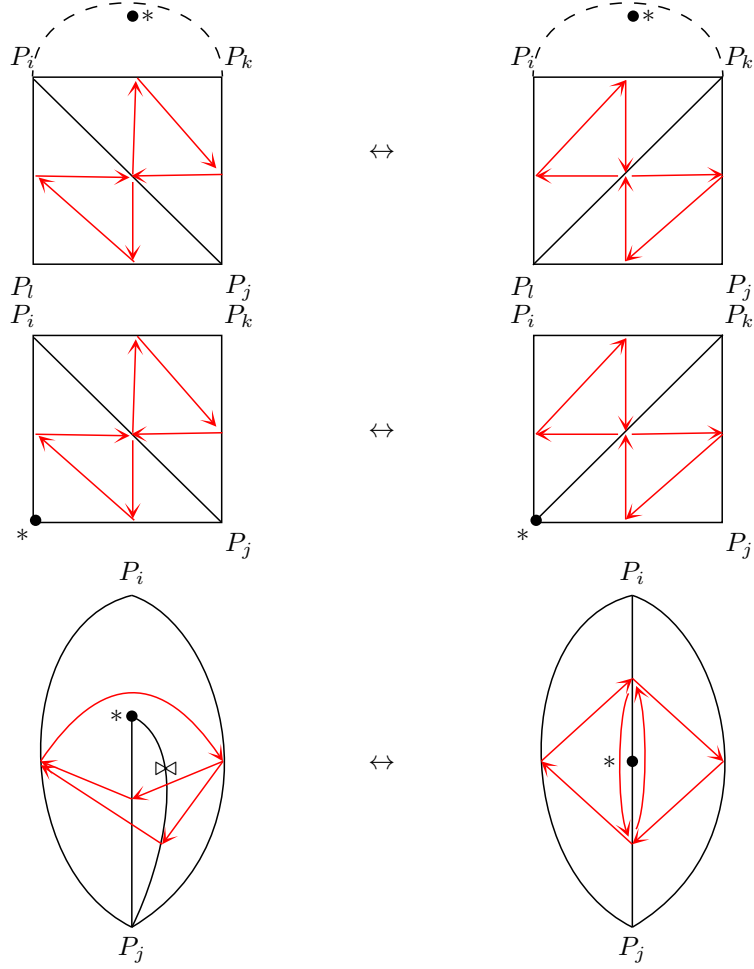

It is sufficient to prove that for any two vertices common to both triangulations, and for any path between them in one triangulation, we can find a path with the same $\theta$-length in the other triangulation. In each case, certain compositions of two arrows in one of the diagrams have to be replaced by one arrow in the other one. We can check case by case that the $\theta$-length of both coincide.

For example, suppose that the triangulations only differ in a square not incident to the puncture as shown at the top of Figure \ref{flip}. Considering the given position of the puncture, the arcs are $(P_i, P_k)$, $(P_i, P_j)$, $(P_i, P_l)$, $(P_j, P_k)$, $(P_l, P_j)$ and $(P_l, P_k)$. And we get that
\begin{align*}& \ell^\theta((P_j, P_k), (P_i, P_j)) + \ell^\theta((P_i, P_j), (P_l, P_j)) \\ =\,& d(j,i) + d(k,j) + d(i,l) = d(j,l) + d(k,j) = \ell^\theta((P_j, P_k), (P_l, P_j))\end{align*}
and the other cases work in the same way. 
\end{proof}

\begin{proposition} \label{pro:pothom}
 The potential $W_\sigma$ (resp. $W'_\sigma$) on $Q_\sigma$ (resp. $Q'_\sigma$) is homogeneous of $\theta$-length $2n$. Thus, $\ell^\theta$ induces a grading on $\Gamma_\sigma$.
\end{proposition}

\begin{proof}
 Consider a true triangle of $\sigma$. Up to cyclic permutation, we can denote its three sides by $a$, $b$, $c$ in clockwise order, satisfying $a_1 = b_2$, $a_2 = c_2$. Moreover, they satisfy either $b_1 = c_1$ and $a_1 \in \ooi{b_1}{a_2}$ (if the triangle is not incident to the puncture) either $b_1 = b_2$ and $c_1 = c_2$ (if the triangle is incident to the puncture). In any case, the $\theta$-length of the clockwise triangle induced by this true triangle is
 \begin{align*}
  &\ell^\theta_{a,b} + \ell^\theta_{b,c} + \ell^\theta_{c,a} \\ =\,& (d(a_1, b_1) + d(b_1, c_1) + d(c_1, a_1)) + (d(a_2, b_2) + d(b_2, c_2) + d(c_2, a_2))
    \\ &+n\left(\left|\delt{a_1 \in \ooi{b_1}{a_2}} - \delt{b_2 \in \ooi{b_1}{a_2}}\right|+ \left|\delt{b_1 \in \ooi{c_1}{b_2}} - \delt{c_2 \in \ooi{c_1}{b_2}}\right| \right. \\ &\left.+ \left|\delt{c_1 \in \ooi{a_1}{c_2}} - \delt{a_2 \in \ooi{a_1}{c_2}}\right|\right) = n + n + n \times 0 = 2 n.
 \end{align*}

 Using flips introduced in the proof of Proposition \ref{pro:independence}, we can transform $\sigma$ to a triangulation $\tau$ that consists of the sides of $P$ and the tagged arcs incident to the puncture (as in Figure \ref{init1}, page \pageref{init1}). Moreover, using the reasoning of Proposition \ref{pro:independence}, flips clearly conserve the $\theta$-length of anticlockwise planar cycles of $Q'_\sigma$ winding around the puncture or vertices of $P^*$. Therefore, it is enough to see that $W'_\tau$\com{ /*(resp. $W'_\tau$)*/} is homogeneous of $\theta$-length $2n$. This is easy to check by calculation. As terms of $W_\sigma$ are terms of $W'_\sigma$, $W_\sigma$ is also homogeneous.
\end{proof}

 \begin{proposition} \label{pro:minthet}
 Let $a$ and $b$ be two edges of $\sigma$ which are not incident to the puncture or tagged in the same way. The minimal $\theta$-length of paths from $a$ to $b$ in $Q_\sigma$ is $\ell^\theta_{a,b}$.
 \end{proposition}
 \begin{proof}
 Let us prove first that there exists a path from $a$ to $b$ with $\theta$-length $\ell^\theta_{a,b}$. Let us do an induction on $\ell^\theta_{a,b}$. If it is $0$, then $a = b$ and the result is obvious. If $a$ and $b$ are both incident to the puncture and not isotopic, the result is immediate (consider the triangulation $\tau$ consisting the sides of $P$ and the plain tagged arcs incident to the puncture). 
 
 Let us suppose that $a$ is not incident to the puncture. We consider four cases:
 \begin{itemize}
  \item Suppose that $b_1, b_2 \neq a_2$. Let us consider an arc $c$ such that $c_1 = a_1$ and $c_2 = a_2+1$ and tagged in the same way as $b$ if $b_1 = b_2$ and $c_1 = c_2$. The arc $c$ is either isotopic to the union of $a$ and a side of the polygon (if $a_2 \neq a_1-1$), either to a part of $a$ (if $a_2 = a_1-1$). In any case, $a$, $b$ and $c$ are compatible so we can choose a triangulation $\sigma'$ containing $a$, $b$ and $c$. In $Q_{\sigma'}$, there is an arrow $\alpha$ from $a$ to $c$ of $\theta$-length $1$, and, as 
 \begin{align*}
  \ell^\theta_{c, b} &= d(c_1, b_1) + d(c_2, b_2) + n\left|\delt{c_1 \in \ooi{b_1}{c_2}}-\delt{b_2 \in \ooi{b_1}{c_2}}\right| \\
   &= d(a_1, b_1) + d(a_2+1, b_2) + n\left|\delt{a_1 \in \ooi{b_1}{a_2+1}}-\delt{b_2 \in \ooi{b_1}{a_2+1}}\right| \\
   &= d(a_1, b_1) + d(a_2, b_2) - 1 + n\left|\delt{a_1 \in \ooi{b_1}{a_2}}-\delt{b_2 \in \ooi{b_1}{a_2}}\right| = \ell^\theta_{a, b} - 1      
 \end{align*}
 we can apply the induction hypothesis: there is a path $w$ from $c$ to $b$ of $\theta$-length $\ell^\theta_{c, b}$ and the path $\alpha w$ has the expected $\theta$-length $\ell^\theta_{a,b}$. Finally, thanks to Proposition \ref{pro:independence}, there is a path of the same $\theta$-length in $\sigma$.
  \item Suppose that $b_1, b_2 \neq a_1$ and $a_2 \neq a_1 + 1$. The same reasoning works with $c_1 = a_1+1$ and $c_2 = a_2$.
  \item Suppose that $b_1 = a_2$ and $b_2 = a_1$. In this case, we construct $\sigma'$ by putting $a$, $b$ and the two tagged arcs connecting the puncture and $a_2$. In $\sigma'$, there is an arrow from $a$ to $b$ which has, by definition, $\theta$-length $\ell^\theta_{a,b}$.
  \item Suppose that $a_2 = a_1 + 1$ and $b_1 \neq a_1$. In this case, put $c = (P_{a_1+1}, P_{a_2+1})$ and complete to a triangulation $\sigma'$ (containing $b$). We have an arrow $\alpha$ from $a$ to $c$ of $\theta$-length $2$. Moreover, 
\begin{align*}
  \ell^\theta_{c, b} &= d(c_1, b_1) + d(c_2, b_2) + n\left|\delt{c_1 \in \ooi{b_1}{c_2}}-\delt{b_2 \in \ooi{b_1}{c_2}}\right| \\
   &= d(a_1+1, b_1) + d(a_2+1, b_2) + n\left|\delt{a_1+1 \in \ooi{b_1}{a_2+1}}-\delt{b_2 \in \ooi{b_1}{a_2+1}}\right| \\
   &= d(a_1, b_1) + d(a_2, b_2) - 2 + n\left|\delt{a_1 \in \ooi{b_1}{a_2}}-\delt{b_2 \in \ooi{b_1}{a_2}}\right| = \ell^\theta_{a, b} - 2      
 \end{align*}
 and the same reasoning as before works.
 \end{itemize}
  The case where $b$ is not incident to the puncture is similar.

 Let us now prove by induction on the $\theta$-length of any path $w$ from $a$ to $b$ that this $\theta$-length is at least $\ell^\theta_{a,b}$. 
 When $a = b$ or if $w$ is an arrow, it is clear. When $a$ and $b$ are different and $w$ is not an arrow, $w$ is the composition of two nonzero paths $w'$ from $a$ to $c$ and $w''$ from $c$ to $b$. 
 By induction hypothesis and Lemma \ref{lem:subadd}, we have 
 \begin{align*}\ell^{\theta}(w) &= \ell^{\theta}(w')+\ell^{\theta}(w'') \geq \ell^\theta_{a,c} + \ell^\theta_{c,b} \geq \ell^\theta_{a,b}. \qedhere\end{align*}
\end{proof} 

In the end of this section, we will prove that $\Gamma_\sigma$ has the structure of an $R$-order and we will specify its structure. More precisely, we prove that the center of $\Gamma_\sigma$ is isomorphic to $R'$ and we realize $\Gamma_\sigma$ as a matrix algebra whose entries are free $R$-submodules of $\cR'$ (see \eqref{cRp}, page \pageref{cRp}).

For each vertex $d$ of $Q'_\sigma$, we will define an element $C_d$ of $e_d \Gamma_\sigma e_d$ as follows.
\begin{itemize}
 \item If $Q'_\sigma$ does not contain Figure \ref{digon} as a subtriangulation, all the planar cycles at $d$ are equivalent (because of commutativity relations). Denote by $C_d$ the common value of these planar cycles in $\Gamma_\sigma$. 
 \item If $Q'_\sigma$ contains Figure \ref{digon} as a subtriangulation, if $d \neq j$, all cycles at $d$ that are planar in the full subquiver $Q'_\sigma \smallsetminus \{j\}$ are equivalent. Denote by $C_d$ their common value in $\Gamma_\sigma$. We define $C_j = 0$.
\end{itemize}

For each frozen vertex $a$ of $Q'_\sigma$, we denote by $E_a$ the big cycle at $a$ passing through each external arrow once. Finally, if $a$ and $b$ are two vertices of $Q'_\sigma$ corresponding to edges which are not incident to the puncture, we denote $a \vdash b$ if $a_2 \in \ooi{b_2}{a_1}$ or $b_1 \in \ooi{b_2}{a_1}$. We have $a \vdash b$ in the following cases:
$$
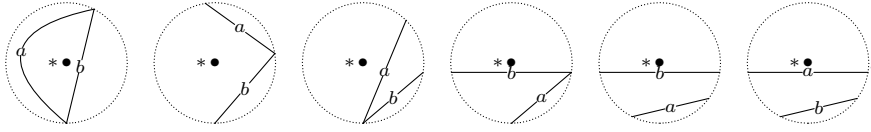
\begin{xy}
{\xypolygon40"A"{~:{(0,-8):}~>.}},
{\ar @{} "A1";"A21" |{\bullet} ^{*}},
{\ar @{-} |{b} "A1";"A18" },
{\ar @/_.8cm/@{-} |{a} "A18";"A1" },
\end{xy}
\quad
\begin{xy}
{\xypolygon40"A"{~:{(0,-8):}~>.}},
{\ar @{} "A1";"A21" |{\bullet} ^{*}},
{\ar @{-} |{b} "A1";"A12" },
{\ar @{-} |{a} "A12";"A22" },
\end{xy}
\quad
\begin{xy}
{\xypolygon40"A"{~:{(0,-8):}~>.}},
{\ar @{} "A1";"A21" |{\bullet} ^{*}},
{\ar @{-} |{a} "A1";"A16" },
{\ar @{-} |{b} "A1";"A10" },
\end{xy}
\quad
\begin{xy}
{\xypolygon40"A"{~:{(0,-8):}~>.}},
{\ar @{} "A1";"A21" |{\bullet} ^{*}},
{\ar @{-} |{a} "A1";"A10" },
{\ar @{-} |{b} "A32";"A10" },
\end{xy}
\quad
\begin{xy}
{\xypolygon40"A"{~:{(0,-8):}~>.}},
{\ar @{} "A1";"A21" |{\bullet} ^{*}},
{\ar @{-} |{a} "A38";"A7" },
{\ar @{-} |{b} "A32";"A10" },
\end{xy}
\quad
\begin{xy}
{\xypolygon40"A"{~:{(0,-8):}~>.}},
{\ar @{} "A1";"A21" |{\bullet} ^{*}},
{\ar @{-} |{b} "A38";"A7" },
{\ar @{-} |{a} "A32";"A10" },
\end{xy}.
$$

Let us fix a grading of $\cR'$ such that $u$ and $v$ have degree $1$. Then, $R'$ is graded as a subalgebra of $\cR'$ ($X$ and $Y$ have degree $2n$). If $a$ and $b$ are two vertices of $Q'_\sigma$, we consider a graded $R'$-submodule $A_{a,b}$ of $\cR'$ (that is also free over $R$) defined in the following way:
\begin{itemize}
 \item $A_{a,b} = 0$ if $a$ and $b$ are incident to the puncture and tagged differently (as $i$ and $j$ in Figure \ref{digon});
 \item $A_{a,b} = u^{\ell^\theta_{a,b}-1} v R'$ if one of $a$ and $b$ is incident to the puncture and plain and the other one either incident to the puncture and plain either not incident to the puncture;
 \item $A_{a,b} = u^{\ell^\theta_{a,b}-1} (u-v) R'$ if one of $a$ and $b$ is incident to the puncture and notched and the other one either incident to the puncture and notched either not incident to the puncture;
 \item $A_{a,b} = u^{\ell^\theta_{a,b}} R' + v^{\ell^\theta_{a,b}} R'$ if $a$ and $b$ are not incident to the puncture and $a \vdash b$;
 \item $A_{a,b} = u^{\ell^\theta_{a,b}} R'$ if $a$ and $b$ are not incident to the puncture and $a \notvdash b$.
\end{itemize}
It is an easy consequence of Lemma \ref{lem:subadd} that $A := (A_{a,b})_{a,b \in Q'_{\sigma,0}}$ is an $R$-subalgebra of the matrix algebra $\M_{Q'_{\sigma,0}}(\cR')$.

\begin{proposition} \label{inductD}
 There exists an isomorphism of graded algebras $\phi_\sigma: R'\rightarrow Z(\Gamma_\sigma)$ ($Z(\Gamma_\sigma)$ is graded by $\theta$-length). 
%
%
%
%
 Moreover, for the induced $R'$-algebra structure of $\Gamma_\sigma$, there is an isomorphism of graded $R'$-algebras $\psi_\sigma: A \to \Gamma_\sigma$ induced by isomorphisms of graded $R'$-modules
 $$\psi_\sigma^{a,b} : A_{a,b} \xrightarrow{\sim} e_{a} \Gamma_\sigma e_{b}$$
 ($\Gamma_\sigma$ is graded by $\theta$-length).

%

  Finally, the following properties are satisfied:
 \begin{enumerate}[$(i)_{\sigma}$]
  \item for each frozen vertex $a$ of $Q_\sigma$, 
   $$e_a \phi_\sigma(X) = \phi_\sigma(X) e_a = E_a$$
  \item for each vertex $a$ of $Q_\sigma$, 
   \begin{align*}
    & e_a \phi_\sigma(Y) = \phi_\sigma(Y) e_a \\ =\, & \left\{\begin{array}{ll}e_a \phi_\sigma(X)-C_a & \text{if $\sigma$ has no plain arc incident to the puncture,} \\ C_a & \text{else;}\end{array} \right. \end{align*}
  \item for any pair of frozen vertices $a$ and $b$, $\psi_\sigma^{a,b}(u^{\ell^\theta_{a,b}})$ is equivalent to the shortest path among paths consisting of external arrows from $a$ to $b$;
  \item for any pair of frozen vertices $a$ and $b$ such that $a$ follows immediately $b$ in the anticlockwise order ($b_2 = a_1$), let $s_{a,b}$ be the path from $a$ to $b$ whose composition with the external arrow $b \rightarrow a$ is the anticlockwise external cycle winding around $a_1$; then we have
   $$\psi_\sigma^{a,b}(v^{\ell^\theta_{a,b}}) = \left\{\begin{array}{ll} \psi_\sigma^{a,b}(u^{\ell^\theta_{a,b}}) - s_{a,b} & \text{if $\sigma$ has no plain arc incident to the puncture} \\ s_{a,b} & \text{else;}\end{array} \right. $$
  \item for any external arrow $\alpha$ of $Q'_\sigma$, and any $w \in \Gamma_\sigma$, if $\alpha w = 0$ then $e_{t(\alpha)} w = 0$ and if $w \alpha = 0$ then  $w e_{s(\alpha)} = 0$.
 \end{enumerate}

\end{proposition}

First of all, when the triangulation has only notched arcs incident to the puncture, the situation is similar to the fully plain one. Then, up to applying the $R$-automorphism of $R'$ given by $Y \mapsto X-Y$ and the $K[u^{\pm 1}]$-automorphism of $\cR'$ given by $v \mapsto u-v$, both results are equivalent (note that this pair of automorphisms commutes with the inclusion $R' \subset \cR'$). From now on, we will only look at cases where triangulations have at most one notched arc incident to the puncture.

Remark that $(v)_\sigma$ is in fact implied by the rest of the proposition. Indeed, suppose that $\alpha$ is an external arrow from a vertex $a$ to a vertex $b$ and suppose that $w \in \Gamma_\sigma$ satisfies $\alpha w = 0$. Without loss of generality, we can suppose that $w = e_{b} w e_{c}$ for some vertex $c$ of $Q'_\sigma$. Thus, there is an element $p \in A_{b,c}$ such that $w = \psi_\sigma^{b,c}(p)$. Then, thanks to $(i)_\sigma$ and the multiplicativity of $\psi_\sigma$, we have
$$0 = E_{b} w = \phi_\sigma(X) e_{b} \psi_\sigma^{b,c}(p) = \phi_\sigma(X) \psi_\sigma^{b,c}(p) =  \psi_\sigma^{b,c}(Xp)$$
and, as $\psi_\sigma^{b,c}$ is injective, we get $Xp = 0$ in $A_{b,c}$. Since $A_{b,c}$ is free over $R \subset R'$, it follows that $p = 0$ and therefore $w = 0$. The other case is dealt in the same way.

If we supposed that the existence of $\phi_\sigma$ as a morphism of algebra and $(ii)_\sigma$ are proved, we can easily deduced that $\phi_\sigma$ is graded. Indeed, we proved in Proposition \ref{pro:pothom} that $\ell^\theta(C_k) = 2n$. Thus, thanks to $(ii)_\sigma$, if $\sigma$ has at most one notched edge incident to the puncture, $\phi_\sigma(Y)$ is homogeneous of degree $2n$. Moreover, as $\phi_\sigma$ is a morphism of algebra, we have $\phi_\sigma(X) \phi_\sigma(Y) = \phi_\sigma(Y) \phi_\sigma(Y)$ so $\phi_\sigma(X)$ is also homogeneous of degree $2n$. 

It is then automatic that $\psi_\sigma^{a,b}$ is graded for any vertices $a$ and $b$ (under the hypothesis that $\psi_\sigma^{a,b}$ is an isomorphism of $R'$-modules). Indeed, $\ell^{\theta}_{a,b}$ is by definition the minimal $\theta$-length of a path from $a$ to $b$. 

The strategy for the rest of the proposition is to do an induction on $n$. We start by proving the proposition for two families of initial cases.

\begin{lemma} \label{lemdig}
 Suppose that $n = 2$ and the triangulation $\sigma$ consists of one plain and one notched arc incident to the puncture as in Figure \ref{digon}. Then Proposition \ref{inductD} is satisfied.
\end{lemma}

\begin{proof}
 In this case, the quiver $Q'_\sigma$ is
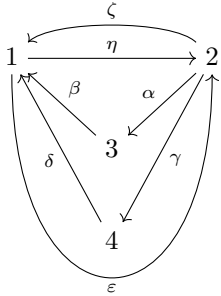
\[{ \begin{tikzcd}[ampersand replacement=\&]
	1 \arrow{rr}{\eta} \arrow[out=-90,in=-90,looseness=3.5]{rr}[swap]{\epsilon}
 \& \& 2 \arrow{dl}[swap]{\alpha} \arrow{ddl}{\gamma} \arrow[out=135,in=45,looseness=.5]{ll}[swap]{\zeta} \\
	\& 3 \arrow{ul}[swap]{\beta} \& \\
        \& 4 \arrow{uul}{\delta} \& \\
  \end{tikzcd}}
  \]
 and the potential is $W'_\sigma = \eta \alpha \beta + \eta \gamma \delta - \gamma \delta \epsilon - \eta \zeta$ so the relations are $\beta \eta = \eta \alpha = 0$, $\delta \eta = \delta \epsilon$, $\eta \gamma = \epsilon \gamma$ and $\zeta = \alpha \beta + \gamma \delta$ ($3$ corresponds to the notched arc and $4$ to the plain one). Using these relations we get 
\begin{claim} Any path different from $\eta$ can be, up to the relations, expressed in a unique way without subpaths of the form $\eta$, $\zeta$, $\beta \epsilon \gamma$ or $\delta \epsilon \alpha$ (the last two are $0$ in $\Gamma_\sigma$). \label{uniqp} \end{claim}
The $\theta$-lengths are given by $\ell^\theta(\alpha) = \ell^\theta(\beta) = \ell^\theta(\gamma) = \ell^\theta(\delta) = 1$ and $\ell^\theta(\epsilon) = \ell^\theta(\zeta) = \ell^\theta(\eta) = 2$.

 Let us prove that there is an isomorphism
 \begin{align*}
  \phi_\sigma: R' & \rightarrow Z(\Gamma_\sigma) \\
       X & \mapsto \epsilon \zeta + \zeta \epsilon + \beta \epsilon \alpha + \delta \epsilon \gamma\\
       Y & \mapsto \epsilon \gamma \delta + \gamma \delta \epsilon + \delta \epsilon \gamma.
 \end{align*}
 It is easy to see that the two elements $\phi_\sigma(X)$ and $\phi_\sigma(Y)$ commute with all arrows so the image is included in the center. Moreover, we have
 \begin{align*}
  \phi_\sigma(X) \phi_\sigma(Y) &= \epsilon \alpha \beta \epsilon \gamma \delta + \epsilon \gamma \delta \epsilon \gamma \delta + \alpha \beta \epsilon \gamma \delta \epsilon + \gamma \delta \epsilon \gamma \delta \epsilon + \delta \epsilon \gamma \delta \epsilon \gamma \\
  &= \epsilon \gamma \delta \epsilon \gamma \delta + \gamma \delta \epsilon \gamma \delta \epsilon + \delta \epsilon \gamma \delta \epsilon \gamma = \phi_\sigma(Y)^2
 \end{align*}
 so $\phi_\sigma$ is a morphism. Notice that
 \begin{align*}\phi_1 = e_1 \phi_\sigma e_1: R'& \rightarrow e_1 Z(\Gamma_\sigma) e_1 \\
	X & \mapsto \epsilon \alpha \beta + \epsilon \gamma \delta \\
        Y & \mapsto \epsilon \gamma \delta
 \end{align*}
 is an isomorphism. Indeed, the surjectivity comes from Claim \ref{uniqp}. For the injectivity, notice that every element of $R'$ can be written (uniquely) in the form $P + Y Q$ where $P$ and $Q$ are polynomials in the variable $X$. Then, as $\epsilon \alpha \beta \epsilon \gamma \delta = \epsilon \gamma \delta \epsilon \alpha \beta = 0$,
 $$\phi_1(P + Y Q) = P(\epsilon \alpha \beta) + P(\epsilon \gamma \delta) - P(0)e_1 + \epsilon \gamma \delta Q(\epsilon \gamma \delta).$$
 If $ \phi_1(P + Y Q) = 0$ then $\beta \phi_1(P + Y Q) = \beta P(\epsilon \alpha \beta) = 0$. As $\Gamma_\sigma / (e_4, \eta)$ is a path algebra (all relations are in the ideal $(e_4, \eta)$ of $K Q'_\sigma$ except $\zeta = \alpha \beta + \gamma \delta$), we get $P = 0$. Then $\epsilon \gamma \delta Q (\epsilon \gamma \delta) = 0$. As $\Gamma_\sigma / (e_3, \eta - \epsilon)$ is a path algebra (all relations are in the ideal $(e_3, \eta - \epsilon)$ of $K Q'_\sigma$ except $\zeta = \alpha \beta + \gamma \delta$), we get $Q = 0$. Thus $\phi_1$ is injective. We deduce immediately that $\phi_\sigma$ is also injective. 

 For the surjectivity of $\phi_\sigma$, take an element $z$ of $Z(\Gamma_\sigma)$. Using Claim \ref{uniqp}, it is immediate that we can write
 $$z = P_1(\epsilon \alpha \beta) + Q_1(\epsilon \gamma \delta) + P_2(\alpha \beta \epsilon) + Q_2(\gamma \delta \epsilon) + P_3(\beta \epsilon \alpha) + Q_4(\delta \epsilon \gamma)$$
 where $Q_1$ and $Q_2$ have no constant terms (we take the convention that $(\epsilon \alpha \beta)^0 = e_1$, $(\alpha \beta \epsilon)^0 = e_2$, $(\beta \epsilon \alpha)^0 = e_3$ and $(\delta \epsilon \gamma)^0 = e_4$).
 Using the identity $\alpha z = z \alpha$, as $\delta \epsilon \alpha = 0$, we get $P_2(\alpha \beta \epsilon) \alpha = \alpha P_3(\beta \epsilon \alpha)$ and, thanks to the grading by $\ell^\theta$, $P_2 = P_3$. In the same way, $\beta z = z \beta$ implies $P_1 = P_3$, $\gamma z = z \gamma$ implies $Q_2 = Q_4$ and $\delta z = z \delta$ implies $Q_1 = Q_4$. So $z = P_1(\epsilon \alpha \beta + \alpha \beta \epsilon + \beta \epsilon \alpha) + Q_1(\epsilon \gamma \delta + \gamma \delta \epsilon + \delta \epsilon \gamma) = \phi_\sigma(P_1(X-Y) + Q_1(Y))$.


 It is an easy observation, using Claim \ref{uniqp}, that $e_1 Z(\Gamma_\sigma) e_i = e_i \Gamma_\sigma e_i$ for every $i \in Q'_{\sigma,0}$. This permits to compute, together with $\theta$-lengths given at the beginning, the following isomorphisms of $R'$-modules (denoted $\psi^{a,b}_\sigma$) from $A_{a,b}$ to $e_a \Gamma_\sigma e_b$ where $a, b \in \cci{1}{4}$:
 $$\begin{array}{|c|c|c|c|c|}
      \hline
      \text{\backslashbox{$a$}{$b$}} & 1  & 2 & 3 & 4 \\
      \hline
     1 & 1 \mapsto e_1 & u^2 \mapsto \epsilon, v^2 \mapsto \eta & (u-v)^3 \mapsto \epsilon \alpha & v^3 \mapsto \epsilon \gamma \\
      \hline
     2 & u^2 \mapsto \zeta, v^2 \mapsto \gamma \delta & 1 \mapsto e_2 & u-v \mapsto \alpha & v \mapsto \gamma \\
      \hline
     3 & u-v \mapsto \beta & (u-v)^3 \mapsto \beta \epsilon & u^{-1} (u-v) \mapsto e_3 & 0 \\
      \hline 
     4 & v \mapsto \delta & v^3 \mapsto \delta \epsilon & 0 & u^{-1}v \mapsto e_4 \\
     \hline
   \end{array} $$
(note that $1 \vdash 2$ and $2 \vdash 1$).
The points $(i)_\sigma$ to $(iv)_\sigma$ are easy to check. The multiplicativity can be checked case by case easily. 
\end{proof}

\begin{lemma} \label{lemtart}
 Suppose that the triangulation consists of the (plain) arcs connecting each vertex of the polygon with the puncture. Then Proposition \ref{inductD} is satisfied.
\end{lemma}

\begin{proof}
 For each $i$ from $1$ to $n$, let us denote by $i$ the arc from $P_i$ to $P_{i+1}$ and by $i'$ the arc from $P_i$ to the puncture. Let us call $\alpha_i$ the arrow of $Q_\sigma$ from $i'$ to $(i+1)'$, $\beta_i$ the arrow from $(i+1)'$ to $i$, $\gamma_i$ the arrow from $i$ to $i'$ and $\delta_i$ the arrow from $i-1$ to $i$ (modulo $n$) (Figure \ref{init1}). Remark that $\ell^\theta(\alpha_i) = \ell^\theta(\delta_i) = 2$, $\ell^\theta(\beta_i) = \ell^\theta(\gamma_i) = n-1$. The relations in $\Gamma_\sigma$ are $\beta_i \gamma_i = \alpha_{i+1} \alpha_{i+2} \cdots \alpha_{i-2} \alpha_i$, $\gamma_i \alpha_i = \delta_{i+1} \gamma_{i+1}$ and $\alpha_i \beta_i = \beta_{i-1} \delta_i$ for all $i \in \cci{1}{n}$.

 \begin{figure}
  \begin{center}
{
\begin{pdfpic}
\begin{pspicture}(0,-1.8)(6.2629104,3.2)
\psdots[dotsize=0.15](2.9610157,0.02)
\psline[linewidth=0.02cm](2.9610157,2.42)(0.96101564,0.82)
\psline[linewidth=0.02cm](2.9610157,2.42)(2.9610157,0.02)
\psline[linewidth=0.02cm](0.96101564,0.82)(2.9610157,0.02)
\psline[linewidth=0.02cm](2.9610157,2.42)(4.9610157,0.82)
\psline[linewidth=0.02cm](4.9610157,0.82)(2.9610157,0.02)
\psbezier[linewidth=0.02cm,linestyle=dashed,dash=0.16cm 0.16cm](0.96101564,0.82)(0.5610156,-3.18)(5.7610154,-2.78)(4.9610157,0.82)
\psline[linewidth=0.02cm,arrowsize=0.05cm 8.0,arrowlength=1.0,arrowinset=0.6]{->}(2.0810156,1.68)(2.9410157,1.18)
\psline[linewidth=0.02cm,arrowsize=0.05cm 8.0,arrowlength=1.0,arrowinset=0.6]{->}(2.9610157,1.14)(1.9810157,0.44)
\psline[linewidth=0.02cm,arrowsize=0.05cm 8.0,arrowlength=1.0,arrowinset=0.6]{->}(1.9410156,0.46)(2.0810156,1.64)
\psline[linewidth=0.02cm,arrowsize=0.05cm 8.0,arrowlength=1.0,arrowinset=0.6]{->}(2.9810157,1.16)(4.0010157,1.56)
\psline[linewidth=0.02cm,arrowsize=0.05cm 8.0,arrowlength=1.0,arrowinset=0.6]{->}(4.0010157,1.54)(3.9210157,0.42)
\psline[linewidth=0.02cm,arrowsize=0.05cm 8.0,arrowlength=1.0,arrowinset=0.6]{->}(3.9210157,0.38)(3.0010157,1.12)
\psbezier[linewidth=0.02cm,arrowsize=0.05cm 8.0,arrowlength=1.0,arrowinset=0.6]{->}(4.017507,1.6)(4.4810157,3.48)(1.5810156,3.46)(2.0644464,1.7322034)
\usefont{T1}{ptm}{m}{n}
\rput(4.4624707,1.645){$i-1$}
\rput(1.8624707,1.745){$i$}
\usefont{T1}{ptm}{m}{n}
\rput(4.0524708,0.065){$(i-1)'$}
\usefont{T1}{ptm}{m}{n}
\rput(1.7824707,0.055){$(i+1)'$}
\usefont{T1}{ptm}{m}{n}
\rput(2.5824707,0.645){$\alpha_i$}
\usefont{T1}{ptm}{m}{n}
\rput(3.3524706,0.505){$\alpha_{i-1}$}
\usefont{T1}{ptm}{m}{n}
\rput(2.6124707,1.155){$i'$}
\usefont{T1}{ptm}{m}{n}
\rput(1.8324707,1.105){$\beta_i$}
\usefont{T1}{ptm}{m}{n}
\rput(3.3324708,1.605){$\beta_{i-1}$}
\usefont{T1}{ptm}{m}{n}
\rput(2.5624706,1.555){$\gamma_i$}
\usefont{T1}{ptm}{m}{n}
\rput(4.3024707,1.000){$\gamma_{i-1}$}
\usefont{T1}{ptm}{m}{n}
\rput(5.3724708,0.825){$P_{i-1}$}
\usefont{T1}{ptm}{m}{n}
\rput(2.9224708,2.585){$P_i$}
\usefont{T1}{ptm}{m}{n}
\rput(0.55247073,0.825){$P_{i+1}$}
\usefont{T1}{ptm}{m}{n}
\rput(3.3824708,3.155){$\delta_i$}
\end{pspicture} 
\end{pdfpic}
}
  \end{center}
  \caption{Initial case}
  \label{init1}
 \end{figure}

 Notice that any path is equivalent to a path containing only arrows of type $\delta$ or to a path containing no arrow of type $\delta$. Then, up to equivalence, a path containing no arrow of type $\delta$ can be supposed not to contain any arrow of type $\gamma$, except maybe at the beginning and not to contain any arrow of type $\beta$ except maybe at the end. To summarize:
 \begin{claim}
  \label{uniq2} Any path of $Q'_\sigma$ is equivalent to a path of the form $$\delta_i \delta_{i+1} \dots \delta_j \quad \text{or} \quad \gamma_i^\mu \alpha_i \alpha_{i+1} \dots \alpha_j \beta_j^\nu,$$ where $\mu, \nu \in \{0,1\}$.
 \end{claim}

For $i \in \cci{1}{n}$, we denote $E_i = \delt{i+1} \delt{i+2} \dots \delta_i$, $C_i = \gamma_i \beta_{i-1} \delta_i = \gamma_i \alpha_i \beta_i = \delt{i+1} \gamma_{i+1} \beta_i$ and $C_{i'} = \beta_{i-1} \delta_i \gamma_i = \beta_{i-1} \gamma_{i-1} \alpha_{i-1} = \alpha_i \alpha_{i+1} \dots \alpha_{i-1} = \alpha_i \beta_i \gamma_i$. Finally, we denote
$$E = \sum_{i = 1}^n (E_i + C_{i'}) \quad \text{and} \quad C = \sum_{i = 1}^n (C_i + C_{i'}).$$

 Let us prove that there is an isomorphism of algebras given by
 \begin{align*}
  \phi_\sigma: R' = K[X, Y]/(YX - Y^2) & \rightarrow Z(\Gamma_\sigma) \\
   X & \mapsto E \\
   Y & \mapsto C.
 \end{align*}

We get easily from the relations that $C_{i'} \alpha_i = \alpha_i C_{(i+1)'}$, $C_{(i+1)'} \beta_i = \beta_i C_i = \beta_i E_i$, $C_i \gamma_i = E_i \gamma_i = \gamma_i C_{i'}$, $\delta_i E_i = E_{i-1} \delta_i$, $\delta_i C_i = C_{i-1} \delta_i$, $\beta_i C_i = \beta_i E_i$, $C_i \gamma_i = E_i \gamma_i$ and $C_i E_i = C_i^2$.
Therefore $C$ and $E$ are in the center of $\Gamma_\sigma$ and $\phi_\sigma$ is a morphism of algebras.

Any element of $R'$ can be written as $P(X) + Y Q(Y)$ where $P$ and $Q$ are two polynomials. If this element is in $\ker \phi_\sigma$, $P(E) + C Q(C) = 0$. Notice now that a path which contains only arrows of type $\delta$ is not related to any other path by the relations. Thus, we should have $P(E) = 0$ and then $P = 0$ as a polynomial (paths appearing in $C Q(C)$ contain arrows which are not of type $\delta$). Then $C Q(C) = 0$. Powers of $C$ have different $\theta$-lengths so $Q = 0$ and finally $\ker \phi_\sigma = 0$.

Let us prove that $\phi_\sigma$ is surjective. Let $z \in Z(\Gamma_\sigma)$. Using Claim \ref{uniq2} and $E_i C_i = C_i^2$, we can write
$$z = \sum_{i = 1}^n \left[P_i(E_i) + Q_i(C_i) + S_i(C_{i'})\right]$$
for some polynomials $P_i$, $Q_i$ and $S_i$ where $Q_i$ has no constant term (we take the convention that $E_i^0 =e_i$ and $C_{i'}^0 = e_{i'}$). For any $i$, $z \alpha_i = \alpha_i z$ implies that $S_i(C_{i'}) \alpha_i = \alpha_i S_{i+1}(C_{(i+1)'})$, so using the grading by $\theta$-length, $S_i = S_{i+1}$. In the same way, $z \beta_i = \beta_i z$ implies $S_{i+1}(C_{(i+1)'}) \beta_i = \beta_i(P_i(E_i) + Q_i(C_i)) = (P_i(C_{(i+1)'}) + Q_i(C_{(i+1)'}) \beta_i$ because $\beta_i E_i = \beta_i C_i = C_{(i+1)'} \beta_i$ so we get $S_{i+1} = P_i + Q_i$. Finally, $z \delta_i = \delta_i z$ implies $(P_{i-1}(E_{i-1}) + Q_{i-1}(C_{i-1})) \delta_i = \delta_i (P_i(E_i) + Q_i(C_i))$. As already observed before, there are no relations between paths containing only $\delta$'s and other paths. Thus, $P_{i-1}(E_{i-1}) \delta_i = \delta_i P_i(E_i)$ and $Q_{i-1}(C_{i-1}) \delta_i = \delta_i Q_i(C_i)$ and using $\theta$-length, $P_{i-1} = P_i$ and $Q_{i-1} = Q_i$. Finally, we get 
 $$z = P_1\left(\sum_{i = 1}^n E_i\right) + Q_1\left(\sum_{i = 1}^n C_i\right) + (P_1+Q_1)\left(\sum_{i = 1}^n C_{i'}\right) = \phi_\sigma(P_1(X) + Q_1(Y))$$
 so $\phi_\sigma$ is surjective.

%

 Let now $i,j$ be two frozen vertices. Notice that $$i \vdash j \Leftrightarrow i+1 \in \ooi{j+1}{i} \text{ or } j \in \ooi{j+1}{i} \Leftrightarrow j = i-1.$$ The following maps are isomorphisms of graded $R'$-modules:
 \begin{align*}
  \psi_\sigma^{i,j}: u^{2d(i,j)}  R' & \rightarrow e_i \Gamma_\sigma e_j & (j \neq i-1)\\
  u^{2d(i,j)} & \mapsto \delt{i+1} \delt{i+2} \dots \delta_j; \\ \\
  \psi_\sigma^{i,i-1}: u^{2(n-1)} R' + v^{2(n-1)} R' & \rightarrow e_i \Gamma_\sigma e_{i-1} \\
   u^{2(n-1)} & \mapsto \delt{i+1} \delt{i+2} \dots \delt{i-1} \\
   v^{2(n-1)} & \mapsto \gamma_i \beta_{i-1}; \\ \\
  \psi_\sigma^{i,j'}: v^{2d(i,j)+n-1} R' & \rightarrow e_i \Gamma_\sigma e_{j'} \\
   v^{2d(i,j)+n-1} & \mapsto \gamma_i \alpha_i \dots \alpha_{j-1}; \\ \\
  \psi_\sigma^{i',j}: v^{2d(i,j+1)+n-1} R' & \rightarrow e_{i'} \Gamma_\sigma e_j \\
   v^{2d(i,j+1)+n-1} & \mapsto \alpha_i \dots \alpha_j \beta_j; \\ \\
  \psi_\sigma^{i',j'}: v^{2d(i,j)} R' & \rightarrow e_{i'} \Gamma_\sigma e_{j'} & (j \neq i)\\
   v^{2d(i,j)} & \mapsto \alpha_i \dots \alpha_{j-1}; \\ \\
  \psi_\sigma^{i',i'}: u^{-1} v R' & \rightarrow e_{i'} \Gamma_\sigma e_{i'} \\
   u^{-1} v & \mapsto e_{i'}.
 \end{align*}

 The argument mainly relies on Claim \ref{uniq2}. Let us for example look at the second case.
 It is easy to check that $$\psi_\sigma^{i,i-1}\left(v^{2(n-1)}\right)(E - C) = 0 \quad \text{and} \quad \psi_\sigma^{i,i-1}\left(u^{2(n-1)}-v^{2(n-1)}\right) C = 0,$$ so $\psi_\sigma^{i,i-1}$ is a morphism. Moreover, if an element $u^{2(n-1)} P(X) + v^{2(n-1)} Q(Y)$ is mapped to $0$ by this map, using the same kind of analysis than before, we prove that $P = Q = 0$ so the map is injective. For the surjectivity, notice that, for any path that does not contain $\delta$ from $i$ to $i-1$, different from $\gamma_i \beta_{i-1}$, in the form given by Claim \ref{uniq2}, we can write
 \begin{align*} \gamma_i \alpha_i \dots \alpha_{i-1} \beta_{i-1} &= \delt{i+1} \gamma_{i+1} \alpha_{i+1} \dots \alpha_{i-1} \beta_{i-1} = \dots \\ &= \delt{i+1} \dots \delt{i-1} \gamma_{i-1} C_{i-1}^k \alpha_{i-1} \beta_{i-1} = \delt{i+1} \dots \delt{i-1} C_{i-1}^{k+1}\end{align*}
 so the map is surjective.

 The multiplicativity can be checked case by case. For example, if $i, j, k \in \cci{1}{n}$ and $j \in \ooi{i}{k}$, we have
 \begin{align*}
  \psi_\sigma^{i,j'}(v^{2d(i,j)+n-1}) \psi_\sigma^{j',k}(v^{2d(j,k+1)+n-1}) &= \gamma_i \alpha_i \alpha_{i+1} \dots \alpha_{j-1} \alpha_j \dots \alpha_{k} \beta_k \\ &= \delt{i+1} \gamma_{i+1} \alpha_{i+1} \dots \alpha_k \beta_k \\ &= \dots = \delt{i+1} \dots \delta_k \gamma_k \alpha_k \beta_k \\ &= \delt{i+1} \dots \delta_k C = \psi_\sigma^{i,k}(u^{2d(i,k)} Y) \\ &= \psi_\sigma^{i,k}(v^{2d(i,k) + 2n}) \\ &= \psi_\sigma^{i,k}(v^{2d(i,j)+n-1}v^{2d(j,k+1)+n-1})
 \end{align*}

 The points $(i)_\sigma$ to $(iv)_\sigma$ are easy to check (and recall that $(v)_\sigma$ is a consequence of them). 
\end{proof}

\begin{proof}[Proof of Proposition \ref{inductD}]
 Let us suppose that the result is proved for all triangulations of polygons with $n-1$ vertices and that there is a \emph{corner triangle} $P_l P_{l+1} P_{l+2}$ in the triangulation as follows:
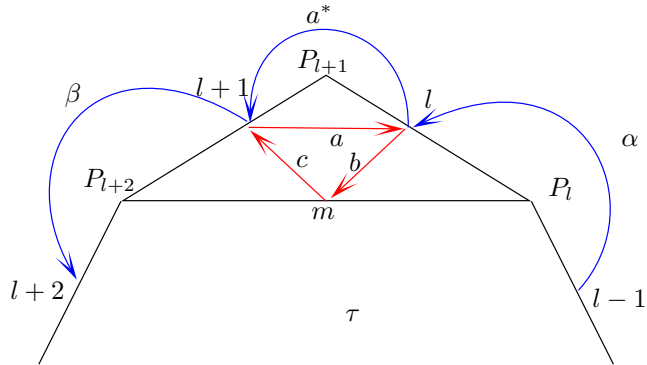
\begin{center}
{
\begin{pdfpic}
\begin{pspicture}(0,-2.48)(8.881894,2.48)
\definecolor{color42}{rgb}{1.0,0.0,0.2}
\pstriangle[linewidth=0.016,dimen=outer](4.22,-0.3)(5.4,1.68)
\psline[linewidth=0.016cm](1.52,-0.3)(0.44,-2.46)
\psline[linewidth=0.016cm](6.92,-0.3)(8.0,-2.46)
\usefont{T1}{ptm}{m}{n}
\rput(4.181455,1.55){$P_{l+1}$}
\usefont{T1}{ptm}{m}{n}
\rput(7.3114552,-0.135){$P_{l}$}
\usefont{T1}{ptm}{m}{n}
\rput(1.3714551,-0.035){$P_{l+2}$}
\psline[linewidth=0.016cm,linecolor=red,arrowsize=0.07cm 8.0,arrowlength=2.0,arrowinset=0.6]{->}(3.2,0.68)(5.24,0.66)
\psline[linewidth=0.016cm,linecolor=red,arrowsize=0.07cm 8.0,arrowlength=2.0,arrowinset=0.6]{->}(5.26,0.66)(4.28,-0.26)
\psline[linewidth=0.016cm,linecolor=red,arrowsize=0.07cm 8.0,arrowlength=2.0,arrowinset=0.6]{->}(4.22,-0.3)(3.22,0.64)
\usefont{T1}{ptm}{m}{n}
\rput(5.571455,1.065){$l$}
\usefont{T1}{ptm}{m}{n}
\rput(2.861455,1.205){$l+1$}
\usefont{T1}{ptm}{m}{n}
\rput(4.361455,0.505){$a$}
\usefont{T1}{ptm}{m}{n}
\rput(4.6,0.225){$b$}
\usefont{T1}{ptm}{m}{n}
\rput(3.900000,0.225){$c$}
\psbezier[linewidth=0.016,linecolor=blue,arrowsize=0.07cm 8.0,arrowlength=2.0,arrowinset=0.6]{->}(5.3,0.68)(5.34,2.28)(3.24,2.46)(3.22,0.74)
\usefont{T1}{ptm}{m}{n}
\rput(8.081455,-1.595){$l-1$}
\usefont{T1}{ptm}{m}{n}
\rput(0.4214551,-1.495){$l+2$}
\psbezier[linewidth=0.016,linecolor=blue,arrowsize=0.07cm 8.0,arrowlength=2.0,arrowinset=0.6]{->}(3.18,0.76)(1.1945455,2.02)(0.0,0.31333333)(0.94,-1.34)
\psbezier[linewidth=0.016,linecolor=blue,arrowsize=0.07cm 8.0,arrowlength=2.0,arrowinset=0.6]{->}(7.54,-1.48)(8.56,-0.35757226)(7.6512055,1.76)(5.36,0.68)
\usefont{T1}{ptm}{m}{n}
\rput(8.221455,0.505){$\alpha$}
\usefont{T1}{ptm}{m}{n}
\rput(4.131455,2.185){$a^*$}
\usefont{T1}{ptm}{m}{n}
\rput(0.89145505,1.105){$\beta$}
\usefont{T1}{ptm}{m}{n}
\rput(4.181455,-0.455){$m$}
\usefont{T1}{ptm}{m}{n}
\rput(4.571455,-1.815){$\tau$}
\end{pspicture} 
\end{pdfpic}
}
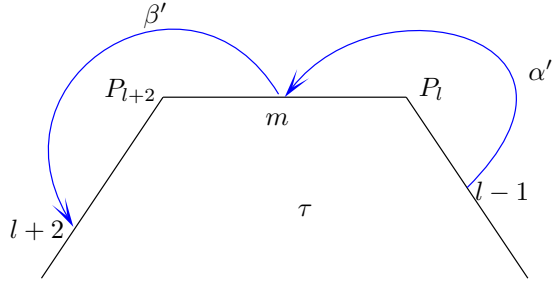
{
\begin{pdfpic}
\begin{pspicture}(0,-2.21)(8.36,2.21)
\psline[linewidth=0.016cm](2.34,0.21)(5.54,0.21)
\psline[linewidth=0.016cm](2.34,0.21)(0.74,-2.19)
\psline[linewidth=0.016cm](5.54,0.21)(7.14,-2.19)
\usefont{T1}{ptm}{m}{n}
\rput(5.871455,0.315){$P_{l}$}
\usefont{T1}{ptm}{m}{n}
\rput(1.911455,0.295){$P_{l+2}$}
\usefont{T1}{ptm}{m}{n}
\rput(6.801455,-1.045){$l-1$}
\usefont{T1}{ptm}{m}{n}
\rput(3.841455,-0.085){$m$}
\usefont{T1}{ptm}{m}{n}
\rput(0.6814551,-1.545){$l+2$}
\usefont{T1}{ptm}{m}{n}
\rput(4.211455,-1.285){$\tau$}
\psbezier[linewidth=0.016,linecolor=blue,arrowsize=0.07cm 8.0,arrowlength=2.0,arrowinset=0.6]{->}(6.34,-0.99)(8.34,1.01)(5.1844444,1.81)(3.94,0.21)
\psbezier[linewidth=0.016,linecolor=blue,arrowsize=0.07cm 8.0,arrowlength=2.0,arrowinset=0.6]{->}(3.86,0.25)(2.7,2.19)(0.0,0.51)(1.16,-1.51)
\usefont{T1}{ptm}{m}{n}
\rput(7.3114552,0.575){$\alpha'$}
\usefont{T1}{ptm}{m}{n}
\rput(2.241455,1.295){$\beta'$}
\end{pspicture} 
\end{pdfpic}
}
\end{center}
(if there is no corner triangle, we are either in the case of Lemma \ref{lemdig} either in the case of Lemma \ref{lemtart}). By induction hypothesis, the expected results hold for $\tau$. As the $\theta$-length depends on the size of the polygon, we will denote by $\ell^{\theta,\tau}$ (resp. $\ell^\theta$) the $\theta$-length in $\tau$ (resp. $\sigma$). In the same way, as the inclusion $R' \subset \cR'$ depends on the triangulation, we will call $u_\tau$ and $v_\tau$ the generators of $\cR'$ when we consider this inclusion for $\tau$.

Then there is a non-unital monomorphism
\begin{align*}
\xi: KQ'_{\tau} & \hookrightarrow KQ'_{\sigma}\\
\alpha' &\mapsto \alpha b,\\
\beta'  &\mapsto c\beta,\\
\gamma  &\mapsto \gamma \text{ for } \gamma \in (Q'_{\tau})_1\smallsetminus\{\alpha', \beta'\}.
\end{align*}

We have $W'_\sigma = \xi(W'_\tau) + abc - aa^*$, so for any $\gamma \in Q'_\tau$ which is not external, $\partial_\gamma W'_\sigma = \xi(\partial_\gamma W'_\tau)$. Thus $\xi$ induces a morphism $\bar \xi: \Gamma_\tau \rightarrow \Gamma_\sigma$. Notice that using the relations, any path of $Q'_\sigma$ is equivalent to an element which does not contain $a^*$, $ca$ or $ab$. It is then easy to see that:

\begin{claim} Any path of $Q'_\sigma$ is equivalent to a path of one of the following forms where $\omega$ is a path of $Q'_\tau$:
\NumTabs{5}
\noindent\begin{inparaenum}
 \item $\xi(\omega)$;
 \tab\item $\xi(\omega) \alpha$;
 \tab\item $\xi(\omega) c$;
 \tab\item $b \xi(\omega)$;
 \tab\item $\beta \xi(\omega)$;
 \tab\item $b \xi(\omega) \alpha$;
 \tab\item $b \xi(\omega) c$;
 \tab\item $a$;
 \tab\item $\beta \xi(\omega) \alpha$;
 \tab\item $\beta \xi(\omega) c$.
\end{inparaenum} \label{uniq3}
\end{claim}

Let us prove that $\bar \xi$ is in fact injective. Let $\pi_\tau$ be the canonical projection $K Q'_\tau \rightarrow \Gamma_\tau$. Thanks to $(iv)_\tau$,  we get: 
\begin{align*}\ker \bar \xi &= \pi_\tau(\xi^{-1}(\langle \partial_b W'_\sigma, \partial_c W'_\sigma\rangle + \xi(\mathcal{J'}(W'_\tau)))) \\ &= \pi_\tau(\xi^{-1}(\langle\partial_b W'_\sigma, \partial_c W'_\sigma\rangle) + \xi^{-1}(\xi(\mathcal{J'}(W'_\tau)))) \\ &= \pi_\tau(\xi^{-1}(\langle\partial_b W'_\sigma, \partial_c W'_\sigma\rangle)) + \pi_\tau(\mathcal{J'}(W'_\tau)) \\&= \pi_\tau\left(\xi^{-1}\left(\left\langle ca - \omega_{m,l-1} \alpha, ab - \beta \omega_{l+2,m}\right\rangle\right)\right)\end{align*}
where $\omega_{m,l-1} = \psi_\tau^{m,l-1}\left(v_\tau^{\ell^{\theta,\tau}_{m,l-1}}\right)$ and $\omega_{l+2,m}= \psi_\tau^{l+2,m}\left(v_\tau^{\ell^{\theta,\tau}_{l+2,m}}\right)$.

Up to equivalence, a path of $Q'_\sigma$ can always be supposed not to contain $a^*$. Moreover, $e_{m} \im(\xi) e_{l} = 0 = e_{l+1} \im(\xi) e_{m}$, so
$$ \xi^{-1}\left(\left\langle ca - \omega_{m,l-1} \alpha, ab - \beta \omega_{l+2,m}\right\rangle\right) = \xi^{-1}\left(\left\langle cab - \omega_{m,l-1} \alpha b, cab - c \beta \omega_{l+2,m}\right\rangle\right)$$
\begin{align*}
\text{and} \quad \ker \bar \xi &=  \pi_\tau\left(\xi^{-1}\left(\left\langle cab - \omega_{m,l-1} \alpha b, cab - c \beta \omega_{l+2,m}\right\rangle\right)\right) \\
      &= \pi_\tau\left(\xi^{-1}\left(\left\langle c\beta\omega_{l+2,m} - \omega_{m,l-1} \alpha b, cab - c\beta \omega_{l+2,m}\right\rangle\right)\right) \\
      &= \pi_\tau\left(\left\langle \beta' \omega_{l+2,m} - \omega_{m,l-1} \alpha'\right\rangle + \xi^{-1}\left(\left\langle cab - c\beta\omega_{l+2,m}\right\rangle\right)\right) \\
      &= \pi_\tau\left(\xi^{-1}\left(\left\langle cab - c\beta \omega_{l+2,m}\right\rangle\right)\right) = \pi_\tau(0) = 0. \\
\end{align*}

Let us prove that the following map is an isomorphism:
\begin{align*}
 \phi_\sigma: Z(\Gamma_\tau) \cong R' & \rightarrow Z(\Gamma_\sigma) \\
 X & \mapsto \bar \xi(\phi_\tau(X)) + E_l + E_{l+1} \\
 Y & \mapsto \bar \xi(\phi_\tau(Y)) + C_l + C_{l+1}
\end{align*}

First of all, $\phi_\sigma(X), \phi_\sigma(Y) \in Z(\Gamma_\sigma)$. Indeed, for each arrow $\gamma$ in $Q_\tau$, we have, by using induction hypothesis,
$$\gamma \phi_\sigma(X) = \gamma \bar \xi(\phi_\tau(X)) = \bar \xi(\gamma \phi_\tau(X)) = \bar \xi(\phi_\tau(X) \gamma) = \bar \xi(\phi_\tau(X)) \gamma = \phi_\sigma(X) \gamma $$
and the same for $Y$. For arrows which are not in $Q_\tau$ we have
\begin{align*}
 a \phi_\sigma(X) &= a E_l = a a^* \beta \overbrace{\dots}^{\text{ext.}} \alpha = C_{l+1} \beta \overbrace{\dots}^{\text{ext.}} \alpha \\ &= \beta C_{l+2} \overbrace{\dots}^{\text{ext.}} \alpha = \dots = \beta \overbrace{\dots}^{\text{ext.}} C_l = E_{l+1} a = \phi_\sigma(X) a
\end{align*}
where ext. denote products of external arrows and, thanks to $(i)_\tau$,
$$b \phi_\sigma(X) = b c \beta \overbrace{\dots}^{\text{ext.}} \alpha b = E_l b  = \phi_\sigma(X) b;$$
$$c \phi_\sigma(X) = c E_{l+1} = c \beta \overbrace{\dots}^{\text{ext.}} \alpha b c = \phi_\sigma(X) c;$$
$$\beta \phi_\sigma(X) = \beta E_{l+2} = E_{l+1} \beta = \phi_\sigma(X) \beta;$$ 
$$\alpha \phi_\sigma(X) = \alpha E_l = E_{l-1} \alpha = \phi_\sigma(X) \alpha.$$
For $\phi_\sigma(Y)$, using $(ii)_\tau$, all the computations are immediate. 

Notice now that 
\begin{align*}
 C_l^2 &= bc\left(ab\right)\left(ca\right) = b c \beta \bar \xi\left(\psi_\tau^{l+2,m}\left(v_\tau^{\ell^{\theta,\tau}_{l+2,m}}\right)\right) \bar \xi\left(\psi_\tau^{m,l-1}\left(v_\tau^{\ell^{\theta,\tau}_{m,l-1}}\right)\right) \alpha \\ &= b c \beta \bar \xi\left(\psi_\tau^{l+2,l-1}\left(v_\tau^{\ell^{\theta,\tau}_{l+2,l-1}+2n}\right)\right) \alpha        
       = b c \beta \bar \xi\left(\psi_\tau^{l+2,l-1}\left(Y u_\tau^{\ell^{\theta,\tau}_{l+2,l-1}}\right)\right) \alpha \\ &= b c \beta C_{l+2} \bar \xi\left(\psi_\tau^{l+2,l-1}\left(u_\tau^{\ell^{\theta,\tau}_{l+2,l-1}}\right)\right) \alpha \\ &= C_l b c \beta  \bar \xi\left(\psi_\tau^{l+2,l-1}\left(u_\tau^{\ell^{\theta,\tau}_{l+2,l-1}}\right)\right) \alpha = C_l E_l
\end{align*}
and $C_{l+1}^2 = C_{l+1} E_{l+1} $ by the same method. Therefore
$$\phi_\sigma(YX) = \bar \xi(\phi_\tau(Y)) \bar \xi(\phi_\tau(X)) + C_l E_l + C_{l+1} E_{l+1} = \bar \xi(\phi_\tau(Y^2)) + C_l^2 + C_{l+1}^2 = \phi_\sigma(Y^2)$$ 
so $\phi_\sigma$ is a morphism. As $\bar \xi$ and $\phi_\tau$ are injective, $\phi_\sigma$ is also injective. The last thing to show is that $\phi_\sigma$  is surjective.

 Using Claim \ref{uniq3} and that, as it commutes with idempotents, every element $z \in Z(\Gamma_\sigma)$ is a linear combination of cycles, we can write
$$z = \bar \xi(z') + b \bar \xi(z'') \alpha +  \beta \bar \xi(z''') c.$$
Then, as $z$ is in the center, for any $x \in \Gamma_\tau$, 
$$\bar \xi(x z') = \bar \xi(x) \bar \xi(z') = \bar \xi(x) z = z \bar \xi(x) = \bar \xi(z') \bar \xi(x) = \bar \xi(z' x)$$
and, as $\bar \xi$ is injective, $xz' = z'x$ and $z' \in Z(\Gamma_\tau)$. Therefore, up to subtracting $\phi_\sigma(z')$, we can suppose that $z = b \bar \xi(z'') \alpha + \beta \bar \xi(z''') c$. Hence, we have
$$0 = z \alpha b = \alpha z b = \bar \xi (\alpha' z'' \alpha') $$
and, as $\bar \xi$ is injective, $\alpha' z'' \alpha' = 0$. Finally, $e_m z'' e_{l-1} = 0$ thanks to $(v)_\tau$. In the same way $e_{l+2} z''' e_m = 0$ so $z = 0$. Therefore $\phi_\sigma$ is surjective.

We will prove the existence of $\psi_\sigma$ as a family of morphisms of $R'$-modules. Then, these morphisms are automatically graded by looking at homogeneous generators. If $i,j \notin \{l, l+1\}$, $\bar \xi$ induces an isomorphism of $R'$-modules from $e_i \Gamma_\tau e_j$ to $e_i \Gamma_\sigma e_j$. It proves the existence of $\psi_\sigma^{i,j}$ in this case.

Let us recall that $l = (P_l, P_{l+1})$, $l+1 = (P_{l+1}, P_{l+2})$ and $m = (P_l, P_{l+2})$. Thus, $l+1 \vdash l$, $l \notvdash l+1$, and for any $i \in Q'_{\sigma,0} \smallsetminus \{l, l+1\}$ which is not incident to the puncture, we have $i_1, i_2 \neq l+1$, so
\begin{align*} i \vdash l &\Leftrightarrow i_2 \in \ooi{l+1}{i_1} \text{ or } l \in \ooi{l+1}{i_1} \\ &\Leftrightarrow i_2 \in \ooi{l}{i_1} \text{ or } l-1 \in \ooi{l}{i_1} \Leftrightarrow i \vdash l-1\end{align*}
\begin{align*}l \vdash i &\Leftrightarrow l+1 \in \ooi{i_2}{l} \text{ or } i_1 \in \ooi{i_2}{l} \\ &\Leftrightarrow l+2 \in \ooi{i_2}{l} \text{ or } i_1 \in \ooi{i_2}{l} \Leftrightarrow m \vdash i\end{align*}
\begin{align*}i \vdash l+1 &\Leftrightarrow i_2 \in \ooi{l+2}{i_1} \text{ or } l+1 \in \ooi{l+2}{i_1} \\ &\Leftrightarrow i_2 \in \ooi{l+2}{i_1} \text{ or } l \in \ooi{l+2}{i_1} \Leftrightarrow i \vdash m\end{align*}
\begin{align*}l+1 \vdash i &\Leftrightarrow l+2 \in \ooi{i_2}{l+1} \text{ or } i_1 \in \ooi{i_2}{l+1} \\ &\Leftrightarrow l+3 \in \ooi{i_2}{l+2} \text{ or } i_1 \in \ooi{i_2}{l+2} \Leftrightarrow l+2 \vdash i.\end{align*}

Let $i \notin \{l, l+1\}$. There are isomorphisms of $R'$-modules
$$\begin{aligned}
 e_i \Gamma_\tau e_{l-1} &\rightarrow e_i \Gamma_\sigma e_l \\
  \omega & \mapsto \bar \xi(\omega) \alpha
\end{aligned} \quad ; \quad
\begin{aligned}
 e_{l+2} \Gamma_\tau e_i &\rightarrow e_{l+1} \Gamma_\sigma e_i \\
  \omega & \mapsto \beta \bar \xi(\omega) 
\end{aligned} \quad ;$$
$$\begin{aligned}
 e_m \Gamma_\tau e_i  &\rightarrow e_l \Gamma_\sigma e_i \\
  \omega & \mapsto b \bar \xi(\omega)
\end{aligned} \quad ; \quad
\begin{aligned}
 e_i \Gamma_\tau e_m  &\rightarrow e_i \Gamma_\sigma e_{l+1} \\
  \omega & \mapsto \bar \xi(\omega) c
\end{aligned}.$$
The injectivity comes from $(v)_\tau$. For example, if $\bar \xi(\omega) \alpha = 0$ then $\bar \xi(\omega) \alpha b = 0$ and therefore $\bar \xi (\omega \alpha') = 0$ so $\omega \alpha' = 0$ and finally $\omega = 0$. For the surjectivity, it is enough to use Claim \ref{uniq3}. In the same way, there is an isomorphism of $R'$-modules
\begin{align*}
 e_m \Gamma_\tau e_m  &\rightarrow e_l \Gamma_\sigma e_{l+1} \\
  \omega & \mapsto b \bar \xi(\omega) c.
\end{align*}

Thus we get the expected $R'$-module structure for $e_i \Gamma_\sigma e_j$ except when $i = j \in \{l, l+1\}$ or $i = l+1$ and $j = l$. 

Suppose that $i = j = l$. The elements of $e_l \Gamma_\sigma e_l$ are of the form $\lambda e_l + b \omega \alpha$ for $\lambda \in K$ and $\omega \in e_m \Gamma_\sigma e_{l-1}$. We already know that $$e_m \Gamma_\sigma e_{l-1} \cong u^{\ell^{\theta}_{m,l-1}} R' + v^{\ell^{\theta}_{m,l-1}} R'$$ and we get the following isomorphism of $R'$-modules:
\begin{align*}
 R' \cong_K K \oplus u^3(u^{\ell^{\theta}_{m,l-1}} R' + v^{\ell^{\theta}_{m,l-1}} R') & \rightarrow e_l \Gamma_\sigma e_l \\
 (\lambda, u^3 p) & \mapsto \lambda e_l + b \psi_\sigma^{m, l-1}(p) \alpha
\end{align*}
(the injectivity comes from $(v)_\tau$ and the injectivity of $\bar \xi$).

In the same way, if $i = j = l+1$, there is an isomorphism of $R'$-modules
\begin{align*}
 R' \cong_K K \oplus u^3(u^{\ell^{\theta}_{l+2,m}} R' + v^{\ell^{\theta}_{l+2,m}} R') & \rightarrow e_{l+1} \Gamma_\sigma e_{l+1} \\
 (\lambda, u^3 p) & \mapsto \lambda e_{l+1} + \beta \psi_\sigma^{l+2, m}(p) c.
\end{align*}

Finally, suppose that $i = l+1$ and $j = l$. The elements of $e_{l+1} \Gamma_\sigma e_l$ are of the form $\lambda a + \beta \omega \alpha$ and there is an isomorphism of $R'$-modules
\begin{align*}
 u^{\ell^{\theta}_{l+1,l}} R' + v^{\ell^{\theta}_{l+1,l}} R'  \cong_K K v^{\ell^{\theta}_{l+1,l}} \oplus u^4 u^{\ell^{\theta}_{l+2,l-1}} R' & \rightarrow e_{l+1} \Gamma_\sigma e_l \\
 (\lambda v^{\ell^{\theta}_{l+1,l}}, u^4 p) & \mapsto \lambda a + \beta \psi_\sigma^{l+2, l-1}(p) \alpha.
\end{align*}
Indeed, the only non-immediate thing to check is $a E_l = a C_l$. By induction hypothesis (in particular $(v)_\tau$), 
\begin{align*} a E_l &= abc \beta \bar \xi\left(\psi_\tau^{l+2, l-1}\left(u_\tau^{\ell^{\theta,\tau}_{l+2,l-1}}\right)\right) \alpha \\ &= \beta \bar \xi\left(\psi_\tau^{l+2, l-1}\left(Y u_\tau^{\ell^{\theta,\tau}_{l+2,l-1}}\right)\right) \alpha \\ &= \beta \bar \xi\left(\psi_\tau^{l+2, l-1}\left(Y v_\tau^{\ell^{\theta,\tau}_{l+2,l-1}}\right)\right) \alpha \\ &= \beta \bar \xi\left(\psi_\tau^{l+2, m}\left(v_\tau^{\ell^{\theta,\tau}_{l+2,m}}\right)\right) \bar \xi\left(\psi_\tau^{m, l-1}\left(v_\tau^{\ell^{\theta,\tau}_{m,l-1}}\right)\right) \alpha \\ &= ab \bar \xi\left(\psi_\tau^{m, l-1}\left(v_\tau^{\ell^{\theta,\tau}_{m,l-1}}\right)\right) \alpha = a C_l.\end{align*}

For the multiplicativity of the $\psi_\sigma$, let us start by noticing, thanks to the beginning of the proof of Lemma \ref{lem:subadd} that for any vertices $i$, $j$ and $k$ of $Q'_\tau$, we have the following identity:
$$\frac{\ell^\theta_{i,j} + \ell^\theta_{j,k} - \ell^\theta_{i,k}}{n} = \frac{\ell^{\theta, \tau}_{i,j} + \ell^{\theta, \tau}_{j,k} - \ell^{\theta, \tau}_{i,k}}{n-1}.$$

It permits to prove that $\psi_\sigma^{i,j}(w) \psi_\sigma^{j,k}(w') = \psi_\sigma^{i,k}(ww')$ for any $(w,w') \in A_{i,j} \times A_{j,k}$. Indeed this is enough to prove that if $w$ and $w'$ are generators as $R'$-modules. Suppose for example that $w = u^{\ell^\theta_{i,j}}$ and $w' = u^{\ell^\theta_{j,k}}$. Then we have, using the induction hypothesis,

\begin{align*}
 \psi_\sigma^{i,j}\left(u^{\ell^\theta_{i,j}}\right) \psi_\sigma^{j,k}\left(u^{\ell^\theta_{j,k}}\right) &= \bar \xi\left(\psi_\tau^{i,j}\left(u_\tau^{\ell^{\theta,\tau}_{i,j}}\right) \psi_\tau^{j,k}\left(u_\tau^{\ell^{\theta,\tau}_{j,k}}\right)\right) \\ &= \bar \xi\left(\psi_\tau^{i,k}\left(u_\tau^{\ell^{\theta,\tau}_{i,j}+\ell^{\theta,\tau}_{j,k}}\right)\right) \\ &=  \bar \xi\left(\phi_\tau\left(X^{\left(\ell^{\theta,\tau}_{i,j}+\ell^{\theta,\tau}_{j,k}-\ell^{\theta,\tau}_{i,k}\right)/2\left(n-1\right)}\right)\psi_\tau^{i,k}\left(u_\tau^{\ell^{\theta,\tau}_{i,k}}\right)\right) \\&= \phi_\sigma\left(X^{\left(\ell^{\theta,\tau}_{i,j}+\ell^{\theta,\tau}_{j,k}-\ell^{\theta,\tau}_{i,k}\right)/2\left(n-1\right)}\right)\psi_\sigma^{i,k}\left(u_\tau^{\ell^\theta_{i,k}}\right) \\&= \phi_\sigma\left(X^{\left(\ell^\theta_{i,j}+\ell^\theta_{j,k}-\ell^\theta_{i,k}\right)/2n}\right)\psi_\sigma^{i,k}\left(u^{\ell^\theta_{i,k}}\right) \\&= \psi_\sigma^{i,k}\left(u^{\ell^\theta_{i,j}+\ell^\theta_{j,k}}\right).
\end{align*}

The multiplicativity for paths starting or ending at $l$ or $l+1$ can be deduced easily from that. For example, if $i$ and $k$ are vertices of $\tau$.
\begin{align*}
 \psi_\sigma^{i,l}\left(u^{\ell^\theta_{i,l}}\right) \psi_\sigma^{l,k}\left(u^{\ell^\theta_{l,k}}\right) &= \psi_\sigma^{i,l-1}\left(u^{\ell^\theta_{i,l-1}}\right) \alpha b \psi_\sigma^{m,k}\left(u^{\ell^\theta_{m,k}}\right) \\ &=
 \psi_\sigma^{i,l-1}\left(u^{\ell^\theta_{i,l-1}}\right) \psi_\sigma^{l-1,m}(u^3) \psi_\sigma^{m,k}\left(u^{\ell^\theta_{m,k}}\right) \\ &=
 \psi_\sigma^{i,k}\left(u^{\ell^\theta_{i,l-1}+3+\ell^\theta_{m,k}}\right)  =
 \psi_\sigma^{i,k}\left(u^{\ell^\theta_{i,l}+\ell^\theta_{l,k}}\right).
\end{align*}

The last thing to check are the five additional conditions. Points $(i)_\sigma$ to $(iv)_\sigma$ are easy to check and $(v)_\sigma$ is a consequence of them. 
\end{proof}

\begin{theorem} \label{cor:lambdamod}
 There is an isomorphism of $R$-orders (and $R'$-algebras):
 $$e_F \Gamma_\sigma e_F = \left[u^{2 d(i,j)} R' + v^{2 d(i,j)} R'^{\delt{j = i-1}}\right]_{i,j \in \cci{1}{n}} \cong \Lambda$$
 where the entries of $e_F \Gamma_\sigma e_F$ are $R'$-submodules of $\cR'$ and $\Lambda$ is defined at \eqref{the order}.

 For each edge $a$ of $\sigma$, the $e_F \Gamma_\sigma e_F$ module
 $$M_a := e_F \Gamma_\sigma e_a = \tr{\begin{bmatrix} A_{1,a}, A_{2,a}, \cdots, A_{n,a} \end{bmatrix} }$$
 is, as a $\Lambda$-module, isomorphic to
 \spovermat{a_1}
 $$\tr{\begin{bmatrix} \bovermat{a_1}{(Y) \cdots (Y)} & \bovermat{n-a_1}{(Y^2) \cdots (Y^2)} \end{bmatrix} } \quad \text{if } a = (P_{a_1},*);$$
 \spovermat{a_1}
 $$\tr{ \begin{bmatrix} \bovermat{a_1}{(X-Y) \cdots (X-Y)} & \bovermat{n-a_1}{(X^2-Y^2) \cdots (X^2-Y^2)} \end{bmatrix} } \quad \text{if } a = (P_{a_1},\bowtie );$$
 \spovermat{a_1}
 $$\tr{\begin{bmatrix} \bovermat{a_1}{R' \cdots R'} & \bovermat{a_2-a_1}{(X,Y) \cdots (X, Y)} & \bovermat{n-a_2}{(X) \cdots (X)} \end{bmatrix}} \quad \text{if } a_1 < a_2;$$
 \spovermat{a_1}
 $$\tr{\begin{bmatrix} \bovermat{a_2}{(X,Y) \cdots (X,Y)} & \bovermat{a_1-a_2}{(X) \cdots (X)} & \bovermat{n-a_1}{(X^2,Y^2) \cdots (X^2,Y^2)} \end{bmatrix}} \quad \text{if } a_1 > a_2.$$
\end{theorem}

\begin{proof}
 If we order the sides of the polygon in the order $(P_1, P_2)$, $(P_2,P_3)$, \dots, $(P_{n-1}, P_n)$, $(P_n, P_1)$, the first equality is a direct application of Proposition \ref{inductD}. Notice that for sides $i$ and $j$ of the polygon, we can rewrite $A_{i,j}$ in the following way:
 $$A_{i,j} = u^{2 d(i,j)} R' + v^{2 d(i,j)} R'^{(\delt{j = i-1})}.$$
 We conjugate by the diagonal matrix with diagonal entries $u^{2 d(1,i)}$ for $i \in \cci{1}{n}$ and the matrix we obtain has entries
 \begin{align*}
  & u^{2 d(1,i)} \left(u^{2 d(i,j)} R' + v^{2 d(i,j)} R'^{(\delt{j = i-1})}\right) u^{-2 d(1,j)} \\ =\,& u^{2 d(i,j)-1 + 2 d(1,i) - 2 d(1,j)} \left( u R' + v R'^{(\delt{j = i-1})}\right) \\ =\,& u^{2 n \delt{i \in \ooi{j}{1}} -1} \left( u R' + v R'^{(\delt{j = i-1})}\right) = X^{(\delt{i > j})} R' + X^{-(\delt{i \leq j})} Y R'^{(\delt{j = i-1})}
 \end{align*}
 for $i, j \in \cci{1}{n}$. It is $\Lambda$. 

 We obtain the structure of $M_a$, up to some degree shift, by multiplying on the left by the same diagonal matrix.
\end{proof}

\begin{remark}
 Notice that in Theorem \ref{cor:lambdamod}, the module $M_a$ depends only on the edge $a$ and not on the triangulation $\sigma$. 
\end{remark}

\subsection{Counterexample with more than one puncture} \label{cex}

In this subsection, we show that we can not expect to generalize these results for polygons with more than one puncture. 

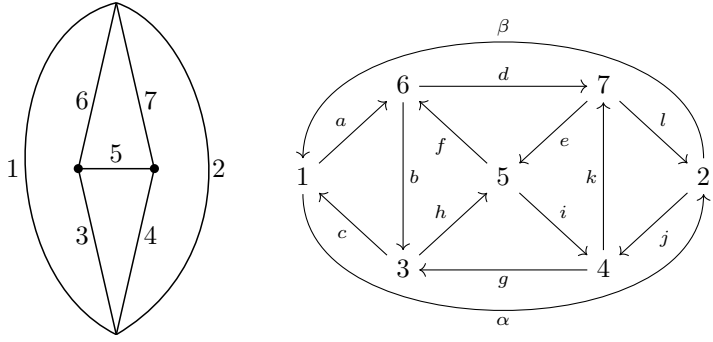
\begin{figure}
\begin{minipage}{4cm}
{
\begin{pdfpic}
\begin{pspicture}(0,-2.2)(3.61,2.28)
\psdots[dotsize=0.12](1.1,0)
\psdots[dotsize=0.12](2.1,0)
\psbezier[linewidth=0.02](1.6,2.2)(2.8,1.8)(3.6,-1.0)(1.6,-2.2)
\psbezier[linewidth=0.02](1.6,2.2)(0.0,1.8)(0.0,-1.4)(1.6,-2.2)

\psline[linewidth=0.02cm](1.1,0)(1.6,-2.2)
\psline[linewidth=0.02cm](2.1,0)(1.6,-2.2)
\psline[linewidth=0.02cm](1.1,0)(1.6,2.2)
\psline[linewidth=0.02cm](2.1,0)(1.6,2.2)
\psline[linewidth=0.02cm](1.1,0)(2.1,0)
\usefont{T1}{ptm}{m}{n}
\rput(0.24,0){$1$}
\usefont{T1}{ptm}{m}{n}
\rput(2.95,0){$2$}
\usefont{T1}{ptm}{m}{n}
\rput(1.15,-0.895){$3$}
\usefont{T1}{ptm}{m}{n}
\rput(2.05,-0.895){$4$}
\usefont{T1}{ptm}{m}{n}
\rput(1.6,0.2){$5$}
\usefont{T1}{ptm}{m}{n}
\rput(1.15,0.895){$6$}
\usefont{T1}{ptm}{m}{n}
\rput(2.05,0.895){$7$}

\end{pspicture} 
\end{pdfpic}
}
\end{minipage}
\begin{minipage}{6.5cm}\begin{tikzcd}[ampersand replacement=\&]
\& 6 \arrow{dd}{b} \arrow{rr}{d} \& \& 7 \arrow{dl}{e} \arrow{dr}{l} \& \\
1 \arrow[bend right=90,swap]{rrrr}{\alpha} \arrow{ur}{a} \& \& 5 \arrow{ul}{f} \arrow{dr}{i} \& \& 2 \arrow{dl}{j} \arrow[bend right=90,swap]{llll}{\beta}  \\
\& 3 \arrow{ul}{c} \arrow{ur}{h} \& \& 4 \arrow{ll}{g} \arrow{uu}{k} \& \\
\end{tikzcd}\end{minipage}
\caption{twice-punctured digon}
\label{tdigon}
\end{figure}

We take the triangulation $\sigma$ of a twice-punctured digon of Figure \ref{tdigon}. It induces the quiver $Q_\sigma$ on the right and using the same definition than in Section \ref{Ice QP}, we get that the natural analogue of $\Gamma_\sigma$ is the path algebra of the quiver modulo all obvious commutativity relations. We still call it $\Gamma_\sigma$. Suppose that $\Gamma_\sigma$ is a $K[U]$-order, for $U$ in the center of $\Gamma_\sigma$. We can write $e_1 U = P(\alpha \beta) + a \omega$ where $P$ is a polynomial and $\omega \in e_6 \Gamma_\sigma e_1$. 
 As, for $\ell > 0$, $c (\alpha \beta)^{\ell}$ is clearly not divisible by $c$ on the right, and as $c U = U c$, we get that $P$ is a constant polynomial. So, if we denote $\pi: e_1 \Gamma_\sigma e_1 \rightarrow K[\alpha \beta]$ the canonical projection, we get $\pi(U) = P(0) \in K$. So $K[\alpha \beta]$ is not a finitely generated module over $K[\pi(U)] = K$ and therefore $e_1 \Gamma_\sigma e_1$ is not a finitely generated module over $K[U]$. It is a contradiction.

This counterexample is easy to generalize to any polygon with at least two punctures.

\section{Cohen-Macaulay Modules over $\Lambda$}
\label{sec:CM}
The aim of this section is to study the representation theory of $\Lambda$ and its connection to tagged triangulations of the punctured polygon $P^*$ and the cluster category of type $D_n$.
In particular, we classify all Cohen-Macaulay $\Lambda$-modules and construct a bijection between the set of the isomorphism classes of all indecomposable Cohen-Macaulay $\Lambda$-modules and the set of all sides and tagged arcs of $P^*$. We then  show that the stable category $\uCM \Lambda$ of Cohen-Macaulay $\Lambda$-modules is $2$-Calabi-Yau and $\uCM \Lambda$ is triangle-equivalent with the cluster category of type $D_n$. To summarize, we will prove that $\CM \Lambda$ admits the Auslander-Reiten quiver of Figures \ref{CMLe}, \ref{CMLo}.

\newcommand{\smp}{{+}}\newcommand{\smm}{{-}}
\begin{figure}
{\tiny \[
\begin{tikzcd}[column sep=0ex, row sep=2ex]
\phantom{X}\arrow[dashed,-]{dd} && & & & & &&& &\phantom{X}\\
 & (m\smp2,\bowtie) \drar & & (m\smp3,*) \arrow[dotted,-]{rrrr} && && (m,\bowtie) \drar & & (m\smp1,*) \drar &  \\
(m\smp2,m\smp1) \urar \drar\arrow[dashed,-]{dd} \arrow{r} & (m\smp2,*) \arrow{r} &(m\smp3,m\smp2)\drar \urar \arrow{r} & (m\smp3,\bowtie) \arrow[dotted,-]{rrrr} && && (m,*) \arrow{r} &(m\smp1,m) \urar \drar \arrow{r} & (m\smp1,\bowtie) \arrow{r} & (m\smp2,m\smp1) \arrow[dashed,-]{uu} \arrow[dashed,-]{dd} \\
&(m\smp3,m\smp1)\urar\drar & &(m\smp4,m\smp2) \arrow[dotted,-]{rrrr} &&&&(m\smp1, m\smm1) \urar\drar & &(m\smp2,m) \urar \drar \\
(m\smp3,m)\urar\drar \arrow[dashed,-]{ddddd}&&(m\smp4,m\smp1)\urar\drar\arrow[dotted,-]{ddddd}\arrow[dotted,-]{rrrrrr}&&&& & &(m\smp2,m\smm1)\drar\urar \arrow[dotted,-]{ddddd}& & (m\smp3,m) \arrow[dashed,-]{ddddd}\\ 
&(m\smp4,m)\urar\arrow[dotted,-]{ddddd}& &(m\smp5,m\smp1) \arrow[dotted,-]{ddddd}\arrow[dotted,-]{rrrr} &&&&(m\smp2,m\smm2)  \urar\arrow[dotted,-]{ddddd} & &(m\smp3,m\smm1)\arrow[dotted,-]{ddddd} \urar \\
& \quad \\
 & \quad \\
 & \quad \\
(2m,3)\drar \arrow[dashed,-]{dd} & &
(1,4)\drar\arrow[dotted,-]{rrrrrr} &&& && &(2m\smm1,2)\drar  & & (2m,3) \arrow[dashed,-]{dd}\\
&(1,3)\drar\urar &&(2,4)\arrow[dotted,-]{rrrr}&& &&(2m\smm1,1)\urar\drar&&(2m,2) \urar\drar \\
(1,2)\urar & &
(2,3)\urar \arrow[dotted,-]{rrrrrr} &&& && &(2m,1) \urar & & (1,2)\\
\phantom{X}\arrow[dashed,-]{u}& & & & & &&&&& \phantom{X}\arrow[dashed,-]{u}
\end{tikzcd}
\]}
\caption{$\CM \Lambda$ for $n = 2m$}
\label{CMLe}
\end{figure}
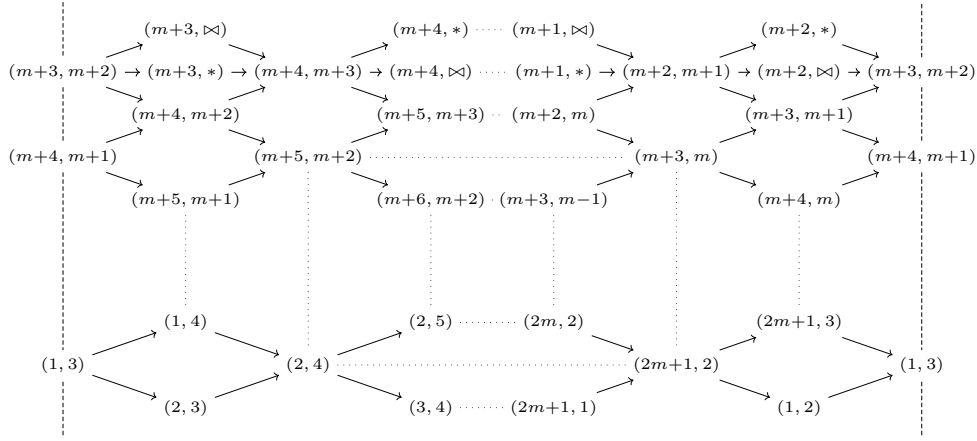
\begin{figure}
{\tiny \[
\begin{tikzcd}[column sep=0ex, row sep=2ex]
\phantom{X}\arrow[dashed,-]{dd} && & & & & &&& &\phantom{X}\\
 & (m\smp3,\bowtie) \drar & & (m\smp4,*) \arrow[dotted,-]{rrrr} && && (m\smp1,\bowtie) \drar & & (m\smp2,*) \drar &  \\
(m\smp3,m\smp2) \urar \drar\arrow[dashed,-]{dd} \arrow{r} & (m\smp3,*) \arrow{r} &(m\smp4,m\smp3)\drar \urar \arrow{r} & (m\smp4,\bowtie) \arrow[dotted,-]{rrrr} && && (m\smp1,*) \arrow{r} &(m\smp2,m\smp1) \urar \drar \arrow{r} & (m\smp2,\bowtie) \arrow{r} & (m\smp3,m\smp2) \arrow[dashed,-]{uu} \arrow[dashed,-]{dd} \\
&(m\smp4,m\smp2)\urar\drar & &(m\smp5,m\smp3) \arrow[dotted,-]{rrrr} &&&&(m\smp2, m) \urar\drar & &(m\smp3,m\smp1) \urar \drar \\
(m\smp4,m\smp1)\urar\drar \arrow[dashed,-]{dddddd}&&(m\smp5,m\smp2)\urar\drar\arrow[dotted,-]{dddddd}\arrow[dotted,-]{rrrrrr}&&&& & &(m\smp3,m)\drar\urar \arrow[dotted,-]{dddddd}& & (m\smp4,m\smp1) \arrow[dashed,-]{dddddd}\\ 
&(m\smp5,m\smp1)\urar\arrow[dotted,-]{dddd}& &(m\smp6,m\smp2) \arrow[dotted,-]{dddd}\arrow[dotted,-]{rrrr} &&&&(m\smp3,m\smm1)  \urar\arrow[dotted,-]{dddd} & &(m\smp4,m)\arrow[dotted,-]{dddd} \urar \\
& \quad \\
 & \quad \\
 & \quad \\
&(1,4)\drar&&(2,5)\arrow[dotted,-]{rrrr}&& &&(2m,2)\drar&&(2m\smp1,3) \drar \\
(1,3)\urar\drar & &
(2,4)\urar\drar\arrow[dotted,-]{rrrrrr} &&& && &(2m\smp1,2)\drar \urar & & (1,3)\\
&(2,3)\urar&&(3,4)\arrow[dotted,-]{rrrr} &&& &(2m\smp1,1)\urar& &(1,2) \urar\\
\phantom{X}\arrow[dashed,-]{uu}& & & & & &&&&& \phantom{X}\arrow[dashed,-]{uu}
\end{tikzcd}
\]}
\caption{$\CM \Lambda$ for $n = 2m+1$}
\label{CMLo}
\end{figure}

\subsection{Classification of Cohen-Macaulay $\Lambda$-modules}

Let $\cs$ be the set of tagged arcs and sides of the once-punctured polygon $P^*$. In this subsection, we prove the following theorem:

\begin{theorem} \label{thm:classif}
 \begin{enumerate}
  \item There is a bijection between $\cs$ and the set of isomorphism classes of indecomposable Cohen-Macaulay $\Lambda$-modules given by $a \mapsto M_a$ ($M_a$ is defined in Theorem \ref{cor:lambdamod}).
  \item Any Cohen-Macaulay $\Lambda$-module is isomorphic to $\bigoplus_{a \in \cs} M_a^{l_a}$ for some non-negative integers $l_a$. Moreover, the $l_a$ are uniquely determined.
 \end{enumerate} 
\end{theorem}

\begin{remark}
Theorem \ref{thm:classif} shows that Krull-Schmidt-Azumaya property is valid in this case. This is interesting by itself since our base ring $R=K[X]$ is not even local, and in such a case, Krull-Schmidt-Azumaya property usually fails.
\end{remark}

First of all, it is immediate that the $M_a$ are non-isomorphic indecomposable Cohen-Macaulay $\Lambda$-modules. To prove this theorem, let us define the following elements of $\Lambda$:
$$\alpha_i = E_{i,i+1} \quad \alpha_n = X E_{n,1} \quad \beta_i = Y E_{i+1,i} \quad \beta_n = X^{-1}Y E_{1,n}.$$
Together with the idempotents $E_{ii}$, they generate $\Lambda$ as an $R'$-algebra and satisfy the relations
$$\left\{ \begin{array}{l} \alpha_i \alpha_{i+1} \cdots  \alpha_{i-1} = X E_{ii} \\ \beta_{i-1} \beta_{i-2} \cdots  \beta_{i} = Y^{n-1} E_{ii} \\ \alpha_i \beta_i =  \beta_{i-1} \alpha_{i-1} = Y E_{ii}\end{array}\right.$$
for $i \in \cci{1}{n}$. In fact, the quiver of $\Lambda$ is given as follows:
 
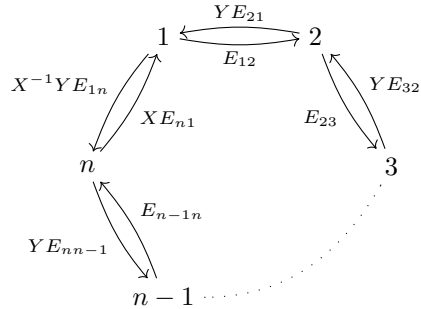
\begin{center}  
\label{remark:quiver}
  \begin{tikzpicture}[commutative diagrams/every diagram]
  \node (a) at (0:2)      {$3$};
  \node (b) at (0+60:2)   {$2$};
  \node (c) at (0+60*2:2) {$1$};
  \node (d) at (0+60*3:2) {$n$};
  \node (e) at (0+60*4:2) {$n-1$};
  \path[commutative diagrams/.cd, every arrow, every label]
  (a) edge[out=110,in=-50,commutative diagrams/rightarrow] node[swap] {$YE_{32}$} (b)
  (b) edge[out=170,in=10,commutative diagrams/rightarrow] node[swap] {$YE_{21}$}  (c)
  (c) edge[out=230,in=70,commutative diagrams/rightarrow] node[swap] {$X^{-1}YE_{1n}$}  (d)
  (d) edge[out=-70,in=130,commutative diagrams/rightarrow] node[swap] {$YE_{nn-1}$} (e)
  (e) edge[out=110,in=-50,commutative diagrams/rightarrow]  node[swap] {$E_{n-1n}$}(d)
  (d) edge[out=50,in=-110,commutative diagrams/rightarrow]  node[swap] {$XE_{n1}$} (c)
  (c) edge[out=-10,in=190,commutative diagrams/rightarrow] node[swap] {$E_{12}$}  (b)
  (b) edge[out=-70,in=130,commutative diagrams/rightarrow] node[swap] {$E_{23}$} (a)
  (a) edge[out=-120,in=0,loosely dotted,commutative diagrams/path] (e);
\end{tikzpicture}
\end{center}

\begin{lemma} \label{arbase}
 Let $r \in (Y)^{\oplus m}$, $s \in R'^{\oplus p}$ and $t \in (X-Y)^{\oplus q}$ be three vectors such that the ideal $I$ generated by their entries includes the ideal $(X,Y)$. Then there exists an invertible $(m+p+q) \times (m+p+q)$ matrix
 $$G = \left( \begin{matrix}
           A & B & 0 \\
	   C & D & E \\
	   0 & F & G
          \end{matrix} \right)$$
 with coefficients in $R'$ where $B$ has coefficients in $(Y)$ and $F$ has coefficients in $(X-Y)$ such that
 \begin{itemize}
  \item either $G \tr{\vecht{r}{s}{t}}$ contains one $1$ in its second block and $0$ everywhere else;
  \item either $G \tr{\vecht{r}{s}{t}}$ contains one $Y$ in the first or second block, one $X-Y$ in the second or third block and $0$ everywhere else.
 \end{itemize}
\end{lemma}

\begin{proof}
 The proof mainly relies on Euclidean algorithm. We can write $r = r' Y$, $s = s' + s'' Y$ and $t = t' (X-Y)$ where $r' \in R^{\oplus m}$, $s', s'' \in R^{\oplus p}$ and $t' \in R^{\oplus q}$. Up to applying Euclidean algorithm on the entries of $s'$ and then on the entries of $s''$ (which is multiplying by an invertible matrix on the left), we can suppose that 
 $$s = \tr{\vech{Q_1 + Q_2 Y}{Q_3 Y}}$$
 for some $Q_1, Q_2, Q_3 \in R$ (we can ignore $0$ entries). With the same method, we can suppose that $r$ has only one (nonzero) entry $P Y$ and $t$ has only one (nonzero) entry $S(X - Y)$. Using a sequence of (authorized) matrix multiplications
 \begin{align*} \mat{1}{0}{T(X-Y)}{1} \vecv{Q_1+Q_2 Y}{S(X-Y)} &= \vecv{Q_1+Q_2 Y}{(S + T Q_1) (X-Y)} \\
      \mat{1}{T}{0}{1} \vecv{Q_1+Q_2 Y}{S(X-Y)} &= \vecv{(Q_1+TS X)+(Q_2-TS) Y}{S (X-Y)}, 
 \end{align*}
 thanks to the Euclidean algorithm, we can ensure that $Q_1 = 0$ or $S = 0$. In the same way, we can ensure that $Q_3 = 0$ or $P = 0$. If $Q_1 = 0$, we can ensure that at most one of $Q_2$, $Q_3$ and $P$ is nonzero. To summarize, and forgetting about $0$ entries, we can ensure to be in one of the following four cases:
 \begin{enumerate}
  \item $r = PY$, $s = 0$, $t = S(X-Y)$: In this case, by our assumption, $$(Y) \oplus (X-Y) = (X, Y) \subset I = (PY, S(X-Y)) = (PY) \oplus (S(X-Y)),$$ 
   so we deduce that $(Y) = (PY)$ and $(X-Y) = (S(X-Y))$ so up to scalar multiplication $P = S = 1$. This is one of the expected cases.
  \item $r = 0$, $s = Q Y$, $t = S(X-Y)$: In this case, by the same reasoning than before, up to scalar multiplication $Q = S = 1$. This is one of the expected cases.
  \item $r = PY$, $s = Q_1 + Q_2 Y$, $t = 0$: In this case, let us rewrite $s = Q_1' + Q_2' (X-Y)$ where $Q_2' = - Q_2$ and $Q_1' = Q_1 + Q_2 X$. Up to a sequence of authorized matrix multiplications
  \begin{align*} 
      \mat{1}{0}{T}{1} \vecv{P Y}{ Q_1' + Q_2' (X-Y)} &= \vecv{P Y}{(Q_1'+TPX) + (Q_2' - TP) (X-Y)} \\
      \mat{1}{T Y}{0}{1} \vecv{P Y}{ Q_1' + Q_2' (X-Y)}  &= \vecv{(P+T Q_1') Y}{ Q_1' + Q_2' (X-Y)} 
  \end{align*}
  we can suppose that $P = 0$ or $Q_1' = 0$. If $Q_1' = 0$, we are in the previous case and we can conclude. If $P = 0$, $I$ is a principal ideal containing $(X, Y)$ so $I = R'$ and, up to scalar multiplication, $Q_1' + Q_2' (X-Y) = 1$. We are again in an expected case.
  \item $r = 0$, $s = \vech{Q_1 + Q_2 Y}{Q_3 Y}$, $t = 0$: this case in similar to the previous one. \qedhere
 \end{enumerate}
\end{proof}

\begin{lemma} \label{drt}
 Let $M = R' \oplus M_2$ and $N$ be Cohen-Macaulay $R'$-modules and $f: M \rightarrow N$, $g,g' : N \rightarrow M$ be morphisms satisfying $gf = X\Id_M$, $fg = X \Id_N$, $g' f = Y  \Id_M$ and $f g' = Y \Id_N$. There exists an isomorphism $\phi: N \rightarrow N_1 \oplus N_2$ such that 
 $$\phi f = \mat{\psi_{11}}{\psi_{12}}{0}{\psi_{22}} \quad ; \quad g \phi^{-1} = \mat{\chi_{11}}{\chi_{12}}{0}{\chi_{22}} \quad ; \quad g' \phi^{-1} = \mat{\chi'_{11}}{\chi'_{12}}{0}{\chi'_{22}}$$
 where 
 \begin{itemize}
  \item either $N_1 = R'$, $\psi_{11} = 1$, $\chi_{11} = X$ and $\chi_{11}' = Y$;
  \item or $N_1 = (X, Y)$, $\psi_{11} = X$, $\chi_{11}$ is the inclusion and $\chi'_{11}$ maps both $X$ and $Y$ to $Y$.
 \end{itemize}
\end{lemma}

\begin{proof}
 Let $f_1: R' \rightarrow N$ be  
 $$f_1 = f \circ \vecv{\Id_{R'}}{0}.$$ 
 As $(Y)$, $R'$ and $(X-Y)$ are the only isomorphism classes of indecomposable Cohen-Macaulay $R'$-modules, we can decompose, up to isomorphism of $N$:
 $$N = (Y)^{\oplus m} \oplus R'^{\oplus p} \oplus (X-Y)^{\oplus q} \quad ; \quad f_1 = \vecvt{r}{s}{t}$$
 where $r$ is a vector with entries in $(Y)$, $s$ is a vector with entries in $R'$ and $t$ is a vector with entries in $(X-Y)$. Using $g f =X \Id_N$ and $g' f = Y \Id_N$, we obtain that the ideal generated by the entries of $r$, $s$ and $t$ contains $(X,Y)$, so, thanks to Lemma \ref{arbase}, up to multiplying $f$ on the left by an invertible matrix and reordering the rows, we can suppose that we are in one of the following cases:

\begin{enumerate}
 \item $N = R' \oplus N_2$ and
  $$f = \mat{1}{*}{0}{*}.$$
  In this case, we can write 
  $$g = \mat{X}{*}{0}{*} \quad \text{and} \quad g' = \mat{Y}{*}{0}{*}$$
  using the identities $gf = X \Id_N$ and $g' f = Y \Id_N$. We are in the first expected case.
 \item $N = (Y) \oplus R' \oplus N_2$ and 
  $$f = \mattd{Y}{*}{X-Y}{*}{0}{*}.$$
  In this case, we can write 
  $$g = \matdt{\iota_Y}{1}{*}{*}{*}{*} \quad \text{and} \quad g' = \matdt{\iota_Y}{0}{*}{*}{*}{*}$$ 
  Up to column operations on $g$ and corresponding row operations on $f$, we can write
  $$f = \mattd{Y}{*}{X}{0}{0}{*} \quad g = \matdt{0}{1}{0}{*}{*}{*} \quad g' = \matdt{\iota_Y}{0}{*}{*}{*}{*}$$ 
  (the $0$ on the second column of $f$ comes from $gf = X \Id_N$). It is now easy to see that we cannot get $f g' = Y \Id_M$. So this case is excluded.
 \item $N = (Y) \oplus (X-Y) \oplus N_2$ and 
  $$f = \mattd{Y}{*}{X-Y}{*}{0}{*}.$$
  In this case, we can write 
  $$g = \matdt{\iota_Y}{\iota_{X-Y}}{*}{0}{0}{*} \quad \text{and} \quad g' = \matdt{\iota_Y}{0}{*}{0}{0}{*}$$ 
  (once again, we use that $(Y)$ and $(X-Y)$ are in direct sum). Using the equality $(Y) \oplus (X-Y) = (X,Y)$, we are in the second expected case.
 \item $N = R' \oplus R' \oplus N_2$ and 
  $$f = \mattd{Y}{*}{X-Y}{*}{0}{*}.$$
  We can write 
  $$g = \matdt{1}{1}{*}{*}{*}{*} \quad \text{and} \quad g' = \matdt{1}{0}{*}{*}{*}{*}$$ 

  As in (ii), using column operations on $g$, we can rewrite:
  $$f = \mattd{Y}{*}{X}{0}{0}{*} \quad g = \matdt{0}{1}{0}{*}{*}{*} \quad g' = \matdt{1}{0}{*}{*}{*}{*}$$ 
  and it contradicts $fg' = Y \Id_M$.

 \item $N = R' \oplus (X-Y) \oplus N_2$ and 
  $$f = \mattd{Y}{*}{X-Y}{*}{0}{*}.$$
  We can write 
  $$g = \matdt{1}{\iota_{X-Y}}{*}{*}{*}{*} \quad \text{and} \quad g' = \matdt{1}{0}{*}{*}{*}{*}$$ 

  Using column operations on $g$ and $g'$, we can rewrite:
  $$f = \mattd{Y}{0}{X-Y}{*}{0}{*} \quad g = \matdt{1}{\iota_{X-Y}}{*}{*}{*}{*} \quad g' = \matdt{1}{0}{0}{*}{*}{*}$$ 
  (the $0$ in the second column of $f$ comes from $g'f = Y \Id_N$). We cannot have $fg = X \Id_M$ so this case is excluded. \qedhere
 \end{enumerate}
\end{proof}

We can easily dualize (over $R'$) the previous lemma:
\begin{lemma} \label{gch}
 Let $M = R' \oplus M_2$ and $N$ be Cohen-Macaulay $R'$-modules and $f: N \rightarrow M$, $g,g' : M \rightarrow N$ be morphisms satisfying $gf = X\Id_N$, $fg = X \Id_M$, $g' f = Y  \Id_N$ and $f g' = Y \Id_M$. There exists an isomorphism $\phi: N \rightarrow N_1 \oplus N_2$ such that 
 $$f \phi^{-1} = \mat{\psi_{11}}{0}{\psi_{21}}{\psi_{22}} \quad ; \quad \phi g = \mat{\chi_{11}}{0}{\chi_{21}}{\chi_{22}} \quad ; \quad \phi g' = \mat{\chi'_{11}}{0}{\chi'_{21}}{\chi'_{22}}$$
 where 
 \begin{itemize}
  \item either $N_1 = R'$, $\psi_{11} = 1$, $\chi_{11} = X$ and $\chi_{11}' = Y$;
  \item or $N_1 = (X, Y)$, $\psi_{11}$ is the inclusion, $\chi_{11} = X$ and $\chi'_{11} = Y$.
 \end{itemize}
\end{lemma}

\begin{lemma} \label{clas1}
 Let $M$ be a Cohen-Macaulay $\Lambda$-module. If $M$, as an $R'$-module, has a direct summand isomorphic to $R'$, then $M$ has a direct summand isomorphic to $M_a$ for some tagged arc or side $a$ of $P^*$ which is not incident to the puncture.
\end{lemma}

\begin{proof}
 For $i \in \cci{1}{n}$, let $M_i = E_{ii} M$. By abuse of notation, we call $\alpha_i: M_{i+1} \rightarrow M_i$ and $\beta_i: M_i \rightarrow M_{i+1}$ the morphisms of $R'$-modules corresponding to the elements with same names in $\Lambda$.

 Let $i, j \in \cci{1}{n}$ such that $\alpha_i \alpha_{i+1} \dots \alpha_{j-1} \alpha_j$ has a direct summand isomorphic to $$R' \xrightarrow{X} R'$$ (such a pair exists as $M$ contains $R'$ as a direct summand, and $\alpha_i \alpha_{i+1} \dots \alpha_{i-1} = X \Id_{M_i}$ for any $i \in \cci{1}{n}$). If $j < i$, note that $\alpha_n \alpha_1$ appears in the previous composition. The number $k$ of factors of this composition is $d(i,j)+1$. We make the additional assumption that $k$ is as small as possible. Without loss of generality, we can suppose that $i = 1 \leq j$ (the problem is invariant by cyclic permutation). Using Lemma \ref{drt} for $f=\alpha_j$, $g=\alpha_{j+1} \alpha_{j+2} \dots \alpha_{j-1}$ and $g'=\beta_j$, we get actually that $j > 1$ and we can suppose that $M_{j+1} = R' \oplus M_{j+1}'$, $M_j = (X, Y) \oplus M_j'$ and
 $$\alpha_j = \mat{X}{\alpha_{j,12}}{0}{\alpha_{j,22}}$$
 (the other possibility of Lemma \ref{drt} would contradicts the minimality of $k$). Then, we get easily that we can write $M_1 = R' \oplus M_1'$ and
 $$\gamma = \alpha_1 \dots \alpha_{j-1} = \mat{\iota_{(X, Y)}}{\gamma_{12}}{0}{\gamma_{22}}$$
 where $\iota_{(X, Y)}: (X, Y) \rightarrow R'$ is the inclusion. 
 Note that by hypothesis
 $$\gamma \alpha_j = \mat{X}{0}{0}{*}.$$

 As all morphisms to $R'$ which are in the radical of $\CM R'$ factor through $\iota_{(X,Y)}: (X,Y) \rightarrow R'$, by column operations on $\gamma$ which do not affect the previous shapes, we can then suppose that one of the following holds:
 \begin{itemize}
  \item $\gamma_{12} = 0$;
  \item $M_j' = R' \oplus M_j''$ and 
   $$\gamma = \matdt{\iota_{(X,Y)}}{1}{0}{0}{*}{*} \quad \text{and} \quad \alpha_j = \mattd{X}{\alpha_{j,12}}{0}{*}{0}{*}$$ 
   then by a column operation on $\gamma$, we get
   $$\gamma = \matdt{0}{1}{0}{*}{*}{*} \quad \text{and} \quad \alpha_j = \mattd{X}{\alpha_{j,12}}{X}{0}{0}{*}$$ 
   (the $0$ in the second column of $\alpha_j$ comes from the shape of $\gamma \alpha_j$).
   But it contradicts the existence of $\beta_j$ such that $\alpha_j \beta_j = Y \Id_{M_j}$.
 \end{itemize}
 Finally, we get the situation
 $$\gamma = \mat{\iota_{(X,Y)}}{0}{0}{\gamma_{22}} \quad \text{and} \quad \alpha_j = \mat{X}{0}{0}{\alpha_{j,22}}$$ 
 
 Now, using Lemma \ref{gch} if $j > 2$ permits to suppose that
 $$\alpha_1 = \mat{\iota_{(X,Y)}}{0}{\alpha_{1, 21}}{\alpha_{1, 22}}$$
 Then we get easily that
 $$\gamma' = \alpha_2 \dots \alpha_{j-1} = \mat{1}{0}{\gamma'_{21}}{\gamma'_{22}}.$$
 By row operations on $\gamma'$ (and the corresponding ones on $\alpha_1$, we can suppose that
 $$\alpha_1 = \mat{\iota_{(X,Y)}}{0}{0}{\alpha_{1,22}} \quad \text{and} \quad \gamma' = \mat{1}{0}{0}{\gamma'_{22}}.$$
 We also get easily that
 $$\gamma'' = \alpha_{j+1} \dots \alpha_n = \mat{1}{0}{0}{\gamma''_{22}}.$$
 Acting by automorphism on $M_3$, \dots, $M_{j-1}$ if $j > 3$ and on $M_{j+2}$, \dots, $M_n$ if $j < n-1$ permits easily to suppose that
 $$\alpha_\ell = \mat{1}{0}{0}{\alpha_{\ell, 22}}$$
 for any $\ell \in \cci{2}{j-1} \cup \cci{j+1}{n}$. Then, we conclude that $M$ has a direct summand isomorphic to
 \spovermat{n} 
 \begin{align*}& \tr{\begin{bmatrix} R' & \bovermat{j-1}{(X,Y) \cdots (X,Y)} & \bovermat{n-j}{(X) \cdots (X)} \end{bmatrix} } \cong M_{(P_1, P_j)}. \qedhere\end{align*}
\end{proof}

\begin{lemma} \label{clas2}
 Let $M$ be a Cohen-Macaulay $\Lambda$-module. If $M$, as an $R'$-module, has no direct summand isomorphic to $R'$ then $M$ has a direct summand isomorphic to some $M_a$ where $a$ is a tagged arc of $P^*$ incident to the puncture.
\end{lemma}

\begin{proof}
 Denote as before $M_i = E_{ii} M$. As an $R'$-module, $M$ is a direct sum of copies of $(Y)$ and $(X-Y)$. As there are no morphisms between $(Y)$ and $(X-Y)$, we can suppose that only one of them appear as a summand of $M$ and therefore the matrix coefficients of the $\alpha_i$ are just elements of $R$. Up to circular permutation, we can suppose that $\alpha_n$ is not invertible. Choose an $R$-basis $\{e_1, \dots, e_\ell\}$ of $M_n$ such that $e_1$ is not in the image of $\alpha_n$. By usual Euclidean algorithm applied on the right of $\alpha_n$, we can suppose that
 $$\alpha_n = \mat{\lambda}{0}{*}{*} \quad \text{and} \quad \alpha_1 \dots \alpha_{n-1} = \mat{\lambda'}{0}{*}{*}$$
 As $\lambda' \lambda = X$ and $e_1$ is not reached by $\alpha_n$, we can suppose up to scalar change of basis that $\lambda = X$ and $\lambda' = 1$. Hence, by row operations on $\alpha_1 \dots \alpha_{n-1}$, we can suppose that
 $$\alpha_n = \mat{X}{0}{0}{*} \quad \text{and} \quad \alpha_1 \dots \alpha_{n-1} = \mat{1}{0}{0}{*}$$
 (the lower left $0$ of $\alpha_n$ comes from $\alpha_n \alpha_1 \alpha_2 \dots \alpha_{n-1} = X \Id_{M_n}$). Therefore, by changes of basis of $M_2$, \dots, $M_{n-1}$, we can suppose that
 $$\alpha_\ell = \mat{1}{0}{0}{*}$$
 for $\ell \in \cci{1}{n-1}$. Finally, $M$ as a direct summand isomorphic to 
\spovermat{n} 
 \begin{align*}& \tr{\begin{bmatrix} \bovermat{n}{(Y) \cdots (Y)}  \end{bmatrix} } \cong M_{(P_{n}, *)} \\
  \\
 \text{or} \quad & \tr{\begin{bmatrix} \bovermat{n}{(X-Y) \cdots (X-Y)}  \end{bmatrix} } \cong M_{(P_{n}, \bowtie)} .
 \qedhere\end{align*}
\end{proof}

\begin{proof}[Proof of Theorem \ref{thm:classif}]
 First of all, thanks to Lemmas \ref{clas1} and \ref{clas2}, any Cohen-Macaulay $\Lambda$-module can be decomposed as expected. 

 For the uniqueness of the decomposition, we need to use Proposition \ref{morph} (notice that we do not use Theorem \ref{thm:classif} in its proof). The endomorphism algebra of $M_a$ is isomorphic to $R'$ if $a$ is not incident to the puncture and isomorphic to $R$ if $a$ is incident to the puncture. Moreover, any endomorphism factorizing through another indecomposable is in the ideal $(X,Y)$ in the first case and $(X)$ in the second case. Thus, if we denote $\hat \Lambda = K\llbracket X \rrbracket \otimes_{R} \Lambda$, and consider the functor $ K\llbracket X \rrbracket \otimes_{R} -: \CM \Lambda \rightarrow \CM \hat \Lambda$, non-isomorphic indecomposable objects are mapped to non-isomorphic objects, which are also indecomposable. Moreover, the endomorphism rings of the objects $K \llbracket X \rrbracket \otimes_{R} M_a$ are local so we get the uniqueness of the decomposition of objects in the essential image of the functor $K\llbracket X \rrbracket \otimes_{R} -$. It permits to conclude.
\end{proof}

\subsection{Homological structure of $\CM \Lambda$}

The aim of this subsection is to compute spaces of morphisms and extensions in the category $\CM \Lambda \cong \CM \Lambda'$ where $\Lambda' = e_F \Gamma_\sigma e_F$. For convenience of notation, we will work on $\Lambda'$. Notice that the definitions of $a \vdash b$ and $A_{a,b}$ given for two tagged arcs $a$ and $b$ before Proposition \ref{inductD} make sense even when $a$ and $b$ are not compatible (there are cases where $a \vdash b$ other than the one depicted there).

\begin{proposition} \label{morph}
 Let $a$ and $b$ be two tagged arcs or sides of $P^*$. Following the notation of Proposition \ref{inductD}, we have the following isomorphism:
 $$\Hom_{\Lambda'}(M_a, M_b) \cong A_{a,b}.$$
 Moreover, these morphisms are realized by right multiplication in $\cR'$ and therefore, the composition of morphisms corresponds to the multiplication in $\cR'$.
\end{proposition}

\begin{proof}
 First of all, for $i \in \cci{1}{n}$, using Theorem \ref{cor:lambdamod}, recall that $E_{ii} M_a \cong A_{i,a}$ and $E_{ii} M_b \cong A_{i,b}$ (in a compatible way with the $\Lambda'$-module structure). Thus, we know that for any $i$, $\Hom_{R'}(A_{i,a}, A_{i,b})$ can be realized as an $R'$-submodule of $\cR'$ through multiplication. Namely, 
  $$\begin{array}{|c|c|c|c|}
      \hline
      \text{\backslashbox{$A_{i,a}$}{$A_{i,b}$}} & u^{j'-1} v R' & u^{j'-1} (u-v) R' & u^{j'} R' \\
      \hline
     u^{j-1} v R' & u^{j'-j-1} v R' & 0 & u^{j'-j+2n-1} v R' \\
      \hline
     u^{j-1} (u-v) R' & 0 & u^{j'-j-1} (u-v) R' & u^{j'-j+2n-1} (u-v) R'\\
      \hline
     u^j R' & u^{j'-j-1} v R' & u^{j'-j-1} (u-v) R' & u^{j'-j} R' \\
      \hline
   \end{array} $$
 where $j=\ell^\theta_{i,a}$ and $j' = \ell^\theta_{i,b}$ (the only other kind of $A_{i,a}$ or $A_{i,b}$ which can appear is $u^j R' + u^{j-1} v R' = u^{j-1} ((u-v) R' \oplus v R')$ which can be realized as the direct sum of the two first rows, and two first columns. In any of these cases, the sum is direct inside $\cR'$).

 If $f \in \Hom_{\Lambda'} (M_a, M_b)$, let $f_i \in \Hom_{R'}(A_{i,a}, A_{i,b})$ be its $i$-th component. As $u^2 E_{i, i+1} \in \Lambda'$, for any $m \in M_a$, $f(u^2 E_{i,i+1} m) = u^2 E_{i,i+1} f(m)$ holds. It can be rewritten as $f_i(u^2 m_{i+1}) = u^2 f_{i+1} (m_{i+1})$ or again $f_i u^2 m_{i+1} = u^2 f_{i+1} m_{i+1}$ if $f_i, f_{i+1}$ are considered as elements of $\cR'$. As $u^2$ is invertible in $\cR'$, we get $f_i m_{i+1} = f_{i+1} m_{i+1}$. This is true for any $m_{i+1} \in A_{i+1, a}$ so $f_i - f_{i+1}$ is in the annihilator of $A_{i+1,a}$. Looking at Theorem \ref{cor:lambdamod}, the annihilators of $A_{i+1,a}$ and $A_{i,a}$ are the same and included in
 \begin{itemize}
  \item $(u-v) \cR'$ if $a$ is incident to the puncture and plain;
  \item $v \cR'$ if $a$ is incident to the puncture and notched;
  \item $0$ if $a$ is not incident to the puncture.
 \end{itemize}
 Moreover, looking at the previous table, 
 \begin{itemize}
  \item $\Hom_{R'}(A_{i,a}, A_{i,b}) + \Hom_{R'}(A_{i+1,a}, A_{i+1,b}) \subset v \cR'$  if $a$ is incident to the puncture and plain;
  \item $\Hom_{R'}(A_{i,a}, A_{i,b}) + \Hom_{R'}(A_{i+1,a}, A_{i+1,b}) \subset (u-v) \cR'$ if $a$ is incident to the puncture and notched.
 \end{itemize}
 so $\Hom_{R'}(A_{i,a}, A_{i,b}) + \Hom_{R'}(A_{i+1,a}, A_{i+1,b})$ intersects the annihilator of $A_{i,a}$ at $0$ and we obtain that $f_i = f_{i+1}$.

 Finally, we get
 $$\Hom_{\Lambda'} (M_a, M_b) = \bigcap_{i=1}^n \Hom_{R'} (A_{i,a}, A_{i,b})$$
 as $R'$-submodules of $\cR'$.

 (a) Suppose now that none of $a$ and $b$ is incident to the puncture. For $i \in \cci{1}{n}$, we have
 $$A_{i,a} = \left\{ \begin{array}{ll} u^{\ell^\theta_{(P_i,P_{i+1}), a}} R' & \text{if } (P_i,P_{i+1}) \notvdash a \\ u^{\ell^\theta_{(P_i,P_{i+1}), a}} (R' + u^{-1}vR') & \text{if } (P_i,P_{i+1}) \vdash a \end{array} \right.$$
 and 
 $$A_{i,b} = \left\{ \begin{array}{ll} u^{\ell^\theta_{(P_i,P_{i+1}), b}} R' & \text{if } (P_i,P_{i+1}) \notvdash b \\ u^{\ell^\theta_{(P_i,P_{i+1}), b}} (R' + u^{-1}vR') & \text{if } (P_i,P_{i+1}) \vdash b. \end{array} \right.$$
 Therefore, we obtain
 $$\Hom_{R'} (A_{i,a}, A_{i,b}) = u^{\ell^\theta_{(P_i, P_{i+1}), b} - \ell^\theta_{(P_i, P_{i+1}), a}} R'$$
  if $(P_i,P_{i+1}) \notvdash a$ and $(P_i,P_{i+1}) \notvdash b$, 
 $$\Hom_{R'} (A_{i,a}, A_{i,b}) = u^{\ell^\theta_{(P_i, P_{i+1}), b} - \ell^\theta_{(P_i, P_{i+1}), a}+2n} (R' + u^{-1}vR')$$
  if $(P_i,P_{i+1}) \vdash a$ and $(P_i,P_{i+1}) \notvdash b$, 
 $$\Hom_{R'} (A_{i,a}, A_{i,b}) = u^{\ell^\theta_{(P_i, P_{i+1}), b} - \ell^\theta_{(P_i, P_{i+1}), a}} (R' + u^{-1}vR')$$
  if $(P_i,P_{i+1}) \vdash b$.

 Using Lemma \ref{caldif}, we obtain
 $$\Hom_{R'} (A_{i,a}, A_{i,b}) = u^{\ell^{\theta}_{a,b} - 2 n \left( \delt{b_2 \in \ooi{b_1}{a_2}}\delt{a_1 \in \cci{a_2}{b_1}} + \delt{i \in \oci{a_1}{b_1}}  \right)} R'$$
  if $(P_i,P_{i+1}) \notvdash a$ and $(P_i,P_{i+1}) \notvdash b$, 
 $$\Hom_{R'} (A_{i,a}, A_{i,b}) = u^{\ell^{\theta}_{a,b} - 2 n \left( \delt{b_2 \in \ooi{b_1}{a_2}}\delt{a_1 \in \cci{a_2}{b_1}} + \delt{i \in \oci{a_1}{b_1}} - 1 \right)} (R' + u^{-1}vR')$$
  if $(P_i,P_{i+1}) \vdash a$ and $(P_i,P_{i+1}) \notvdash b$, 
 $$\Hom_{R'} (A_{i,a}, A_{i,b}) = u^{\ell^{\theta}_{a,b} - 2 n \left( \delt{b_2 \in \ooi{b_1}{a_2}}\delt{a_1 \in \cci{a_2}{b_1}} + \delt{i \in \oci{a_1}{b_1}} \right)} (R' + u^{-1}vR')$$
  if $(P_i,P_{i+1}) \vdash b$.

 Notice that $(P_i,P_{i+1}) \vdash a$ if and only if $i \in \oci{a_1}{a_2}$.

 (a-1) Suppose that $a \notvdash b$. It means that $a_2 \notin \ooi{b_2}{a_1}$ and $b_1 \notin \ooi{b_2}{a_1}$. In this case, $\delt{b_2 \in \ooi{b_1}{ a_2}}\delt{a_1 \in \cci{a_2}{b_1}} = 0$. Taking $i = a_1$, we have $(P_i, P_{i+1}) \notvdash a$ and $(P_i, P_{i+1}) \notvdash b$ and an easy computation gives $\Hom_{R'}(A_{i,a}, A_{i,b}) = u^{\ell^\theta_{a,b}} R'$.
 The only way to get a smaller module would be in the case $(P_i,P_{i+1}) \vdash a$ and $(P_i,P_{i+1}) \notvdash b$ that is $i \in \oci{a_1}{a_2} \cap \oci{b_2}{b_1}$. With the current hypotheses, we get that $\oci{a_1}{a_2} \cap \oci{b_2}{b_1} \subset \oci{a_1}{b_1}$, so actually, we cannot get a smaller module.

 (a-2) Suppose now that $a \vdash b$. It means that $a_2 \in \ooi{b_2}{a_1}$ or $b_1 \in \ooi{b_2}{a_1}$. Let us consider two cases:
 \begin{itemize}
  \item When $b_2 \in \ooi{b_1}{a_2}$ and $a_1 \in \cci{a_2}{b_1}$. Taking $i = b_2$, we have $(P_i, P_{i+1}) \vdash a$ and $(P_i, P_{i+1}) \notvdash b$ and an easy computation gives
   $$\Hom_{R'}(A_{i,a}, A_{i,b}) = u^{\ell^\theta_{a,b}} (R'+u^{-1}vR').$$
 Thanks to the term $\delt{b_2 \in \ooi{b_1}{a_2}}\delt{a_1 \in \cci{a_2}{b_1}}$, submodules that appear for any other $i$ are bigger. 
  \item When $b_2 \notin \ooi{b_1}{a_2}$ or $a_1 \notin \cci{a_2}{b_1}$. In this case, as $b_1 \neq b_2$, we get that in fact $b_1 \in \ooi{b_2}{a_1}$. Taking $i = a_1$, we obtain $(P_i, P_{i+1}) \vdash b$ and 
   $$\Hom_{R'}(A_{i,a}, A_{i,b}) = u^{\ell^\theta_{a,b}} (R'+u^{-1}vR').$$
 Other submodules are bigger as, if $(P_i, P_{i+1}) \notvdash b$, that is $i \in \oci{b_2}{b_1}$, we would have $i \in \oci{a_1}{b_1}$. 
 \end{itemize}

 We finished the case where none of $a$ and $b$ are incident to the puncture. 

 (b) Suppose now that both $a$ and $b$ are incident to the puncture. For $i \in \cci{1}{n}$, we have
 $$A_{i,a} = \left\{ \begin{array}{ll} u^{\ell^\theta_{(P_i,P_{i+1}), a}-1} v R' & \text{if $a$ is plain} \\ u^{\ell^\theta_{(P_i,P_{i+1}), a}-1} (u-v) R' & \text{if $a$ is notched}\end{array} \right.$$
 and 
 $$A_{i,b} = \left\{ \begin{array}{ll} u^{\ell^\theta_{(P_i,P_{i+1}), b}-1} v R' & \text{if $b$ is plain} \\ u^{\ell^\theta_{(P_i,P_{i+1}), b}-1} (u-v) R' & \text{if $b$ is notched.} \end{array} \right.$$
                         
 As there are no morphisms if the tags are different, we can suppose that both $a$ and $b$ are plain, and we obtain, for any $i \in \cci{1}{n}$,
  $$\Hom_{R'} (A_{i,a}, A_{i,b}) = u^{\ell^\theta_{(P_i, P_{i+1}), b} - \ell^\theta_{(P_i, P_{i+1}), a}-1} v R'.$$          

 Using Lemma \ref{caldif} and the fact that $b_1 = b_2$, we deduce that
  $$\Hom_{R'} (A_{i,a}, A_{i,b}) = u^{\ell^\theta_{a, b} - 2n \delt{i \in \oci{a_1}{b_1}}-1} v R'$$
 and for $i = a_1$, $\Hom_{R'} (A_{i,a}, A_{i,b}) = u^{\ell^\theta_{a, b} -1} v R'$ which is of course the smallest possible.

 (c) Suppose now that $a$ is incident to the puncture and $b$ is not. Without loss of generality, we can suppose $a$ is plain. We have $(X-Y) M_a = 0$. Therefore, for any $f \in \Hom_{\Lambda'}(M_a, M_b)$, $(X-Y) \im f = 0$. Notice now that
 $$M_b' = \{m \in M_b \,|\, (X-Y) m = 0\}$$
 satisfies
 \begin{align*} E_{ii} M'_b &= \left\{ \begin{array}{ll} u^{\ell^\theta_{(P_i,P_{i+1}), b}+2n-1}v R' & \text{if } (P_i,P_{i+1}) \notvdash b \\ u^{\ell^\theta_{(P_i,P_{i+1}), b}-1} vR' & \text{if } (P_i,P_{i+1}) \vdash b \end{array} \right. \\ &= u^{\ell^\theta_{(P_i,P_{i+1}), b}+2n\delt{i \in \oci{b_2}{b_1}}-1}v R'. \end{align*}
 Let $b' = (P_{b_2}, *)$. Thanks to Lemma \ref{caldif}, we can rewrite 
 \begin{align*}
  & \ell^\theta_{(P_i,P_{i+1}), b}+2n\delt{i \in \oci{b_2}{b_1}} \\ =\, & \ell^\theta_{(P_i,P_{i+1}), b'} + \ell^\theta_{b', b} - 2n \delt{i \in \oci{b_2}{b_1}}   +2n\delt{i \in \oci{b_2}{b_1}} = \ell^\theta_{(P_i,P_{i+1}), b'} + \ell^\theta_{b', b}.
 \end{align*}
 Thus $M_b' = u^{\ell^\theta_{b', b}}M_{b'}$ and 
 $$\Hom_{\Lambda'}(M_a, M_b) = u^{\ell^\theta_{b', b}} \Hom_{\Lambda'}(M_a, M_{b'}) = u^{\ell^\theta_{b', b}+\ell^\theta_{a, b'} -1} v R' $$
 and we have (using $a_1 = a_2$),
 \begin{align*} 
  \ell^\theta_{b', b}+\ell^\theta_{a, b'} &= d(b_2, b_1) + 2d(a_1, b_2)  \\
    &= d(a_1, b_1) + n \delt{b_2 \in \ooi{b_1}{a_1}} +  d(a_1, b_2) = \ell^\theta_{a, b}
 \end{align*}
 which concludes this case.

 (d) Finally, suppose that $b$ is incident to the puncture and $a$ is not. Without loss of generality, we can suppose $b$ is plain. As $(X-Y) M_b = 0$, 
 $$\Hom_{\Lambda'}(M_a, M_b) = \Hom_{\Lambda'}(M_a', M_b)$$
 where $M_a' = M_a/(X-Y)M_a$. Using the same idea as before, 
 $$E_{ii} M'_a = u^{\ell^\theta_{(P_i,P_{i+1}), a}-1} vR'$$
 and thanks to Lemma \ref{caldif}, if $a' = (P_{a_1}, *)$,
 \begin{align*}\ell^\theta_{(P_i,P_{i+1}), a} &= \ell^\theta_{(P_i,P_{i+1}), a'} + \ell^\theta_{a',a} -2n \left(\delt{a_2 \in \ooi{a_1}{a_1}} \delt{a_1 \in \cci{a_1}{a_1}} + \delt{i\in \oci{a_1}{a_1}}\right) \\ &=\ell^\theta_{(P_i,P_{i+1}), a'} + \ell^\theta_{a',a} -2n \end{align*}
 and therefore, $M'_a = u^{\ell^\theta_{a',a} - 2n} M_{a'}$. Thus
 $$\Hom_{\Lambda'}(M_a, M_b) = u^{2n - \ell^\theta_{a', a}} \Hom_{\Lambda'}(M_{a'}, M_b) = u^{2n - \ell^\theta_{a', a} + \ell^\theta_{a',b}-1} v R' $$
 and, using $b_1 = b_2$, 
 \begin{align*} 
  2n - \ell^\theta_{a', a} + \ell^\theta_{a',b} &= 2n - (d(a_1,a_2) + n) + 2 d(a_1, b_1) \\
      &= 2n - (n - d(a_2,a_1) + n) + 2 d(a_1, b_1) \\
      &= d(a_1, b_1) + d(a_2, b_1) + n \delt{a_1 \in \ooi{b_1}{a_2}} = \ell^\theta_{a,b}.
 \end{align*}

 It concludes the proof. 
\end{proof}

\begin{proposition} \label{stabmorph} 
 Let $a$ and $b$ be two tagged arcs or sides. Let $M_a$ and $M_b$ be the corresponding indecomposable $\Lambda'$-modules. We have the following isomorphism of graded $R'$-modules:

 \begin{itemize}
  \item $\uHom_{\Lambda'}(M_a, M_b) = 0$ if $a$ and $b$ are both incident to the puncture with different tags;
  \item $\uHom_{\Lambda'}(M_a, M_b) = u^{\ell^\theta_{a,b}} (R'/(X,Y))^{\oplus \epsilon}$, where
 $$\epsilon = \delt{a_1-1 \in \ooi{b_1}{b_2}} \delt{b_2+1 \in \ooi{a_1}{a_2}} = \delt{a_2 - 1 \in \ooi{b_1}{b_2}} \delt{b_1+1\in \ooi{a_1}{a_2}}$$  if either $a$ and $b$ are both incident to the puncture with the same tag, or exactly one of them is incident to the puncture.
  \item $\uHom_{\Lambda'}(M_a, M_b) = u^{\ell^\theta_{a,b}} (R'/(X,Y))^{\oplus \epsilon}$, where
 $$\epsilon = \delt{a_1-1 \in \ooi{b_1}{b_2}} \delt{b_2+1 \in \ooi{a_1}{a_2}} + \delt{a_2 - 1 \in \ooi{b_1}{b_2}} \delt{b_1+1\in \ooi{a_1}{a_2}}$$  if none of $a$ and $b$ is incident to the puncture.
 \end{itemize}
\end{proposition}

\begin{proof}
 (a) Suppose first that none of $a$ and $b$ is incident to the puncture. For any $i \in \cci{1}{n}$, let $P_i$ be the projective module corresponding to the arc $(P_i,P_{i+1})$ of the polygon. Thanks to Proposition \ref{morph}, we have
 $$\Hom_{\Lambda'}(M_a,P_i) = \left\{ \begin{array}{ll} u^{\ell^\theta_{a, (P_i,P_{i+1})}} R', & \text{if } a \notvdash (P_i,P_{i+1}); \\ u^{\ell^\theta_{a, (P_i,P_{i+1})}} (R' + u^{-1}v R'), & \text{if } a \vdash (P_i,P_{i+1}); \end{array} \right.$$
 and
 $$\Hom_{\Lambda'}(P_i,M_b) = \left\{ \begin{array}{ll} u^{\ell^\theta_{(P_i,P_{i+1}),b}} R', & \text{if } (P_i,P_{i+1}) \notvdash b; \\ u^{\ell^\theta_{(P_i,P_{i+1}),b}} (R' + u^{-1}v R'), & \text{if } (P_i,P_{i+1}) \vdash b. \end{array} \right.$$
 As a consequence, we get
 $$\Hom_{\Lambda'}(P_i,M_b) \circ \Hom_{\Lambda'}(M_a,P_i) = u^{\ell^\theta_{a, (P_i,P_{i+1})} + \ell^\theta_{(P_i,P_{i+1}),b}} R'$$
 if $a \notvdash (P_i,P_{i+1})$ and $(P_i,P_{i+1}) \notvdash b$ and
 $$\Hom_{\Lambda'}(P_i,M_b) \circ \Hom_{\Lambda'}(M_a,P_i) = u^{\ell^\theta_{a, (P_i,P_{i+1})} + \ell^\theta_{(P_i,P_{i+1}),b}} (R'+u^{-1}vR')$$
 if $a \vdash (P_i,P_{i+1})$ or $(P_i,P_{i+1}) \vdash b$.

 Using Lemma \ref{caldif},
 $$\ell^\theta_{a, (P_i,P_{i+1})} + \ell^\theta_{(P_i,P_{i+1}),b} = \ell^\theta_{a,b} + 2n \left(\delt{a_1\in\oci{a_2-1}{b_1}} \delt{b_2\in\oci{a_2-1}{b_1}} + \delt{i \in \ooi{b_1}{a_2-1}} \right)$$
 The minimum is reached for $i \in \cci{a_2-1}{b_1}$ and is
 $$\ell^\theta_{a,b} + 2n \delt{b_2 \in \oci{a_2-1}{b_1}} \delt{a_1 \in \oci{a_2-1}{b_1}} = \ell^\theta_{a,b} + 2n \delt{a_2-1 \in \ooi{b_1}{b_2}} \delt{b_1+1 \in \ooi{a_1}{a_2}}.$$

 Recall now that $a \vdash (P_i,P_{i+1})$ if and only if $a_2 \in \ooi{i+1}{a_1}$ if and only if $i \in \coi{a_1-1}{a_2-1}$ and $(P_i,P_{i+1}) \vdash b$ in and only if $b_1 \in \ooi{b_2}{i}$ if and only if $i \in \oci{b_1}{b_2}$. So $a \vdash (P_i,P_{i+1})$ or $(P_i,P_{i+1}) \vdash b$ if and only if $i \in \coi{a_1-1}{a_2-1} \cup \oci{b_1}{b_2}$. If $\coi{a_1-1}{a_2-1} \cup \oci{b_1}{b_2}$ intersects $\cci{a_2-1}{b_1}$, we deduce that 
 \begin{align*}\mathcal{P}(M_a,M_b) &=\sum_{i=1}^n \Hom_{\Lambda'}(P_i,M_b) \circ \Hom_{\Lambda'}(M_a,P_i) \\ &= u^{\ell^\theta_{a,b} + 2n \delt{a_2-1 \in \ooi{b_1}{b_2}} \delt{b_1+1 \in \ooi{a_1}{a_2}}} (R'+u^{-1}vR')\end{align*}
 and otherwise
 \begin{align*}\mathcal{P}(M_a,M_b) &= \sum_{i=1}^n \Hom_{\Lambda'}(P_i,M_b) \circ \Hom_{\Lambda'}(M_a,P_i) \\ & = u^{\ell^\theta_{a,b} + 2n \delt{a_2-1 \in \ooi{b_1}{b_2}} \delt{b_1+1 \in \ooi{a_1}{a_2}}} R'.\end{align*}
 
 Notice that $\coi{a_1-1}{a_2-1} \cup \oci{b_1}{b_2}$ intersects $\cci{a_2-1}{b_1}$ if and only if $b_1+1 \in \coi{a_1}{a_2}$ or $a_2-1 \in \oci{b_1}{b_2}$, if and only if $b_1+1 \in \ooi{a_1}{a_2}$ or $a_2-1 \in \ooi{b_1}{b_2}$ or $a_1 = b_1 + 1$ or $a_2 = b_2 + 1$.

 Then, we can simplify $\mathcal{P}(M_a, M_b)$ in the following way:
 \begin{enumerate}[{Case} 1:]
  \item $\mathcal{P}(M_a, M_b) = u^{\ell^\theta_{a,b} + 2n} (R' + u^{-1} v R')$, if $b_1+1 \in \ooi{a_1}{a_2}$ and $a_2-1 \in \ooi{b_1}{b_2}$;
  \item $\mathcal{P}(M_a, M_b) = u^{\ell^\theta_{a,b}} R'$, if $b_1+1 \notin \coi{a_1}{a_2}$ and $a_2-1 \notin \oci{b_1}{b_2}$;
  \item $\mathcal{P}(M_a, M_b) = u^{\ell^\theta_{a,b}} (R' + u^{-1} v R')$, else.
 \end{enumerate}

 Recall also that 
 $$\mathcal{P}(M_a, M_b) \subset \Hom_{\Lambda'}(M_a, M_b) = \left\{ \begin{array}{ll} u^{\ell^\theta_{a,b}} R' & \text{if } a \notvdash b \\
                                           u^{\ell^\theta_{a,b}} (R' + u^{-1} v R') & \text{if } a \vdash b \\
                                          \end{array} \right.$$
 and therefore, in case 3, we will always get $\uHom_{\Lambda'}(M_a, M_b) = 0$. In case 1, if $a \notvdash b$, we get $\uHom_{\Lambda'}(M_a, M_b) \cong u^{\ell^\theta_{a,b}} R'/(X,Y)$ (as graded $R'$-modules); if $a \vdash b$, we get $\uHom_{\Lambda'}(M_a, M_b) \cong u^{\ell^\theta_{a,b}} (R'/(X,Y) \oplus R'/(X,Y))$. In case 2, if $a \notvdash b$, we get $\uHom_{\Lambda'}(M_a, M_b) = 0$; if $a \vdash b$, we get $\uHom_{\Lambda'}(M_a, M_b) \cong u^{\ell^\theta_{a,b}} R'/(X,Y)$. 
 
 Notice that $a \vdash b$ if and only if $a_1-1 \in \coi{b_1}{b_2}$ or $b_2+1 \in \oci{a_1}{a_2}$. Then, by an easy case by case analysis, it concludes the case where none of $a$ and $b$ are incident to the puncture.

 (b) Suppose now that at least one of $a$ and $b$ is incident to the puncture. Without loss of generality, we can suppose that no notched tag appears. An easy computation shows us that, in any case,
 $$\Hom_{\Lambda'}(P_i,M_b) \circ \Hom_{\Lambda'}(M_a,P_i) = u^{\ell^\theta_{a, (P_i,P_{i+1})} + \ell^\theta_{(P_i,P_{i+1}),b}-1} v R'.$$
 
 Using Lemma \ref{caldif},
 $$\ell^\theta_{a, (P_i,P_{i+1})} + \ell^\theta_{(P_i,P_{i+1}),b} = \ell^\theta_{a,b} + 2n \left(\delt{a_1\in\oci{a_2-1}{b_1}} \delt{b_2\in\oci{a_2-1}{b_1}} + \delt{i \in \ooi{b_1}{a_2-1}} \right)$$
 which gives, as before, that the minimum is reached for $i \in \cci{a_2-1}{b_1}$ and is 
 $$\ell^\theta_{a,b} + 2n \delt{a_2-1 \in \ooi{b_1}{b_2}} \delt{b_1+1 \in \ooi{a_1}{a_2}}.$$
 
 Therefore,  $\uHom_{\Lambda'}(M_a, M_b) \cong u^{\ell^\theta_{a,b}} R'/(X,Y)$ if $a_2-1 \in \ooi{b_1}{b_2}$ and $b_1+1 \in \ooi{a_1}{a_2}$ and  $\uHom_{\Lambda'}(M_a, M_b) = 0$ else.
\end{proof}

\begin{proposition} \label{ARseq}
 The category $\CM \Lambda'$ admits the following Auslander-Reiten sequences:
 $$0 \rightarrow M_{(P_i, P_j)} \xrightarrow{\svecv{u}{u}} M_{(P_{i+1},P_j)} \oplus M_{(P_i, P_{j+1})} \xrightarrow{\svech{-u}{u}} M_{(P_{i+1},P_{j+1})} \rightarrow 0\quad (j \neq i, i+1);$$
 $$0 \rightarrow M_{(P_i, *)} \xrightarrow{v} M_{(P_{i+1},P_i)} \xrightarrow{u-v} M_{(P_{i+1}, \bowtie)} \rightarrow 0;$$
 $$0 \rightarrow M_{(P_i, \bowtie)} \xrightarrow{u-v} M_{(P_{i+1},P_i)} \xrightarrow{v} M_{(P_{i+1}, *)} \rightarrow 0.$$
 Notice that $M_{(P_i,P_i)}$, if it appears, has to be interpreted as $M_{(P_i, *)} \oplus M_{(P_i, \bowtie)}$. Thus, $\CM \Lambda'$ admits an Auslander-Reiten translation $\tau$ defined by 
 \begin{align*}
  \tau(M_{(P_i, P_j)}) = M_{(P_{i-1}, P_{j-1})} & \quad \text{if } j \neq i, i+1;\\
  \tau(M_{(P_i, *)}) = M_{(P_{i-1}, \bowtie)};\\
  \tau(M_{(P_i, \bowtie)}) = M_{(P_{i-1}, *)}.
 \end{align*}
\end{proposition}

\begin{proof}
 (a) Consider the first case. Let $a$ be a side or an arc of $P^*$ which is not incident to the puncture or a formal sum $ (P_{a_1}, *) \oplus (P_{a_1}, \bowtie)$ and $f: M_{(P_i,P_j)} \rightarrow M_a$ be a morphism which is not a split monomorphism. According to Proposition \ref{morph}, the degree $\deg(f)$ of $f$ is at least 
 $$\ell^\theta_{(P_i, P_j), a} + 2n \delt{a = (P_i, P_j)}.$$
 Moreover, using the beginning of the proof of Lemma \ref{caldif}, we get the equalities
 \begin{align*}
  &\ell^\theta_{(P_i, P_j), (P_{i+1}, P_j)} + \ell^\theta_{(P_{i+1}, P_j), a} - \ell^\theta_{(P_i, P_j), a} \\
   = \, & n \left( \delt{i = a_1} + 0 + 0 + \left| \delt{i+1 \in \ooi{a_1}{j}} - \delt{a_2 \in \ooi{a_1}{j}}\right| - \left| \delt{i \in \ooi{a_1}{j}} - \delt{a_2 \in \ooi{a_1}{j}}\right| \right) \\
   = \, & n \left( \delt{i = a_1} +  \left(\delt{a_2 \in \cci{j}{a_1}} - \delt{a_2 \in \ooi{a_1}{j}} \right) \left( \delt{i+1 \in \ooi{a_1}{j}} - \delt{i \in \ooi{a_1}{j}} \right) \right)  \\
   = \, & n \left( \delt{i = a_1} +  \left(2 \delt{a_2 \in \cci{j}{a_1}} - 1 \right) \delt{i = a_1} \right) = 2 n \delt{i = a_1} \delt{a_2 \in \cci{j}{a_1}}
 \end{align*}
 so $\deg(f) \geq \ell^\theta_{(P_i, P_j), (P_{i+1}, P_j)} + \ell^\theta_{(P_{i+1}, P_j), a} + 2n(\delt{a = (P_i, P_j)} - \delt{i = a_1} \delt{a_2 \in \cci{j}{a_1}})$
 and
 \begin{align*}
  &\ell^\theta_{(P_i, P_j), (P_i, P_{j+1})} + \ell^\theta_{(P_i, P_{j+1}), a} - \ell^\theta_{(P_i, P_j), a} \\
   = \, & n \left( 0 + \delt{j = a_2} + 0 + \left| \delt{i \in \ooi{a_1}{j+1}} - \delt{a_2 \in \ooi{a_1}{j+1}}\right| - \left| \delt{i \in \ooi{a_1}{j}} - \delt{a_2 \in \ooi{a_1}{j}}\right| \right) \\
   = \, & n \left( \delt{j = a_2} - \left( \delt{i \in \ooi{a_1}{j}} - \delt{i \in \cci{j}{a_1}} \right) \delt{j = a_2} \right) = 2n \delt{j = a_2} \delt{i \in \cci{j}{a_1}}
 \end{align*}
 so $\deg(f) \geq \ell^\theta_{(P_i, P_j), (P_i, P_{j+1})} + \ell^\theta_{(P_i, P_{j+1}), a} + 2n(\delt{a = (P_i, P_j)} -\delt{j = a_2} \delt{i \in \cci{j}{a_1}})$.
 As at least one of $\delt{a = (P_i, P_j)} - \delt{i = a_1} \delt{a_2 \in \cci{j}{a_1}}$ and $\delt{a = (P_i, P_j)} -\delt{j = a_2} \delt{i \in \cci{j}{a_1}}$ is non-negative, we get
  $$\deg(f) \geq \min\left( \ell^\theta_{(P_i, P_j), (P_{i+1}, P_j)} + \ell^\theta_{(P_{i+1}, P_j), a} , \ell^\theta_{(P_i, P_j), (P_i, P_{j+1})} + \ell^\theta_{(P_i, P_{j+1}), a} \right).$$

  Suppose that $a_1 \neq i$ and $a_2 \neq j$. In this case, $$\ell^\theta_{(P_i, P_j), (P_{i+1}, P_j)} + \ell^\theta_{(P_{i+1}, P_j), a} = \ell^\theta_{(P_i, P_j), (P_i, P_{j+1})} + \ell^\theta_{(P_i, P_{j+1}), a}.$$  Notice that if $(P_i,P_j) \vdash a$, \emph{i.e.} $j \in \ooi{a_2}{i}$ or $a_1 \in \ooi{a_2}{i}$ then we have $j \in \ooi{a_2}{i+1}$ or $a_1 \in \ooi{a_2}{i+1}$ or $j+1 \in \ooi{a_2}{i}$ or $a_1 \in \ooi{a_2}{i}$, \emph{i.e.} $(P_{i+1},P_j) \vdash a$ or $(P_i,P_{j+1}) \vdash a$. From that fact and easy observations, we get that 
   $$f \in \Hom_{\Lambda'}\left(M_{(P_{i+1}, P_j)} \oplus M_{(P_i,P_{j+1})}, M_a\right) u.$$

 Suppose that $a_1 = i$. Then $\deg(f) \geq \ell^\theta_{(P_i, P_j), (P_i, P_{j+1})} + \ell^\theta_{(P_i, P_{j+1}), a}$.

 If $(P_i,P_j) \vdash a$, \emph{i.e.} $j \in \ooi{a_2}{i}$ then we have $j+1 \in \ooi{a_2}{i}$ or $j+1 = i$, \emph{i.e.} $(P_i,P_{j+1}) \vdash a$ or $(P_i, P_{j+1}) = (P_i, *) \oplus (P_i, \bowtie)$. From that fact and easy observations, we get that $$f \in \Hom_{\Lambda'}\left(M_{(P_i,P_{j+1})}, M_a\right) u.$$ 

 Suppose that $a_2 = j$. Then $\deg(f) \geq \ell^\theta_{(P_i, P_j), (P_{i+1}, P_j)} + \ell^\theta_{(P_{i+1}, P_j), a}$.

 If $(P_i,P_j) \vdash a$, \emph{i.e.} $a_1 \in \ooi{a_2}{i}$ then we have $a_1 \in \ooi{a_2}{i+1}$, \emph{i.e.} $(P_{i+1},P_j) \vdash a$. From that fact and easy observations, we get that $$f \in \Hom_{\Lambda'}\left(M_{(P_{i+1},P_j)}, M_a\right) u.$$ 

 We finished to prove that in the first case, we have an almost split sequence.

 (b) Let us consider the second case (the third case is similar to the second case). Let $f: M_{(P_i,*)} \rightarrow M_a$ be a morphism which is not a split monomorphism. As before, $\deg(f) \geq \ell^\theta_{(P_i, *), a} + 2 n \delt{a=(P_i, *)}$.

 Notice that, thanks to the beginning of the proof of Lemma \ref{caldif}, 
 \begin{align*}
  &\ell^\theta_{(P_i, *), (P_{i+1}, P_i)} + \ell^\theta_{(P_{i+1}, P_i), a} - \ell^\theta_{(P_i, *), a} \\ = \,& n \left(\delta_{i = a_1} + 0 + 0 + \left| \delt{i = a_1} - \delt{a_2 \in \ooi{a_1}{i}} \right| - \delt{a_2 \in \ooi{a_1}{i}} \right) \\ = \,& 2n \delt{i = a_1} \delt{a_2 \in \cci{i}{a_1}} = 2n \delt{a=(P_i, *)}
 \end{align*}
 so  $f \in (\Hom_{\Lambda'}(M_{(P_{i+1},P_i)}, M_a) v$.
\end{proof}

\begin{proposition} \label{nonsplit}
 Denoting $v' = u-v$, the non-split extensions between indecomposable objects of $\CM \Lambda'$ are, up to isomorphism,
 $$0 \rightarrow M_{(P_i,P_j)} \xrightarrow{\svecv{u^{d(j,l)}}{u^{d(i,k)}}} M_{(P_i,P_l)} \oplus M_{(P_k,P_j)} \xrightarrow{\svech{u^{d(i,k)}}{-u^{d(j,l)}}} M_{(P_k,P_l)} \rightarrow 0$$
 if $k \in \ooi{i}{j}$ and $l \in \oci{j}{i}$; 
 $$0 \rightarrow M_{(P_i,P_j)} \xrightarrow{\svecv{u^{d(j,k)+n}}{u^{d(i,l)}}} M_{(P_k,P_i)} \oplus M_{(P_l,P_j)} \xrightarrow{\svech{u^{d(i,l)}}{-u^{d(j,k)+n}}} M_{(P_k,P_l)} \rightarrow 0$$
 if $k \in \coi{j}{i}$ and $l \in \ooi{i}{j}$; 
 $$0 \rightarrow M_{(P_i,P_j)} \xrightarrow{\svecv{u^{d(j,k)}}{u^{d(i,l)}}} M_{(P_k,P_i)} \oplus M_{(P_l,P_j)} \xrightarrow{\svech{u^{d(i,l)}}{-u^{d(j,k)}}} M_{(P_k,P_l)} \rightarrow 0$$ $$ 0 \rightarrow M_{(P_i,P_j)} \xrightarrow{\svecv{u^{d(j,l)}}{u^{d(i,k)}}} M_{(P_l,P_i)} \oplus M_{(P_k,P_j)} \xrightarrow{\svech{u^{d(i,k)}}{-u^{d(j,l)}}} M_{(P_k,P_l)} \rightarrow 0 $$
 if $l \in \ooi{i}{k}$ and $j \in \ooi{k}{i}$;
 $$ 0 \rightarrow M_{(P_i,P_j)} \xrightarrow{\svecv{u^{d(i,k)}}{v^{d(j,i)}}} M_{(P_k,P_j)} \oplus M_{(P_i,*)} \xrightarrow{\svech{v^{d(j,k)}}{-v^{2d(i,k)}}} M_{(P_k,*)} \rightarrow 0 $$ $$ 0 \rightarrow M_{(P_i,P_j)} \xrightarrow{\svecv{u^{d(i,k)}}{v'^{d(j,i)}}} M_{(P_k,P_j)} \oplus M_{(P_i,\bowtie)} \xrightarrow{\svech{v'^{d(j,k)}}{-v'^{2d(i,k)}}} M_{(P_k,\bowtie)} \rightarrow 0 $$
 if $k \in \ooi{i}{j}$;
 $$ 0 \rightarrow M_{(P_i,*)} \xrightarrow{\svecv{v^{d(i,k)}}{v^{2d(i,l)}}} M_{(P_k,P_i)} \oplus M_{(P_l,*)} \xrightarrow{\svech{u^{d(i,l)}}{-v^{d(l,k)}}} M_{(P_k,P_l)} \rightarrow 0 $$ $$ 0 \rightarrow M_{(P_i,\bowtie)} \xrightarrow{\svecv{v'^{d(i,k)}}{v'^{2d(i,l)}}} M_{(P_k,P_i)} \oplus M_{(P_l,\bowtie)} \xrightarrow{\svech{u^{d(i,l)}}{-v'^{d(l,k)}}} M_{(P_k,P_l)} \rightarrow 0 $$
 if $i \in \ooi{k}{l}$;
 $$ 0 \rightarrow M_{(P_i,*)} \xrightarrow{v^{d(i,k)}} M_{(P_k,P_i)} \xrightarrow{v'^{d(i,k)}} M_{(P_k,\bowtie)} \rightarrow 0 $$ $$ 0 \rightarrow M_{(P_i,\bowtie)} \xrightarrow{v'^{d(i,k)}} M_{(P_k,P_i)} \xrightarrow{v^{d(i,k)}} M_{(P_k,*)} \rightarrow 0 $$
 if $i \neq k$.

 Moreover, in each case, fixing representatives of these isomorphism classes of short exact sequences induces a basis of the corresponding extension group.
\end{proposition}

\begin{proof}
 First of all, this is easy to check that all these non-split extensions exist (to prove the exactness, the easiest is to project the sequence at each idempotent) and they are non-split and not isomorphic to each other (and therefore linearly independent). Let $i,j,k,l \in \cci{1}{n}$ such that $j \neq i, i+1$ and $l \neq k, k+1$. Thanks to Proposition \ref{ARseq}, we know that $\CM \Lambda'$ admits an Auslander-Reiten duality
 $$\Ext^1_{\Lambda'} (X, Y) \cong \Hom_K (\uHom_{\Lambda'}(Y, \tau(X)), K).$$

 Then, using Proposition \ref{stabmorph}, we get
  \begin{align*}
    \dim \Ext^1_{\Lambda'}\left((M_{(P_k,P_l)}, M_{(P_i,P_j)}\right) &= \dim \uHom_{\Lambda'} \left(M_{(P_i,P_j)}, M_{(P_{k-1}, P_{l-1})}\right) \\
      &= \delt{i -1 \in \ooi{k-1}{l-1}} \delt{l \in \ooi{i}{j}} + \delt{j -1 \in \ooi{k-1}{l-1}} \delt{k \in \ooi{i}{j}} \\
      &= \delt{i \in \ooi{k}{l}} \delt{l \in \ooi{i}{j}} + \delt{j \in \ooi{k}{l}} \delt{k \in \ooi{i}{j}} \\
      &= \begin{cases} 2, & \text{if } l \in \ooi{i}{k} \text{ and } j \in \ooi{k}{i}, \\ 1, & \text{if } k \in \ooi{i}{j} \text{ and } l \in \oci{j}{i}, \\ &\text{or } k \in \coi{j}{i} \text{ and } l \in \ooi{i}{j}, \\ 0, & \text{else.} \end{cases}
  \end{align*}
   We also get
  \begin{align*}
    \dim \Ext^1_{\Lambda'}\left(M_{(P_k,*)}, M_{(P_i,P_j)}\right) &= \dim \uHom_{\Lambda'} \left(M_{(P_i,P_j)}, M_{(P_{k-1}, \bowtie)}\right) \\
      &= \delt{i -1 \in \ooi{k-1}{k-1}} \delt{k \in \ooi{i}{j}}  = \delt{k \in \ooi{i}{j}};
  \end{align*}
  \begin{align*}
    \dim \Ext^1_{\Lambda'}\left(M_{(P_k,P_l)}, M_{(P_i,*)}\right) &= \dim \uHom_{\Lambda'} \left(M_{(P_i,*)}, M_{(P_{k-1}, P_{l-1})}\right) \\
      &= \delt{i -1 \in \ooi{k-1}{l-1}} \delt{l \in \ooi{i}{i}}  = \delt{i \in \ooi{k}{l}};
  \end{align*}
  \begin{align*}
    \dim \Ext^1_{\Lambda'}\left(M_{(P_k,*)}, M_{(P_i,*)}\right) &= \dim \uHom_{\Lambda'} \left(M_{(P_i,*)}, M_{(P_{k-1}, \bowtie)}\right) = 0;
  \end{align*}
  \begin{align*}
    \dim \Ext^1_{\Lambda'}\left(M_{(P_k,*)}, M_{(P_i,\bowtie)}\right) &= \dim \uHom_{\Lambda'} \left(M_{(P_i,\bowtie)}, M_{(P_{k-1}, \bowtie)}\right) \\
      &= \delt{i - 1 \in \ooi{k-1}{k-1}} \delt{k \in \ooi{i}{i}}  = \delt{i \neq k}.
  \end{align*}
  The other cases are realized by swapping $*$ and $\bowtie$. In any case, we exhausted the dimensions with the provided short exact sequences. 
\end{proof}

\begin{corollary} \label{cor:extensioncross}
 If $a$ and $b$ are two tagged arcs of $P^*$, $\dim \Ext^1_{\Lambda'}(M_a, M_b)$ is the minimal number of intersection points between representatives of their isotopy classes (where $(P_i, *)$ and $(P_j, \bowtie)$ are intersecting once for $i \neq j$ by convention).
\end{corollary}

\begin{proof}
 It is an easy case by case study. 
\end{proof}

\subsection{Cluster tilting objects of $\CM \Lambda$ and relation to the cluster category}

Let us recall the definition of cluster tilting objects.
\begin{definition}
Let $\cc$ be a triangulated or exact category. An object $T$ in $\cc$ is said to be \emph{cluster tilting} if $$\add T=\{Z\in \cc \mid \Ext^1_\cc(T,Z)=0 \}=\{Z\in \cc \mid \Ext^1_\cc(Z,T)=0 \},$$ where $\add T$ is the set of finite direct sums of direct summands of $T$.
\end{definition}

For any tagged triangulation $\sigma$ of the once-punctured polygon $P^*$, we denote $T_\sigma = \bigoplus_{a \in \sigma} M_a \cong e_F \Gamma_\sigma$.

\begin{theorem}
  \label{thm:cluster tilting}
 The map $\sigma \to T_\sigma$ gives a one-to-one correspondence between the set of tagged triangulations of $P^*$ and the set of isomorphism classes of basic cluster tilting objects in $\CM \Lambda$. 

 Moreover, for any tagged triangulation $\sigma$, $\End_\Lambda(T_\sigma) \cong \Gamma_\sigma^{\operatorname{op}}$ through right multiplication.
\end{theorem}
\begin{proof}
 Let $E$ be a set of tagged arcs and sides of $P^*$ and $M_E = \bigoplus_{a \in E} M_a$ the corresponding object in $\CM \Lambda$. 
By Corollary \ref{cor:extensioncross}, any two arcs in $E$ are compatible if and only if $\Ext^1_{\Lambda}(M_E,M_E)= 0$. Thus, $M_E$ is cluster tilting if and only if it is a maximal set of compatible tagged arcs and sides of $P^*$ if and only if $E = \sigma$ is a tagged triangulation of $P^*$. Thus, $M_E = T_\sigma$.

For the second part, thanks to Propositions \ref{inductD} and \ref{morph}, for any $a, b \in \sigma$, $\Hom_{\Lambda'} (M_a, M_b) \cong A_{a,b} \cong e_a \Gamma_\sigma e_b$. Therefore,
$$\End_{\Lambda'} (T_\sigma) = \bigoplus_{a,b \in \sigma} e_a \Gamma_\sigma e_b = \Gamma_\sigma.$$
Moreover, the composition on the left coincides with the multiplication on the right by Propositions \ref{inductD} and \ref{morph}. Notice that we get the opposite algebra because we make endomorphism rings act on the left.
\end{proof}  

Theorem \ref{thm:classif} and \ref{thm:cluster tilting} show that the category $\CM \Lambda$ is very similar to the cluster category of type $D_n$. In the rest to this section, we give an explicit connection. First, we recall some basic facts about cluster categories. The cluster category is defined in \cite{BMRRT} as follows.
 \begin{definition}
 For an acyclic quiver $Q$, the \emph{cluster category} $\cc(KQ)$ is the orbit category $\cd^{\operatorname{b}}(KQ)/F$ of the bounded derived category $\cd^{\operatorname{b}}(KQ)$ by the functor $F = \tau^{-1}[1]$, where $\tau$ denotes the Auslander-Reiten translation and $[1]$ denotes the shift functor.
 The objects in $\cc(KQ)$ are the same as in $\cd^{\operatorname{b}}(KQ)$, and the morphisms are given by
 $$\Hom_{\cc(KQ)}(X,Y)={\bigoplus_{i\in \ZZ}}\Hom_{\cd^{\operatorname{b}}(KQ)}(F^iX,Y),$$
 where $X$ and $Y$ are objects in $\cd^{\operatorname{b}}(KQ)$.
 For $f\in \Hom_{\cc(KQ)}(X,Y)$ and $g\in \Hom_{\cc(KQ)}(Y,Z)$, the composition is defined by 
 $$(g\circ f)_i= \sum_{i_1+i_2=i } g_{i_1} \circ F^{i_1}(f_{i_2})$$
 for all $i \in \ZZ$.
 \end{definition}
In \cite{Happel2}, Happel proved that $\cd^{\operatorname{b}}(KQ)$ has Auslander-Reiten triangles.
For a Dynkin quiver $Q$, he shows in \cite{Happel} that the Auslander-Reiten quiver of $\cd^{\operatorname{b}}(KQ)$ is $\ZZ\Delta$, where $\Delta$ is the underlying Dynkin diagram of $Q$.
Then the Auslander-Reiten quiver of $\cc(KQ)$ is $\ZZ\Delta/\varphi$, where $\varphi$ is the graph automorphism induced by $\tau^{-1}[1]$.
In type $D_n$, the Auslander-Reiten quiver of $\cc$ has the shape of a cylinder with $n$ $\tau$-orbits. As a quiver, it is the same as the quiver of $\uCM \Lambda$ (see Figures \ref{CMLe}, \ref{CMLo}).

Recall that a triangulated category is said to be algebraic if it is the stable category of a Frobenius category. Let us state the following result by Keller and Reiten.
\begin{theorem}[{\cite[Introduction and Appendix]{KR}}]
\label{KR}
If $K$ is a perfect field and $\cc$ an algebraic $2$-Calabi-Yau triangulated category containing a cluster tilting object $T$ with $\End_\cc(T)\cong KQ$ hereditary, then there is a triangle-equivalence $\cc(K Q)  \rightarrow \cc$.
\end{theorem}

   By using the above statements, we show the following triangle-equivalences between cluster categories of type $D$ and stable categories of Cohen-Macaulay modules.
\begin{theorem}
\label{th:cyeq}
\begin{enumerate}
\item
\label{2-CY}
The stable category $\uCM \Lambda$ is $2$-Calabi-Yau.
\item
\label{th:equivalence}
If $K$ is perfect, then there is a triangle-equivalence $\cc(KQ) \cong \uCM \Lambda$ for a quiver $Q$ of type $D_n$. 
 \end{enumerate}   
\end{theorem}
\begin{proof}
We will prove (1) in the next subsection independently.

Let $\sigma$ be the triangulation of $P^*$ whose set of tagged arcs is $$\{(P_1,P_3),(P_1,P_4),\cdots,(P_1,P_n), (P_1,*), (P_1, \bowtie )\}.$$

The full subquiver $Q$ of $Q_\sigma$ with set of vertices $Q_{\sigma,0} \smallsetminus F$ is a quiver of type $D_n$.
Thus, we have  $$\Gamma_\sigma^{\operatorname{op}}/ (e_F) \cong (KQ)^{\operatorname{op}}.$$

By Theorem \ref{thm:cluster tilting}, for the cluster tilting object $T_\sigma$, we have the following isomorphism 
\begin{align*}
\uEnd_{\Lambda}(T_\sigma) \cong \Gamma_\sigma^{\operatorname{op}}/ (e_F) .
\end{align*} 

Then, by Theorem \ref{KR}, we have $\cc\left((KQ)^{\operatorname{op}}\right) \cong \uCM \Lambda$.
\end{proof}

\subsection{Proof of Theorem \ref{th:cyeq} (\ref{2-CY})}
\label{s:2-CY}

In this subsection, we prove that the stable category $\uCM \Lambda$ is $2$-Calabi-Yau. Throughout, we denote $\DK:=\Hom_K(-,K)$, $\DR:=\Hom_{{R}}(-,{R})$ and $(-)^*:=\Hom_{\Lambda}(-,\Lambda)$.

Let us recall some general definitions and facts about Cohen-Macaulay modules. Let $A$ be an $R$-order. 

\begin{definition}
\label{def:gorenstein}
 We say that $X$ is an \emph{injective} {Cohen-Macaulay} $A$-module if $\Ext^1_A(Y,X)=0$ for any $Y \in \CM A$, or equivalently, $X \in \add(\Hom_R(A^{\text{op}}, R))$.
 Denote by $\inj A$ the category of {\emph{injective}} Cohen-Macaulay $A$-modules.
 
 An $R$-order $A$ is \emph{Gorenstein} if $\Hom_{{R}}(A_A, R)$ is projective as a left $A$-module, or equivalently, if $\Hom_R({}_A A, R)$ is projective as a right $A$-module. 
 We have an exact duality $\DR : \CM A^{\operatorname{op}} \to \CM A$.
 \end{definition}

 The Nakayama functor is defined here by $\nu: \proj A \xlongrightarrow{(-)^*} \proj A^{\operatorname{op}} \xlongrightarrow{\DR} \inj A$, which is isomorphic to $(\DR A) \otimes_{A}-$.
For any Cohen-Macaulay $A$-module $X$, consider a projective presentation $$P_1\xrightarrow{f} P_0 \to X \to 0$$
 and apply $(-)^* : \mod A \to \mod A^{\operatorname{op}}$ to it to get the following exact sequence: $$0\lra X^*\lra P_0^*\xlongrightarrow{f^*} P_1^*\lra \coker(f^*) \lra0.$$
 We denote $\coker(f^*)$ by $\Tr X$ and we get $\im(f^*)=  \Omega \Tr X$, where $\Omega$ is the syzygy functor: $\umod A^{\operatorname{op}} \to \umod A^{\operatorname{op}}$.
Then we apply $\DR: \CM A^{\operatorname{op}} \to \CM A$ to $$0\lra X^*\lra P_0^* \xlongrightarrow{f^*} \Omega \Tr X \lra 0$$ and denote $\tau X:=\DR \Omega \Tr X$, we get the exact sequence
\begin{equation}
\label{eq:tau}
  0\lra \tau X \lra \nu P_0 \lra \nu X \lra0.
  \end{equation}
 
  For an $R$-order $A$, if $K(x) \otimes_R A$ is a semisimple $K(x)$-algebra, then we call $A$ an isolated singularity.
  By using the notions above, we have the following well-known results in Auslander-Reiten theory.
\begin{theorem}\cite{FMO,almostsplitseqorder,MR0412223}
 \label{isolated singularity th}
Let $A$ be an $R$-order. If $A$ is an isolated singularity, then 
  \begin{enumerate}
     \item  \cite[Chapter \uppercase\expandafter{\romannumeral1}, Proposition 8.3]{FMO} The construction $\tau$ gives an equivalence $\uCM A \to \oCM A$, where $\oCM A$ is the quotient of $\CM A$ with respect to the subgroup of maps which factor through an injective object.
     \item   \cite[Chapter \uppercase\expandafter{\romannumeral1}, Proposition 8.7]{FMO} For $X,Y\in \uCM A$, there is a functorial isomorphism $$\uHom_{A}(X,Y) \cong \DK\Ext^1_{A}(Y, \tau X).$$
  \end{enumerate}
\end{theorem}

For Gorenstein orders, we have the following nice properties.
\begin{proposition}
 \label{gorenstein prop}
 Assume that $A$ is a Gorenstein isolated singularity, then we have
  \begin{enumerate}
    \item $\CM A$ is a Frobenius category.
    \item $\uCM A$ is a $K$-linear Hom-finite triangulated category.
    \item $\tau=\Omega \nu=[-1]\circ \nu$.
  \end{enumerate}
\end{proposition}
\begin{proof}
(1) The projective objects in $\CM A$ are just projective $A$-modules. They are also injective objects. Since each finitely generated $A$-module is a quotient of a projective $\Lambda$-module, it follows that $\CM A$ is a Frobenius category; (2) is due to \cite{Happel} and \cite[Lemma 3.3]{CM}; (3) is a direct consequence of  \eqref{eq:tau}.
\end{proof}

The order $\Lambda$ is Gorenstein. Indeed, as a graded left $\Lambda$-module,
$$\DR(\Lambda_\Lambda) = \Hom_R\left(
   \begin{sbmatrix}
    R' & R'  & R' & \cdots & R' & R' & X^{-1}(X, Y) \\
    (X,Y) & R' & R' & \cdots & R' & R' & R' \\
    (X) & (X,Y) & R' & \cdots& R' & R' & R' \\
    \vdots&\vdots & \vdots &\ddots&\vdots&\vdots&\vdots\\
    (X) & (X) & (X) & \cdots & R'  & R' & R' \\
    (X) & (X) & (X) & \cdots & (X,Y) & R' & R' \\
    (X) & (X) & (X) & \cdots & (X) & (X,Y) & R'
   \end{sbmatrix}
   , R \right)$$
can be identified with
\begin{align*} & X^{-1} \begin{sbmatrix}
    R' & X^{-1}(X, Y)  & X^{-1} R' & \cdots  &X^{-1} R' & X^{-1} R'& X^{-1} R' \\
    R' & R' & X^{-1}(X,Y)  & \cdots& X^{-1} R' &X^{-1} R'& X^{-1} R'\\
    R' & R' & R' & \cdots& X^{-1} R' &X^{-1} R'& X^{-1} R'\\
    \vdots&\vdots & \vdots &\ddots&\vdots&\vdots&\vdots\\
    R' & R' & R' & \cdots& R'  &X^{-1}(X,Y) &X^{-1} R' \\
    R'&R' & R' &\cdots& R'& R'&  X^{-1}(X,Y) \\
    (X,Y) & R' &R' & \cdots& R' & R'& R'
   \end{sbmatrix} \\ =\, & \Lambda V^{-1} \subset \M_n(R'[X^{-1}]),\end{align*}
 where
  $$V=
   \begin{sbmatrix}
  0           & 0                   & \ldots     & 0              & X              & 0   \\
  0           & 0                      & \ldots     & 0              & 0              & X   \\
   X^2      & 0                    & \ldots     & 0              & 0              & 0    \\
   0          & X^2                 & \ldots     & 0              & 0              & 0   \\
  \vdots   & \vdots       & \ddots    & \vdots      & \vdots     & \vdots \\
  0           & 0                    & \ldots      & X^2          & 0              & 0
  \end{sbmatrix}.$$
  Therefore $\DR(\Lambda_\Lambda)$ is a projective (left) $\Lambda$-module.

According to Theorem \ref{isolated singularity th} and Proposition \ref{gorenstein prop},
we have $$\uHom_{\Lambda}(X,Y) \cong \DK\uHom_{\Lambda}(Y, \nu X)$$ for $X,Y\in \CM \Lambda.$
Thus $\nu= (\DR\Lambda)\otimes_{\Lambda}-$ is a Serre functor.
 We want to prove that $$(\DR\Lambda)\otimes_{\Lambda}- \cong  \Omega^{-2}(-).$$
Thanks to the previous discussion, there is an isomorphism of $\Lambda$-modules:
  \[ f: \Lambda \to \DR(\Lambda_\Lambda) \text{, }  \mu \mapsto \mu V^{-1}. \]

  We define the automorphism $\alpha$ of $\Lambda$ by $\alpha (\lambda)= V^{-1} \lambda V$ for $\lambda \in \Lambda$. 
       The automorphism $\alpha$ corresponds to a $4\pi/n$ counter-clockwise rotation of the quiver of $\Lambda$ {shown} page \pageref{remark:quiver}.
    In fact, if $$\lambda=
  \begin{sbmatrix}
  \lambda_{1,1}             & \lambda_{1,2}           & \ldots     & \lambda_{1,n-2}        &\lambda_{1,n-1}      & \lambda_{1,n} \\
  \lambda_{2,1}          & \lambda_{2,2}           & \ldots     & \lambda_{2,n-2}        &\lambda_{2,n-1}      & \lambda_{2,n}   \\
  \lambda_{3,1}       &  \lambda_{3,2}         & \ldots     & \lambda_{3,n-2}        &\lambda_{3,n-1}      & \lambda_{3,n}   \\
  \vdots                        & \vdots                       & {\ddots}     & \vdots                      & \vdots                    & \vdots     \\
  \lambda_{n-2,1}   & \lambda_{1,2}      & \ldots     & \lambda_{n-2,n-2}     &\lambda_{n-2,n-1}   & \lambda_{n-2,n}  \\
  \lambda_{n-1,1}    & \lambda_{n-1,2}   & \ldots     & \lambda_{n-1,n-2}   &\lambda_{n-1,n-1}   & \lambda_{n-1,n}   \\
  \lambda_{n,1}       & \lambda_{n,2}      & \ldots     & \lambda_{n,n-2}   & \lambda_{n,n-1}  & \lambda_{n,n}    
  \end{sbmatrix}$$ is an element in $\Lambda$, 
    then
  $$\alpha(\lambda) =
    \begin{sbmatrix}
  \lambda_{3,3}             & \lambda_{3,4}           & \ldots     & \lambda_{3,n}        &X^{-1} \lambda_{3,1}      & X^{-1} \lambda_{3,2} \\
  \lambda_{4,3}           & \lambda_{4,4}           & \ldots     & \lambda_{4,n}        &X^{-1}\lambda_{4,1}      & X^{-1} \lambda_{4,2}   \\
  \lambda_{5,3}       & \lambda_{5,4}          & \ldots     & \lambda_{5,n}        &X^{-1}\lambda_{5,1}      & X^{-1} \lambda_{5,2}   \\
  \vdots                        & \vdots                       & {\ddots}     & \vdots                   & \vdots                 & \vdots     \\
   \lambda_{n,3}      & \lambda_{n,4}      & \ldots     & \lambda_{n,n}        & X^{-1}\lambda_{n,1}       & X^{-1}\lambda_{n,2}  \\
  X \lambda_{1,3}       & X \lambda_{1,4}      & \ldots     & X \lambda_{1,n}       &\lambda_{1,1}       & \lambda_{1,2}   \\
 X \lambda_{2,3}       & X \lambda_{2,4}      & \ldots     & X \lambda_{2,n}   & \lambda_{2,1}     & \lambda_{2,2}    
  \end{sbmatrix}.$$ 

Let $A$ and $B$ be two $R$-orders. We define ${}_\vartheta M_{\varsigma}$ for an $(A,B)$-bimodule $M$, $\vartheta \in \Aut(A)$ and $\varsigma \in \Aut(B)$ as follows:
${}_\vartheta M_{\varsigma}:=M$ as a vector space and the $(A,B)$-bimodule structure is given by $$a \times  m \times b = \vartheta(a) m \varsigma(b)$$ for $m \in{}_\vartheta M_{\varsigma}$ and $a\in A$, $b \in B$.
Since $\vartheta \in \Aut(A)$, ${}_\vartheta (-)$ is an automorphism of $\mod A$.

 \begin{proposition}
  \label{pro:bimodule}
    The above $f: \Lambda \to \DR \Lambda$ gives an isomorphism of $\Lambda$-bimodules $${}_1\Lambda_{\alpha} \cong \DR\Lambda.$$    \end{proposition}
    \begin{proof}
    Clearly, $f$ preserves the left action of $\Lambda$.
    Moreover, it preserves the right action since for $\lambda$, $ \mu \in \Lambda$, we have \begin{align*}&f(\mu \alpha(\lambda))= f(\mu (V^{-1} \lambda V) )=\mu (V^{-1} \lambda V)  V^{-1}= \mu V^{-1} \lambda=f(\mu) \lambda.\qedhere \end{align*}
  \end{proof}
  
  By using the isomorphism of Proposition \ref{pro:bimodule}, we find the following description of  the Nakayama functor $\nu$.
  \begin{lemma}
  \label{lem:functoriso}
  We have an isomorphism  $\nu \cong {}_{\alpha^{-1}} (-)$ of endofunctors of $\uCM \Lambda$.
  \end{lemma}
  \begin{proof}
  Since $\DR\Lambda \cong  {}_1\Lambda_{\alpha} $, 
  it follows that $\nu  \cong {}_1\Lambda_{\alpha} \otimes_{\Lambda}-$.
On the other hand, we have an isomorphism $H:  {}_1\Lambda_{\alpha} \otimes_{\Lambda}- \cong {}_{\alpha^{-1}} (-)$ given by $\lambda \otimes - \mapsto \alpha^{-1}(\lambda)(-)$.
  Thus the assertion follows.
  \end{proof}
   
   Let $T=K[x,y]$ and $S:=K[x,y]/(p)$ for some $p \in T$.

We define a $\ZZ/n\ZZ$-grading on $T$ by setting $\text{deg}(x)=1$ and $\text{deg}(y)=-1$. 
   This makes $T$  a $\ZZ/n\ZZ$-graded algebra
   \[
   T=\bigoplus_{\overline{i}\in \ZZ/n\ZZ}T_{\overline{i}}=T_{\overline{0}} \oplus T_{\overline{1}} \oplus \ldots \oplus T_{\overline{n-1}}.
   \]
    Suppose that $p$ is homogeneous of degree $d$ with respect to this grading. Then the quotient ring $$S=K[x,y]/(p)=S_{\overline{0}} \oplus S_{\overline{1}} \oplus \ldots \oplus S_{\overline{n-1}}$$ has a natural structure of a $\ZZ/n\ZZ$-graded algebra.
    The following result can be easily established from classical results about matrix factorization (see \cite[Theorem 3.22]{DeLu} for a detailed proof). 

\begin{theorem}[\cite{CM}]
\label{degree}
In the category $\uCM^{\ZZ/n\ZZ}(S)$, there is an isomorphism of autoequivalences $[2] \cong (-d)$.
\end{theorem}
   Setting $p:=x^{n-1}y-y^2$, we have $S=K[x,y]/(x^{n-1}y-y^2)$. 
   Identifying $R'=K[X,Y]/(XY-Y^2)$ as a subalgebra of $S$ via $X \mapsto x^n$ and $Y \mapsto xy$, we regard $S$ as an $R'$-algebra. We obtain the following lemma.
   
\begin{lemma}
\label{lem:matrixS}
For $i \in \cci{0}{n-1}$, as $R'$-modules, we have $$S_{\overline{i}}
\cong\begin{cases}
  R' x^i, &\text{if } i\in \cci{0}{n-2}, \\
  (1, X^{-1}Y) x^{n-1}, &\text{if } i=n-1.
\end{cases}
$$
\end{lemma}
\begin{proof}
 Let $i \in \cci{0}{n-2}$. Over $R=K[X]$, $S_{\overline{i}}$ is generated by $x^{i}$ and $x^{i+1}y$. Thus, we have $S_{\overline{i}}\cong  R' x^i $. Over $R$, $S_{\overline{n-1}}$ is generated by $x^{n-1}$ and $y$. So $S_{\overline{n-1}} \cong (1, X^{-1}Y) x^{n-1}$.
 \end{proof}
 We define an $R$-order $S^{[n]}$ which is a subalgebra of $\M_n(S)$ from the $\ZZ/n\ZZ$-graded algebra $S$ as follows:  
 $$S^{[n]}=\begin{sbmatrix}
      S_{\overline{0}} & S_{\overline{1}}  & S_{\overline{2}}    &\cdots & S_{\overline{n-2}}    & S_{\overline{n-1}}\\
      S_{\overline{n-1}} & S_{\overline{0}}  & S_{\overline{1}}    &\cdots & S_{\overline{n-3}}    & S_{\overline{n-2}}\\
      S_{\overline{n-2}} & S_{\overline{n-1}}  & S_{\overline{0}}    &\cdots & S_{\overline{n-4}}    & S_{\overline{n-3}}\\
      \vdots           & \vdots            & \vdots              &\ddots & \vdots              & \vdots\\
      S_{\overline{2}} & S_{\overline{3}}  & S_{\overline{4}}    &\cdots & S_{\overline{0}}    & S_{\overline{1}}\\
      S_{\overline{1}} & S_{\overline{2}}  & S_{\overline{3}}    &\cdots & S_{\overline{n-1}}    & S_{\overline{0}}
  \end{sbmatrix}.$$
  
    \begin{proposition}
    \label{pro:iso}
We have an isomorphism $S^{[n]} \cong \Lambda$ of $R'$-algebras.
  \end{proposition}
  \begin{proof}
According to Lemma \ref{lem:matrixS}, we have
   $$S^{[n]}=\begin{sbmatrix}
      R' & R'{x}  & R'{x}^2    &\cdots & R'{x}^{n-2}   & (1, X^{-1}Y){x}^{n-1}\\
      (1, X^{-1}Y){x}^{n-1} & R'  & R'{x}    &\cdots & R'{x}^{n-3}    & R'{x}^{n-2}\\
      R'{x}^{n-2} & (1, X^{-1}Y){x}^{n-1}  & R'    &\cdots & R'{x}^{n-4}    & R'{x}^{n-3}\\
      \vdots           & \vdots            & \vdots              &\ddots & \vdots              & \vdots\\
      R'{x}^{2} & R'{x}^{3}  & R'{x}^4    &\cdots & R'  & R'{x}\\
      R'{x} & R'{x}^{2}  &R'{x}^{3}    &\cdots & (1, X^{-1}Y){x}^{n-1}    & R'
  \end{sbmatrix}.$$
  Taking the conjugation by the diagonal matrix $B = \diag(x^i)_{i \in \cci{0}{n}}$, we get $$B S^{[n]} B^{-1} = \Lambda.$$
  \end{proof}
From now on, we identify $\Lambda$ and $S^{[n]}$.  
 Consider the matrix 
 $$U=
   \begin{sbmatrix}
  0           & 0                    & \ldots     & 0              & 1              & 0   \\
  0           & 0                    & \ldots     & 0              & 0              & 1   \\
   1      & 0                      & \ldots     & 0              & 0              & 0    \\
   0          & 1                   & \ldots     & 0              & 0              & 0   \\
  \vdots   & \vdots      & \ddots    & \vdots      & \vdots     & \vdots \\
  0           & 0           & \ldots                 & 1         & 0              & 0
  \end{sbmatrix}.$$
The automorphism $\beta$ of $S^{[n]}$ which is given by $\beta(s)=U^{-1}s U$ for $s \in S^{[n]}$ corresponds to the automorphism $\alpha$ of $\Lambda$.
  Thus we have an isomorphism ${}_1S^{[n]}_{\beta^{-1}} \cong  {}_1 \Lambda_{\alpha^{-1}}$ of $S^{[n]}$-bimodules.

Using the same notation above, we have the following lemma.
\begin{lemma}
\begin{enumerate}
 \item \cite[Theorem 3.1]{IyamaLerner}
 The functor 
\begin{align*}
F:\mod^{\ZZ/n\ZZ}S &\to \mod S^{[n]}\\
M_{\overline{0}} \oplus M_{\overline{1}} \oplus \ldots \oplus M_{\overline{n-1}} 
& \mapsto 
  \begin{sbmatrix}
  M_{\overline{0}} & M_{\overline{1}} & \ldots & M_{\overline{n-1}}
  \end{sbmatrix}^{\operatorname{t}}
\end{align*} is an equivalence of categories.
\item For $i\in \ZZ$, we denote by
$(i):\mod^{\ZZ/n\ZZ}S \to \mod^{\ZZ/n\ZZ}S$
the grade shift functor defined by $M(i)_{\overline{j}}:=M_{\overline{i+j}}$ for $M \in \mod^{\ZZ/n\ZZ}S$.
The grade shift functor $(i)$ induces an autofunctor (denoted by $\gamma_i$) in $\mod S^{[n]}$ which makes the following diagram commute:
\[
  \begin{tikzcd}
  \mod^{\ZZ/n\ZZ}S \arrow{r}{F} \arrow{d}[swap]{(i)} &\mod S^{[n]} \arrow{d}{\gamma_i}\\
  \mod^{\ZZ/n\ZZ}S \arrow{r}{F} & \mod S^{[n]}.
  \end{tikzcd}
\]
\end{enumerate}
\end{lemma}

More precisely, for any left $S^{[n]}$-module $\begin{sbmatrix}
  M_{\overline{0}} & M_{\overline{1}} & \ldots & M_{\overline{n-1}}
  \end{sbmatrix}^{\operatorname{t}}$, we have $$\gamma_i \left(\begin{sbmatrix}
 M_{\overline{0}} & M_{\overline{1}} & \ldots & M_{\overline{n-1}}
  \end{sbmatrix}^{\operatorname{t}}\right)=\begin{sbmatrix}
  M_{\overline{i}} & M_{\overline{i+1}} & \ldots & M_{\overline{i+n-1}}
  \end{sbmatrix}^{\operatorname{t}}.$$
  
Now we can prove the $2$-Calabi-Yau property of $\uCM \Lambda$.
\begin{proof}[Proof of Theorem \ref{th:cyeq} \eqref{2-CY}]
    The equivalence $\mod^{\ZZ/n\ZZ}S \cong \mod S^{[n]} = \mod\Lambda$ induces an equivalence $$\CM^{\ZZ/n\ZZ}(S) \cong \CM S^{[n]} = \CM \Lambda.$$
     In the category $\uCM^{\ZZ/n\ZZ} S$, according to Theorem \ref{degree}, we have an isomorphism of functors $$[2] \cong (-\text{deg }(x^{n-2}-y^2))=(2).$$
     By Lemma \ref{lem:functoriso}, we have $\nu \cong {}_{\alpha^{-1}}(-)$. Therefore, it is enough to prove ${}_{\alpha^{-1}} (-) \cong (2)$.
     This is equivalent to prove that ${}_\beta M \cong \gamma_{-2}(M)$ holds for any $M \in \CM S^{[n]}$.
    
 Let $s_i$ be the row matrix which has $1$ in the $i$-th column.  Since     
   \begin{align*}
     {}_\beta M & =     \begin{sbmatrix}
       s_{0} \times {}_\beta M & s_{1}\times {}_\beta M & s_2 \times {}_\beta M & \ldots & s_{n-1} \times{{}_\beta M}
      \end{sbmatrix}^{\operatorname{t}}\\
     & =
     \begin{sbmatrix}
     \beta(s_{0})M & \beta(s_{1})M & \beta(s_{2})M & \ldots & \beta(s_{n-1})M
     \end{sbmatrix}^{\operatorname{t}}\\
     & =
     \begin{sbmatrix}
     s_{n-2}M & s_{n-1}M & s_{0} M & \ldots & s_{n-3}M
     \end{sbmatrix}^{\operatorname{t}},
   \end{align*}
   it follows that ${}_\beta M \cong \gamma_{-2}(M)$. Therefore, the category $\uCM \Lambda$ is $2$-Calabi-Yau.
  \end{proof}

\section{Graded Cohen-Macaulay $\Lambda$-modules}
\label{s:graded}
In this section, we prove a graded version of Theorem \ref{th:cyeq} which gives a relationship between the category $\CM^{\ZZ}(\Lambda)$ of graded Cohen-Macaulay $\Lambda$-modules and the bounded derived category $\cd^{\operatorname{b}}(KQ)$ of type $D_n$.

Let $Q$ be an acyclic quiver. 
We denote by $\ck^{\operatorname{b}}(\proj KQ)$ the bounded homotopy category of finitely generated projective $KQ$-modules, and by $\cd^{\operatorname{b}}(KQ)$ the bounded derived category of finitely generated $KQ$-modules.
These are triangulated categories and the canonical embedding 
$\ck^{\operatorname{b}}(\proj KQ) \to \cd^{\operatorname{b}}(KQ)$ is a triangle functor.

We define a grading on $\Lambda$ by $\Lambda_i = \Lambda \cap \M_n(K X^i + K X^{i-1} Y)$ for $i \in \ZZ$.
 This makes $\Lambda=\bigoplus_{i\in \ZZ}\Lambda_i$ a $\ZZ$-graded algebra.
The category of graded Cohen-Macaulay $\Lambda$-modules, $\CM^{\ZZ}(\Lambda)$, is defined as follows.
The objects are graded $\Lambda$-modules which are Cohen-Macaulay,
and the morphisms in $\CM^{\ZZ}(\Lambda)$ are $\Lambda$-morphisms preserving the degree.
The category $\CM^{\ZZ}(\Lambda)$ is a Frobenius category. Its stable category is denoted by $\uCM^{\ZZ}(\Lambda)$.
For $i\in \ZZ$, we denote by
$(i):\CM^{\ZZ}(\Lambda) \to \CM^{\ZZ}(\Lambda)$
the grade shift functor: Given a graded Cohen-Macaulay $\Lambda$-module $X$, we define $X(i)$ to be $X$ as a $\Lambda$-module, with the grading $X(i)_j = X_{i+j}$ for any $j \in \ZZ$. 

\begin{remark}
We show that this grading of $\Lambda$ is analogous to the grading of $\Lambda'$ given by the $\theta$-length.
 Let $i,j\in F$. By Theorem \ref{cor:lambdamod}, $e_i \Lambda' e_j \cong e_i \Lambda e_j$ holds.
Let $\lambda \in e_i \Lambda' e_j \cong e_i \Lambda e_j$. Using a similar argument as in the proof of Theorem \ref{cor:lambdamod}, 
we have $${\frac{\ell^{\theta}(\lambda)+2 d(1,i) - 2d(1,j) }{2n}}={\deg}(\lambda)$$
where $\deg(\lambda)$ is the degree of $\lambda$ as a member of $\Lambda$. Consider the two graded algebras $$\Lambda' =\bigoplus_{i=1}^n \Lambda' e_i \and \Lambda'' := \End\left(\bigoplus_{i=1}^n u^{2 d(1,i)}\Lambda' e_i  \right).$$ By graded Morita equivalence, we have $\CM^\ZZ \Lambda'\cong \CM^\ZZ(\Lambda'')$.
Since $\Lambda \cong \Lambda''$ as $R$-orders and $\deg(X)=2n$ in $\Lambda''$, it follows that the Auslander-Reiten quiver of $\CM^\ZZ(\Lambda'')$ has $2n$ connected components each of which is a degree shift of the Auslander-Reiten quiver of $\CM^\ZZ \Lambda$.
\end{remark}

We introduce the properties of $ \CM^\ZZ \Lambda$ in the following theorems.

\begin{theorem}
\label{th:graded modules}
\begin{enumerate}
   \item The set of isomorphism classes of indecomposable objects of $\CM^\ZZ \Lambda$ is $$\{(i,j) \mid i,j \in \ZZ, 0 < j-i < n\} \cup \{(i, *) \mid i \in \ZZ\} \cup \{(i,\bowtie)\mid i\in \ZZ\},$$ where
  \begin{align*} 
   \\
   (i,j) &= \tr{\begin{bmatrix} \bovermat{i}{(X) \cdots (X)} & \bovermat{j-i}{(X^2,Y^2) \cdots (X^2, Y^2)} & \bovermat{n-j}{(X^2) \cdots (X^2)} \end{bmatrix}} \quad \text{if } 0 < i < j \leq n; \\ \\
   (i,j) &=  \tr{\begin{bmatrix} \bovermat{j-n}{(X,Y) \cdots (X,Y)} & \bovermat{n-j+i}{(X) \cdots (X)} & \bovermat{n-i}{(X^2,Y^2) \cdots (X^2, Y^2)} \end{bmatrix}} \quad \text{if } i \leq n < j;\\ \\
   (i,*) &= \tr{\begin{bmatrix} \bovermat{i}{(Y) \cdots (Y)} & \bovermat{n-i}{(Y^2) \cdots (Y^2)} \end{bmatrix} } \quad \text{if } 0 < i \leq n;\\ \\
   (i,\bowtie)&= \tr{\begin{bmatrix} \bovermat{i}{(X-Y) \cdots (X-Y)} & \bovermat{n-i}{(X^2-Y^2) \cdots (X^2-Y^2)} \end{bmatrix} }  \quad \text{if } 0 < i \leq n,
  \end{align*}
  and the other $(i,j)$ are obtained by shift: \begin{align*} (i+kn,j+kn)&=(i,j)(k) \\ (i+kn,*)&=(i,*)(k)\\ (i+kn,\bowtie)&=(i,\bowtie)(k)\end{align*} for $k\in \ZZ$. The projective-injective objects are of the form $(i, i+1)$ for $i \in \ZZ$.
 \item The non-split extensions of indecomposable objects of $\CM^\ZZ \Lambda$ are of the form 
 $$0 \rightarrow (i,j) \rightarrow (i,l) \oplus (k,j) \rightarrow (k,l) \rightarrow 0 \quad \text{if } i<k<j<l<i+n;$$ 
 $$0 \rightarrow (i,j) \rightarrow (k,i+n) \oplus (l,j) \rightarrow (k,l+n) \rightarrow 0 \quad \text{if } i < l < j \leq k < i+n;$$ 
 $$\left. \begin{matrix} 0 \rightarrow (i,j) \rightarrow (k,i+n) \oplus (l,j) \rightarrow (k,l+n) \rightarrow 0 \\ 0 \rightarrow (i,j) \rightarrow (l,i+n) \oplus (k,j) \rightarrow (k,l+n) \rightarrow 0 \end{matrix} \right\} \quad \text{if } i < l < k < j < i+n;$$
 $$\left. \begin{matrix}  0 \rightarrow (i,j) \rightarrow (k,j) \oplus (i,*) \rightarrow (k,*) \rightarrow 0 \\ 0 \rightarrow (i,j) \rightarrow (k,j) \oplus (i,\bowtie) \rightarrow (k,\bowtie) \rightarrow 0 \end{matrix} \right\} \quad \text{if } i < k < j < i+n;$$
 $$\left. \begin{matrix}  0 \rightarrow (i,*) \rightarrow (k,i+n) \oplus (l,*) \rightarrow (k,l+n) \rightarrow 0 \\ 0 \rightarrow (i,\bowtie) \rightarrow (k,i+n) \oplus (l,\bowtie) \rightarrow (k,l+n) \rightarrow 0 \end{matrix} \right\} \quad \text{if } i < l < k < i+n;$$
 $$\left. \begin{matrix}  0 \rightarrow (i,*) \rightarrow (k,i+n) \rightarrow (k,\bowtie) \rightarrow 0 \\ 0 \rightarrow (i,\bowtie) \rightarrow (k,i+n) \rightarrow (k,*) \rightarrow 0 \end{matrix} \right\} \quad \text{if } i < k < i+n.$$

 Moreover, fixing representatives of these isomorphism classes of short exact sequences induces bases of the corresponding extension groups.
 \item The exact category $\CM^\ZZ \Lambda$ admits the Auslander-Reiten sequences
  $$0\lra (i,j) \lra (i,j+1)\oplus(i+1,j) \lra (i+1,j+1)\lra 0$$
  for $i+1 < j < i+n$ (with the convention that $(i, i+n) = (i,*) \oplus (i, \bowtie)$);
  \begin{align*} &0\lra (i,*) \lra (i+1,i+n) \lra (i+1,\bowtie)\lra 0 \\ \text{and} \quad & 0\lra (i,\bowtie) \lra (i+1,i+n) \lra (i+1,*)\lra 0\end{align*}
  for any $i \in \ZZ$.
 \item The Auslander-Reiten quiver of $\CM^{\ZZ}(\Lambda)$ is the repetitive quiver of type $D_{n+1}$ (unfolded version of Figures \ref{CMLe} and \ref{CMLo}).
 \item The syzygy in $\uCM^{\ZZ}(\Lambda)$ is defined on indecomposable objects by:
  $$\Omega((i,j)) = (i+1-n, j+1-n) ;  \quad \Omega((i,*)) = (i+1-n, \bowtie) ;  \quad \Omega((i,\bowtie)) = (i+1-n, *).$$
 \end{enumerate}
\end{theorem}

\begin{proof}
(1) First of all, it is immediate that the graded modules $(i,j)$ for $0 < j-i < n$, $(i, *)$ and $(i,\bowtie)$ for $i \in \ZZ$ are not isomorphic. Therefore, we need to prove that there are no other isomorphism classes. We consider the degree forgetful functor $F: \CM^\ZZ \Lambda \to \CM \Lambda$. Let $X \in \CM^\ZZ \Lambda$ be indecomposable and $M$ be an indecomposable direct summand of $FX$ in $\CM \Lambda$. Using Theorem \ref{thm:classif}, there exists a tagged arc or a side $a$ of $P^*$ such that $M \cong M_a$. Then, it is immediate that $M \cong FY$ where $Y \in \CM^\ZZ \Lambda$ is $(i,j)$ for $0 < j-i < n$ or $(i, *)$ or $(i,\bowtie)$ for $i \in \ZZ$. There are two morphisms $f: FY \rightarrow F X$ and $g: F X \rightarrow F Y$ such that $gf = \Id_{F Y}$. Let us write 
$$f = \sum_{m \in \ZZ} f_m \and g = \sum_{m \in \ZZ} g_m$$
where $f_m$ is a graded morphism from $Y$ to $X(m)$ and $g_m$ a graded morphism from $X(m)$ to $Y$. Thus, we have
$$\sum_{k \in \ZZ} g_{k} f_{k} = \Id_{F Y}$$
and, as the endomorphism ring of $Y$ is $K$, there exists $k \in \ZZ$ such that $g_k f_{k}$ is a non-zero multiple of $\Id_Y$. In other terms we found two graded morphisms $\tilde f: Y \rightarrow X(k)$ and $\tilde g: X(k) \rightarrow Y$ such that $g f = \Id_Y$. Thus, in $\mod^\ZZ(\Lambda)$, we have an isomorphism $X \cong Y(-k) \oplus X'$. As $\mod^\ZZ(R)$ is Krull-Schmidt, $X'$ is necessarily a graded Cohen-Macaulay module. Finally, as $X$ is indecomposable in $\CM^\ZZ \Lambda$, we get that $X \cong Y(-k)$. 

Thus, the set of isomorphism classes of indecomposable graded Cohen-Macaulay $\Lambda$-modules is 
$$\{(i,j) \mid i,j \in \ZZ, 0 < j-i < n\} \cup \{(i, *) \mid i \in \ZZ\} \cup \{(i,\bowtie)\mid i\in \ZZ\}.$$

Statements (2) and (3) are direct consequences through $F$ of the ungraded versions of Propositions \ref{ARseq} and \ref{nonsplit}. The statement (4) is a direct consequence of (1) and (3). 

 For (5), using (2), the short exact sequences constructed from projective covers are:
 $$0 \rightarrow (i+1-n, j+1-n) \rightarrow (i,i+1) \oplus (j-n, j-n+1) \rightarrow (i,j) \rightarrow 0;$$
 $$0 \rightarrow (i+1-n, \bowtie) \rightarrow (i,i+1) \rightarrow (i,*) \rightarrow 0;$$
 \begin{align*} &0 \rightarrow (i+1-n, *) \rightarrow (i,i+1) \rightarrow (i,\bowtie) \rightarrow 0. \qedhere\end{align*}
\end{proof}

For any indecomposable graded Cohen-Macaulay $\Lambda$-module $A$, if $A$ is of the form $(i,j)$ for two integers $i$ and $j$, we write $A_1 = i$ and $A_2 =j$, and if $A$ is of the form $(i,*)$ or $(i,\bowtie)$, we write $A_1 = i$ and $A_2 = i+n$. In this way, all morphisms in $\CM^\ZZ \Lambda$ are going in the increasing direction in term of these pairs of integers.

\begin{definition}
Let $\cc$ be a triangulated category. An object $T$ is said to be \emph{tilting} if $\Hom_\cc(T, T[k])=0$ for any $k \neq 0$ and $\thick(T)=\cc$, where $\thick(T)$ is the smallest full triangulated subcategory of $\cc$ containing $T$ and closed under isomorphisms and direct summands.
\end{definition}

\begin{theorem}[{\cite[Theorem 4.3]{derdgcat}, \cite[Theorem 2.2]{IyTa13}, \cite{BK1991}}]
\label{BK}
Let $\cc$ be an algebraic triangulated Krull-Schmidt category. If $\cc$ has a tilting object $T$, then there exists a triangle-equivalence
$$\cc \rightarrow \ck^{\operatorname{b}}(\proj \End_{\cc}(T)).$$
\end{theorem}

Now we get the following theorem which is analogous to Theorem \ref{thm:cluster tilting} and Theorem \ref{th:cyeq} \eqref{th:equivalence}.
\begin{theorem}
  \label{thm:cluster cat as stable cat for D type}
 Let $Q$ be a quiver of type $D_n$. Then
  \begin{enumerate}
     \item for a tagged triangulation $\sigma$ of the once-punctured polygon $P^*$, the cluster tilting object $e_F \Gamma_\sigma$ can be lifted to a tilting object in $\uCM^\ZZ(\Lambda)$;
     \item then there exists a triangle-equivalence $\cd^{\operatorname{b}}(KQ) \cong \uCM^{\ZZ}(\Lambda)$.
  \end{enumerate}
  \end{theorem}

\begin{proof}
(1) First, we have $e_F\Gamma_\sigma\cong \bigoplus_{a\in \sigma} M_a$. So we need to choose some degree shift of each $M_a$.

Suppose that all tagged arcs of $\sigma$ are incident to the puncture. Suppose without loss of generality that they are tagged plain. We can lift $\sigma$ to the set $\sigma'$ of indecomposable objects of $\CM^\ZZ \Lambda$-modules of the form $(i,*)$ for $1 \leq i \leq n$. Let us prove that the graded module $T'_{\sigma'} = \bigoplus_{A \in \sigma'} A$ is tilting (it is $T_\sigma$ if we forget the grading). Let us check that  $$\Hom_{\uCM^{\ZZ}(\Lambda)}((i,*), \Omega^k (j,*))=0$$ for any $i, j \in \cci{1}{n}$ and $k \neq 0$. Thanks to Theorem \ref{th:graded modules}, it is easy to compute projective covers of modules and we know that $\Omega^{k} (j,*) = (j + k(1-n),*)$ if $k$ is even and $\Omega^{k} (j,*) = (j + k(1-n),\bowtie)$ if $k$ is odd. Therefore, if $k$ is odd, $\Hom_{\uCM^{\ZZ}(\Lambda)}((i,*), \Omega^{k} (j,*))=0$.

Moreover, if $k \geq 2$, we get $j + k(1-n) \leq j + 2 - 2n \leq 0 < i$. So $$\Hom_{\uCM^{\ZZ}(\Lambda)}((i,*), \Omega^{k} (j,*))=0.$$ 

If $k \leq -2$ is even, we want to prove that $$\Hom_{\uCM^{\ZZ}(\Lambda)}((i,*), \Omega^{k} (j,*))= \Ext^1_{\CM^{\ZZ}(\Lambda)}((i,*), \Omega^{k+1} (j,*)) = 0.$$ We have $\Omega^{k+1} (j,*) = (j + (k+1)(1-n),\bowtie)$ and $j + (k+1)(1-n) \geq j + n - 1 \geq n \geq i$ and clearly $\Ext^1_{\CM^{\ZZ}(\Lambda)}((i,*), \Omega^{k+1} (j,*)) = 0$.

Let us now prove that $\thick(T'_{\sigma'})=\uCM^{\ZZ}(\Lambda)$. First of all, for any $i \in \ZZ$ such that $n \leq i < 2n-2$, considering the short exact sequence
$$0 \rightarrow (i-n+1,*) \rightarrow (i,i+1) \oplus (n-1, *) \rightarrow (i,2n-1) \rightarrow 0,$$
as $(i-n+1,*)$ and $(n-1,*)$ are in $\sigma'$ and $(i,i+1)$ is projective, we get that $(i,2n-1) \in \thick(T'_{\sigma'})$. Now, for any $i \in \ZZ$ such that $n < i < 2n-1$, using the short exact sequence
$$0 \rightarrow (n, 2n-1) \rightarrow (i, 2n-1) \oplus (n,*) \rightarrow (i, *) \rightarrow 0$$
we get that $(i,*) \in \thick(T'_{\sigma'})$. Thus, as $\Omega^{2k} (j,*) = (j + 2k(1-n),*)$, all the $(j, *)$ for $j \in \ZZ$ are in $\thick(T'_{\sigma'})$. Consider $i, j \in \ZZ$ such that $1 < j-i < n$. We then have a short exact sequence
$$0 \rightarrow (i+1-n,*) \rightarrow (i,i+1) \oplus (j-n, *) \rightarrow (i,j) \rightarrow 0$$
and, as $(i+1-n,*)$ and $(j-n, *)$ are in $\thick(T'_{\sigma'})$ and $(i,i+1)$ is projective, we get that $(i,j) \in \thick(T'_{\sigma'})$. Finally, as $\Omega(i,*) = (i-n+1, \bowtie)$, all the $(i,\bowtie)$ are in $\thick(T'_{\sigma'})$. We finished to prove in this case that $T'_{\sigma'}$ is tilting.

Suppose now that there is at least a tagged arc of $\sigma$ which is not incident to the puncture. Thus, there exists a vertex $i_0$ of $P$ such that $i_0$ does not have any incident internal edge in $\sigma$. Therefore, we can lift the tagged arcs of $\sigma$ to a set $\sigma'$ of indecomposable objects of $\CM^\ZZ \Lambda$-modules such that for any $A \in \sigma'$, $i_0 < A_1 < i_0+n$ and $i_0+1 < A_2 < i_0 + 2n$. Let us prove that the graded module $T'_{\sigma'} = \bigoplus_{A \in \sigma'} A$ is tilting (it is $T_\sigma$ if we forget the grading). Let us check that  $$\Hom_{\uCM^{\ZZ}(\Lambda)}(A, \Omega^k B)=0$$ for any $A, B\in \sigma'$ and $k \neq 0$. Let $B' = \Omega^{k} B$. Thanks to Theorem \ref{th:graded modules}, $B'_1= B_1 + k(1-n)$ and $B'_2= B_2 + k(1-n)$.

Therefore, if $k > 0$, we get $B_1' \leq B_1 + 1 - n \leq i_0 < A_1$. So $\Hom_{\uCM^{\ZZ}(\Lambda)}(A, B')=0$. 

If $k < -1$, we want to prove that $\Hom_{\uCM^{\ZZ}(\Lambda)}(A, B')= \Ext^1_{\CM^{\ZZ}(\Lambda)}(A, \Omega B') = 0$. If we denote $B'' = \Omega B'$, we have $B''_1 = B'_1+1-n \geq B_1 + n - 1 \geq i_0 + n > A_1$. Then as the morphisms are directed positively, $\Ext^1_{\CM^{\ZZ}(\Lambda)}(A, \Omega B') = 0$. 

For $k = -1$, by Theorem \ref{thm:cluster tilting}, we get $$\Hom_{\uCM^{\ZZ}(\Lambda)}(T'_{\sigma'}, \Omega^{-1} T'_{\sigma'}) {\subset}\Hom_{\uCM \Lambda}(T_{\sigma},  \Omega^{-1} T_{\sigma})= \Ext^1_{\CM \Lambda}(T_{\sigma},  T_{\sigma})= 0.$$

Let us now prove that $\thick(T'_{\sigma'})=\uCM^{\ZZ}(\Lambda)$. Consider an indecomposable object $A \in \CM^\ZZ \Lambda$ with $A_1 = i_0+n$. Let $A' \in \CM \Lambda$ be its image through the forgetful functor. It is a classical lemma about cluster tilting objects that there exists a short exact sequence $$0 \rightarrow T'_1 \rightarrow T'_0 \rightarrow A' \rightarrow 0$$ of Cohen-Macaulay $\Lambda$-modules such that $T'_0, T'_1 \in \add(T_\sigma)$. 

Let $X'$ be an indecomposable summand of $T'_1$. For any lift $X$ of $X'$ such that $\Ext^1_{\CM^\ZZ \Lambda}(A, X) \neq 0$, we have $i_0 < X_1 < i_0 + n$ by Theorem \ref{th:graded modules} (2), so such a lift $X$ is unique and has to be in $\sigma'$. Moreover, in this case, $\Ext^1_{\CM^\ZZ \Lambda}(A, X) = \Ext^1_{\CM \Lambda}(A, X')$. Therefore, the unique lift $T_1$ of $T'_1$ which is in $\add(T'_{\sigma'})$ satisfies $\Ext^1_{\CM^\ZZ \Lambda}(A, T_1) = \Ext^1_{\CM \Lambda}(A, T_1')$ so we can lift the short exact sequence $$0 \rightarrow T'_1 \rightarrow T'_0 \rightarrow A' \rightarrow 0$$ to a short exact sequence $$0 \rightarrow T_1 \rightarrow T_0 \rightarrow A \rightarrow 0$$ of graded Cohen-Macaulay $\Lambda$-modules. As any indecomposable summand $X$ of $T_0$ is between $T_1$ and $A$ in the Auslander-Reiten quiver, we get $i_0 <  X_1 \leq A_1 = i_0+n$ so $X \in \sigma'$. Finally $T_0$ and $T_1$ are in $\add(T'_{\sigma'})$ so $A \in \thick(T'_{\sigma'})$.

For any $i \in \ZZ$ such that $i_0+n < i <i_0 + 2n - 1$, there is a short exact sequence
$$0 \rightarrow (i_0 + n, i+1) \rightarrow (i, i+1) \oplus (i_0 +n, *) \rightarrow (i, *) \rightarrow 0$$
so, as $(i_0 + n, i+1)$ and $(i_0 + n, *)$ are in $\thick(T'_{\sigma'})$ and $(i,i+1)$ is projective, $(i,*)$ is in $\thick(T'_{\sigma'})$. As we already got the result for $(i_0+n,*)$ and $(i_0+n,\bowtie)$ and $\Omega^{-1}((i_0+n,\bowtie)) = (i_0+2n-1, *)$, all the $(i,*)$ for $i_0 +n \leq i \leq i_0 + 2n - 1$ are in $\thick(T'_{\sigma'})$. Up to a shift by $1-i_0-n$, we already saw that these $(i,*)$ generate $\uCM^\ZZ(\Lambda)$. Therefore, $T'_{\sigma'}$ is a tilting object in $\uCM^\ZZ(\Lambda)$.

(2) Take the triangulation $\sigma$ whose set of tagged arcs is $$\{(P_1,P_3),(P_1,P_4),\cdots,(P_1,P_n),(P_1,*),(P_1,\bowtie)\}.$$
The full subquiver $Q$ of $Q_\sigma$ with the set $Q_{\sigma,0} \smallsetminus F$ of vertices is of type $D_n$. 
Thus, we have $$\Gamma_\sigma^{\operatorname{op}} / (e_F) \cong (KQ)^{\operatorname{op}}.$$
Let $\sigma'$ constructed from $\sigma$ as before. Namely, $i_0 = 2$ and 
 $$\sigma' = \{(n+1,n+3), (n+1,n+4), \dots, (n+1, 2n), (n+1,*), (n+1, \bowtie) \}.$$

 For any $A, B \in \sigma'$ and $k \in \ZZ$, we have $B(k)_1 = n+1 + kn$
so $$\Hom_{\uCM^\ZZ(\Lambda)}(A,B(k)) = 0$$ for $k \neq 0$. Indeed, if $k < 0$ this is immediate as morphisms go increasingly and if $k > 0$,
$\Hom_{\uCM^{\ZZ}(\Lambda)}(A, B(k))= \Ext^1_{\CM^{\ZZ}(\Lambda)}(A, \Omega B(k))$.
Moreover, $\Omega B(k)_1 = 2 + kn \geq n+2$ and for the same reason than before $\Ext^1_{\CM^{\ZZ}(\Lambda)}(A, \Omega B(k)) = 0$.

Thus, by Theorem \ref{thm:cluster tilting},  we have the following isomorphism 
\begin{align*}
\End_{\uCM^\ZZ(\Lambda)}(e_F \Gamma_\sigma)& \cong \End_{\uCM \Lambda}(e_F \Gamma_\sigma ) \\
&\cong \Gamma_\sigma^{\operatorname{op}} / (e_F).
\end{align*}

 Because $\uCM^\ZZ(\Lambda)$ is an algebraic triangulated Krull-Schmidt category, and $e_F \Gamma_\sigma$ is a tilting object in $\uCM^\ZZ(\Lambda)$, by Theorem \ref{BK}, there exists a triangle-equivalence $$\uCM^\ZZ(\Lambda) \cong \ck^{\operatorname{b}}\left(\proj\End^{\operatorname{op}}_{\uCM^\ZZ(\Lambda)}(e_F \Gamma_\sigma)\right).$$
Since $\gldim K{Q} < \infty$, we have a triangle-equivalence $$\ck^{\operatorname{b}}\left(\proj\End^{\operatorname{op}}_{\uCM^\ZZ(\Lambda)}(e_F \Gamma_\sigma)\right) \cong\ck^{\operatorname{b}}(\proj KQ) \cong \cd^{\operatorname{b}}(KQ).$$
Therefore there is a triangle-equivalence $\uCM^\ZZ(\Lambda) \cong \cd^{\operatorname{b}}(KQ)$. 
\end{proof}

\bibliographystyle{plain}
\bibliography{biblio}

\begin{thebibliography}{10}

\bibitem{AIRCC}
Claire Amiot, Osamu Iyama, and Idun Reiten.
\newblock Stable categories of {C}ohen-{M}acaulay modules and cluster
  categories.
\newblock {\em Amer. J. Math.}, 137(3):813--857, 2015.

\bibitem{Araya}
Tokuji Araya.
\newblock Exceptional sequences over graded {C}ohen-{M}acaulay rings.
\newblock {\em Math. J. Okayama Univ.}, 41:81--102 (2001), 1999.

\bibitem{FMO}
M.~Auslander.
\newblock Functors and morphisms determined by objects.
\newblock In {\em Representation theory of algebras ({P}roc. {C}onf., {T}emple
  {U}niv., {P}hiladelphia, {P}a., 1976)}, pages 1--244. Lecture Notes in Pure
  Appl. Math., Vol. 37. Dekker, New York, 1978.

\bibitem{BK1991}
A.~I. Bondal and M.~M. Kapranov.
\newblock Framed triangulated categories.
\newblock {\em Mat. Sb.}, 181(5):669--683, 1990.

\bibitem{BIRS}
A.~B. Buan, O.~Iyama, I.~Reiten, and D.~Smith.
\newblock Mutation of cluster-tilting objects and potentials.
\newblock {\em Amer. J. Math.}, 133(4):835--887, 2011.

\bibitem{BMRRT}
A.~B. Buan, R.~Marsh, M.~Reineke, I.~Reiten, and G.~Todorov.
\newblock Tilting theory and cluster combinatorics.
\newblock {\em Adv. Math.}, 204(2):572--618, 2006.

\bibitem{LF11}
Giovanni Cerulli~Irelli and Daniel Labardini-Fragoso.
\newblock Quivers with potentials associated to triangulated surfaces, {P}art
  {III}: tagged triangulations and cluster monomials.
\newblock {\em Compos. Math.}, 148(6):1833--1866, 2012.

\bibitem{RT}
C.~W. Curtis and I.~Reiner.
\newblock {\em Methods of representation theory. {V}ol. {I}}.
\newblock John Wiley \& Sons Inc., New York, 1981.
\newblock With applications to finite groups and orders, Pure and Applied
  Mathematics, A Wiley-Interscience Publication.

\bibitem{Thanvanden}
L.~{de Thanhoffer de V{\"o}lcsey} and M.~{Van den Bergh}.
\newblock {Explicit models for some stable categories of maximal Cohen-Macaulay
  modules}, June 2010.
\newblock arXiv: 1006.2021.

\bibitem{DeLu}
L.~{Demonet} and X.~{Luo}.
\newblock {Ice quivers with potentials associated with triangulations and
  Cohen-Macaulay modules over orders}, July 2013.
\newblock arXiv: 1307.0676, accepted for publication to Trans. Amer. Math. Soc.

\bibitem{DWZ1}
H.~Derksen, J.~Weyman, and A.~Zelevinsky.
\newblock Quivers with potentials and their representations. {I}. {M}utations.
\newblock {\em Selecta Math. (N.S.)}, 14(1):59--119, 2008.

\bibitem{FST08}
S.~Fomin, M.~Shapiro, and D.~Thurston.
\newblock Cluster algebras and triangulated surfaces. {I}. {C}luster complexes.
\newblock {\em Acta Math.}, 201(1):83--146, 2008.

\bibitem{GLS11}
C.~Gei{\ss}, B.~Leclerc, and J.~Schr{\"o}er.
\newblock Kac-{M}oody groups and cluster algebras.
\newblock {\em Adv. Math.}, 228(1):329--433, 2011.

\bibitem{Happel}
D.~Happel.
\newblock {\em Triangulated categories in the representation theory of
  finite-dimensional algebras}, volume 119 of {\em London Mathematical Society
  Lecture Note Series}.
\newblock Cambridge University Press, Cambridge, 1988.

\bibitem{Happel2}
D.~Happel.
\newblock Auslander-{R}eiten triangles in derived categories of
  finite-dimensional algebras.
\newblock {\em Proc. Amer. Math. Soc.}, 112(3):641--648, 1991.

\bibitem{HiNi92}
H.~Hijikata and K.~Nishida.
\newblock Classification of {B}ass orders.
\newblock {\em J. Reine Angew. Math.}, 431:191--220, 1992.

\bibitem{HiNi97}
Hiroaki Hijikata and Kenji Nishida.
\newblock Primary orders of finite representation type.
\newblock {\em J. Algebra}, 192(2):592--640, 1997.

\bibitem{Iya-order}
O.~Iyama.
\newblock Representation theory of orders.
\newblock In {\em Algebra---representation theory ({C}onstanta, 2000)},
  volume~28 of {\em NATO Sci. Ser. II Math. Phys. Chem.}, pages 63--96. Kluwer
  Acad. Publ., Dordrecht, 2001.

\bibitem{IyamaLerner}
O.~{Iyama} and B.~{Lerner}.
\newblock {Tilting bundles on orders on $P^d$}, June 2013.
\newblock arXiv: 1306.5867.

\bibitem{IyTa13}
O.~Iyama and R.~Takahashi.
\newblock Tilting and cluster tilting for quotient singularities.
\newblock {\em Math. Ann.}, 356(3):1065--1105, 2013.

\bibitem{Matrixrepresentation}
Hiroshige Kajiura, Kyoji Saito, and Atsushi Takahashi.
\newblock Matrix factorization and representations of quivers. {II}. {T}ype
  {$ADE$} case.
\newblock {\em Adv. Math.}, 211(1):327--362, 2007.

\bibitem{Matrixweight}
Hiroshige Kajiura, Kyoji Saito, and Atsushi Takahashi.
\newblock Triangulated categories of matrix factorizations for regular systems
  of weights with {$\epsilon=-1$}.
\newblock {\em Adv. Math.}, 220(5):1602--1654, 2009.

\bibitem{derdgcat}
B.~Keller.
\newblock Deriving {DG} categories.
\newblock {\em Ann. Sci. \'Ecole Norm. Sup. (4)}, 27(1):63--102, 1994.

\bibitem{KR}
B.~Keller and I.~Reiten.
\newblock Acyclic {C}alabi-{Y}au categories.
\newblock {\em Compos. Math.}, 144(5):1332--1348, 2008.
\newblock With an appendix by Michel Van den Bergh.

\bibitem{LF09}
D.~Labardini-Fragoso.
\newblock Quivers with potentials associated to triangulated surfaces.
\newblock {\em Proc. Lond. Math. Soc. (3)}, 98(3):797--839, 2009.

\bibitem{LF12}
D.~{Labardini-Fragoso}.
\newblock {Quivers with potentials associated to triangulated surfaces, part
  {IV:} Removing boundary assumptions}, June 2012.
\newblock arXiv: 1206.1798.

\bibitem{Lad}
Sefi {Ladkani}.
\newblock {On Jacobian algebras from closed surfaces}, July 2012.
\newblock arXiv: 1207.3778.

\bibitem{almostsplitseqorder}
K.~W. Roggenkamp.
\newblock The construction of almost split sequences for integral group rings
  and orders.
\newblock {\em Comm. Algebra}, 5(13):1363--1373, 1977.

\bibitem{MR0412223}
K.~W. Roggenkamp and J.~W. Schmidt.
\newblock Almost split sequences for integral group rings and orders.
\newblock {\em Comm. Algebra}, 4(10):893--917, 1976.

\bibitem{POS}
D.~Simson.
\newblock {\em Linear representations of partially ordered sets and vector
  space categories}, volume~4 of {\em Algebra, Logic and Applications}.
\newblock Gordon and Breach Science Publishers, Montreux, 1992.

\bibitem{CM}
Y.~Yoshino.
\newblock {\em Cohen-{M}acaulay modules over {C}ohen-{M}acaulay rings}, volume
  146 of {\em London Mathematical Society Lecture Note Series}.
\newblock Cambridge University Press, Cambridge, 1990.

\end{thebibliography}
\end{document}